\documentclass[10pt]{amsart}
\usepackage{amssymb}
\usepackage{graphicx}
\usepackage{color}

\def\a{{\bf a}} \def\b{{\bf b}} \def\c{{\bf c}} 
\def\d{{\bf d}} \def\e{{\bf e}} \def\f{{\bf f}} 
  
 \def\v{{\bf v}} \def\w{{\bf w}} 
\def\x{{\bf x}} \def\y{{\bf y}} \def\z{{\bf z}}
 
\def\sa{\hbox{\bf{\tiny a}}} \def\sb{\hbox{\bf{\tiny b}}} \def\sc{\hbox{\bf{\tiny c}}} 
\def\sd{\hbox{\bf{\tiny d}}} \def\se{\hbox{\bf{\tiny e}}}

\def\sx{\hbox{\bf{\tiny x}}} \def\sy{\hbox{\bf{\tiny y}}} \def\sz{\hbox{\bf{\tiny z}}}
 
  \def\nn{{\mathbb N}} 
  \def\rr{{\mathbb R}} 
\def\zz{{\mathbb Z}}

\def\AA{{\mathcal A}}  \def\CC{{\mathcal C}} 
 \def\EE{{\mathcal E}}  
\def\GG{{\mathcal G}}   
   
\def\MM{{\mathcal M}}   
\def\PP{{\mathcal P}} \def\QQ{{\mathcal Q}}

\def\barp{\bar{p}}

\def\bAA{\bar{\mathcal A}}
\def\bEE{\bar{\mathcal E}}
\def\bMM{\bar{\mathcal M}}
\def\path{{\rm Path}}
\def\bpath{\overline{\path}}
\def\tpath{\widetilde{\path}}
\def\tdelta{\tilde{\delta}}
  
\def\bw{\bar{\w}}
\def\bv{\bar{\v}}

 \def\bN{\bar{N}}  
  
   \def\bX{\bar{X}}

\def\bphi{\bar{\phi}}  
\def\bmu{\bar{\mu}}

\def\tMM{\widetilde{\MM}} 
\def\tAA{\tilde{\AA}} 
\def\tCC{\tilde{\CC}}
\def\tGG{\tilde{\GG}}
\def\tw{\tilde{\w}}   
\def\tv{\tilde{\v}}

\def\tC{\tilde{C}}
\def\tphi{\tilde{\phi}}
\def\tpsi{\tilde{\psi}}
\def\btp{\beta\tilde{\phi}}

\def\di{\displaystyle} 
\def\bs{\bigskip\noindent} 
\def\ms{\medskip\noindent} 

\def\lfphi{\overleftarrow{\phi}} 
\def\rtphi{\overrightarrow{\phi}} 
\def\cirb{\stackrel{\curvearrowleft}{\leadsto}} 
\def\lds{\leadsto}
\def\cira{\circlearrowleft}  

\def\var{\hbox{\rm var}} 
\def\per{\hbox{\rm Per}} 
\def\spec{\hbox{\rm spec}}
\def\sper{\hbox{\tiny \rm Per}}

\def\top{\hbox{\tiny top}} 
\def\bp{\beta\phi}
\def\bpp{\beta\phi+\psi}
\def\trace{\hbox{\rm tr}} 
\def\clos{\hbox{\rm clos}} 
\def\1{\hbox{\rm Id}}

\newtheorem{theorem}{\em Theorem} 
\newtheorem{definition}{\em Definition} 
\newtheorem{remark}{\em Remark}
\newtheorem{question}{\em Question}  
\newtheorem*{algo}{\em Algorithm for the $\alpha_J$'s} 
\newtheorem*{facts}{\em Facts}  
\newtheorem*{notations}{\em Notations} 
\newtheorem*{renorm}{\em Renormalization Lemma}  
\newtheorem{proposition}{\em Proposition}
 
\newtheorem{Lemma}{\em Lemma}
\newtheorem{corollary}{\em Corollary}

\def\endproof{\hfill{$\Box$}}
 
\begin{document}

\title[Zero-temperature limit of one-dimensional Gibbs states]
{Zero-temperature limit of one-dimensional Gibbs states via renormalization:\\
the case of locally constant potentials}

\author{J.--R. Chazottes, J.--M. Gambaudo \& E. Ugalde}

\date{21st of april 2010} 


\thanks{This work is part of the project {\it CrystalDyn} funded 
by the Agence Nationale de la Recherche (ANR)}

\begin{abstract}
Let  $A$ be a finite set and $\phi:A^\zz\to\rr$ be a locally constant potential.
For each $\beta>0$ (``inverse temperature''), there is a unique
Gibbs measure $\mu_{\beta\phi}$. We prove that, as $\beta\to+\infty$, the
family $(\mu_{\beta\phi})_{\beta>0}$
converges (in weak-$^*$ topology) to a measure we characterize. 
It is concentrated on a certain subshift of finite type 
which is a finite union of transitive subshifts of finite type.
The two main tools are an approximation by periodic orbits and 
the Perron-Frobenius Theorem for matrices \`a la Birkhoff. 
The crucial idea we bring is a ``renormalization'' procedure
which explains convergence and provides a recursive
algorithm to compute the weights of the ergodic decomposition of the limit.
\end{abstract}

\maketitle

\tableofcontents

\bs \section{Introduction}\

A fundamental problem in equilibrium statistical mechanics
is the description of the set of Gibbs states for a given
potential as temperature changes and, more specifically,
as it goes to zero. 
Already for {\em classical lattice systems}, this is a formidable task  \cite{Georgii88, vanEnter&al93}.
There is no {\em a priori} reason why the zero-temperature limit should exist at all, even for
finite-range potentials. Indeed, it was proven recently that this may not be the case \cite{vER,jrmike}. 
We shall prove, however, that for the one-dimensional lattice $\zz$ and finite-range potentials, it does
always exist and can be identified.

The problem we study can be formulated in the following way.
We consider the space $A^\zz$ of two-sided sequences or configurations,
where $A$ is a finite set, and let $\phi:A^\zz\to\rr$ be a locally constant potential (see below for the
precise definition).
Such a potential admits a unique Gibbs state which is none other than
a $r$-step Markov measure with state space $A$.
(When $r=1$, this is a usual Markov measure.)
Ground states are probability
measures achieving the maximum of the mapping
$\nu\mapsto \int \phi d\nu$, where $\nu$ ranges over the set of shift-invariant probability measures
on $A^\zz$. It is not difficult to prove that ground states of locally constant potentials
are necessarily supported on (nonwandering) {\em subshifts of finite type} (henceforth SFT's). In particular, this implies that the set of ground states for a given locally constant potential is a finite-dimensional
simplex whose extreme points are exactly the ergodic ground states of $\phi$.
Such SFT's have no reason to be transitive but are in general
made of a finite number of transitive SFT's, each of them being
possibly a (finite) union of topologically mixing SFT's cyclically permutated by the shift.

The question we are interested in is: 
\begin{question}\label{question}
For a locally constant potential $\phi$, does the limit (in weak-$^*$ topology)
of $(\mu_{\beta\phi})_{\beta>0}$, as $\beta\to+\infty$, exist ? Can it be precisely described ? 
\end{question}

It follows easily from the variational principle that any weak-$^*$ accumulation
point of $(\mu_{\beta\phi})_{\beta>0}$ is a ground state of $\phi$.
In the particular case where the ground states of $\phi$ are supported on a transitive
SFT, in particular a single periodic configuration, then convergence is enforced
by general arguments we recall below. Otherwise, there is no reason why
convergence should arise. And assuming it does, it is not clear how the limiting
measure spreads among the ergodic ground states.
A basic example for which the zero-temperature limit can be easily computed is
the one-dimensional Ising model (see \cite[Section 3.2]{Georgii88} for full details).
In that model, $A=\{-,+\}$ and one can tune two parameters. Depending on their values, the zero-temperature
limit can be, for instance, $\delta_+$, the Dirac measure on the `all-$+$' configuration; or
$(\delta_{-+}+\delta_{+-})/2$, where $\delta_{ab}$ stands for the Dirac measure
on the periodic configuration ($\cdots ababab\cdots$); but it can also be the unique
measure of maximal entropy of the SFT defined as the set of all configurations
in $\{-,+\}^\zz$ with no two consecutive $-$'s. In all cases, one gets either
the centroid of the ergodic ground states (when there are finitely many periodic ground configurations)
or the measure of maximal entropy of a certain transitive SFT. Another example with $A=\{1,2,3\}$ 
appears in \cite[Section 9]{Jenkinson01}. For an appropriate locally constant potential
the ergodic ground states are the Dirac measures $\delta_1,\delta_2,\delta_3$, where $\delta_i$
stands for the all-$i$ configuration. A tedious calculation shows that the limiting Gibbs measure
is $(\delta_1 + \delta_2)/2$. Hence, we do not get the centroid of the ergodic ground
states but the one of two of them. Below, we shall see examples where the barycentric coordinates
are irrational numbers (hence the limiting Gibbs measure is not a centroid of any of the ground states).
We notice that it was proved in \cite{Nekhoroshev04} that, {\em generically}, the zero-temperature limit of
Gibbs states of locally constant potentials exists and is of the form
$(\delta_{\a}+\cdots+\delta_{\sigma^{n-1}\a})/n$ where $\sigma:A^\zz\to A^\zz$
is the shift and $\a$ is a periodic configuration with period $n$ ($\sigma^n \a=\a$).

Question \ref{question} was previously tackled in \cite{Bremont03} and \cite{Leplaideur05}.
In \cite{Bremont03}, convergence is proved as a consequence of a general statement of Analytic Geometry
(especially the theory of subanalytic maps). The limiting measure is not identified.
In \cite{Leplaideur05} a more explicit approach is used to prove convergence, and the
limit is partially identified.
The weights of the barycentric decomposition of the limiting measure are not explicitly identified and
the proof of convergence is somewhat indirect. 
In the present work, we offer an approach differing from \cite{Bremont03, Leplaideur05}.
It is based on two main tools and a new idea. These tools are
approximation by periodic orbits and the contraction-mapping approach
to the Perron-Frobenius theorem for matrices \cite{Birkhoff57, Seneta81}.
The new idea is a {\em renormalization procedure} which has to be
iterated only a finite number of times (this is due to the fact that we
consider locally constant potentials). Not only does this algorithm
explain convergence but it also allows one to determine recursively the
coefficients of the ergodic decomposition of the zero-temperature
limit. Finally, we also obtain that $\mu_{\beta\phi}$  is exponentially
close to its limit for large enough $\beta$.

Colloquially, the renormalization works as follows.
As mentioned above, the ground states of a given locally constant potential $\phi$
are supported on a SFT $\bX(\phi)$ which decomposes as a union of some transitive SFT's $\bX_J(\phi)$.
In fact, only what we call hereinafter the ``heavy components'' do support the accumulation points
of $(\mu_{\beta\phi})_{\beta>0}$. 
The renormalization consists in defining an SFT $\bX'(\phi)$ whose alphabet is obtained by
labeling the heavy components by $\{1,\ldots,N_\phi\}$ and keeping only
certain arrows between them. Correspondingly, one has to ``renormalize'' the potential
$\phi$. Then we prove that the original problem is reduced to a new one where
$\beta\phi$ has to be replaced by $\beta\phi' + \psi'$, where $\phi',\psi'$ are locally
constant potentials living on the SFT $\bX'(\phi)$. The additional potential $\psi'$, independent
of $\beta$, is a ``compensation factor''. This renormalization procedure is of course
approximative and leads to an error term which depends on $\beta$. So to speak,
when $\beta$ is ``large enough'' ({\em i.e.} temperature is ``small enough''),
we have a simpler system which ``mimics'' the real one. After a {\em finite} number
of iterations of this procedure, one ends up with a transitive SFT. Then we obtain
the weights of those heavy components which carry the limiting measure, some
of them having been eliminated. We provide below an example where one has
to renormalize twice to compute the limiting measure.

We emphasize two crucial ingredients in the renormalization, namely
the use of cohomology of potentials and the fact that we do not loose the
Markovian character of the measure when we `collapse' each heavy
component. 

One may ask what happens beyond locally constant potentials. It has been proved very recently
in \cite{jrmike} that $(\mu_{\beta\phi})_{\beta>0}$ may not converge even for Lipschitz potentials
$\phi$. The idea is to construct a minimal subshift $Y$ supporting only two ergodic measures $\mu_1$
and $\mu_2$, and consider the Lipschitz potential $\phi(\a)=-d(\a,Y)$, where $d$ is the usual distance on $A^\zz$.
Then one shows that $\mu_{\beta\phi}$ ``oscillates'' between $\mu_1$ and $\mu_2$ which are
the two ergodic ground states by contruction.
In principle there is room between locally constant potentials and Lipschitz ones for having
convergence, but the class of functions in ``between'' has not a clear relevance to us.

Let us end this introduction with a few remarks on statistical mechanics.
There, ground states are rather those $\nu$'s which minimize $\int \phi d\nu$.
This is due to the fact that different sign conventions are used in the thermodynamic formalism \cite{Ruelle}.
Another point is that a `potential' in statistical mechanics is a family of shift-invariant functions $\Phi=(\Phi_\Lambda)$ indexed by finite subsets $\Lambda$ of $\zz$. What dynamicists call a potential, as we do here, is the function
$\phi:A^\zz\to\rr$ obtained from $\Phi$ by $\phi=\sum_{\Lambda\ni 0} \Phi_\Lambda/|\Lambda|$. 
One can ``average'' $\Phi$ in other ways leading to the same result upon integration
by shift-invariant measures \cite{Ruelle}. The standard setting is that of ``absolutely summable''
potentials $\Phi$.  Let us notice that locally constant $\phi$'s
are exactly the so-called finite-range $\Phi$'s. 
The function $\phi$ naturally appears in statistical mechanics, when one deals with equilibrium
states, but bears no special name. Its physical interpretation is clear: it is the ``energy''
per site of the system.

\bigskip

\noindent {\bf Scope of the article}.
In Section \ref{setup}, we give some definitions and recall some general results about
zero-temperature limits and ground states. In Section \ref{reduc} we explain how
Question \ref{question} can be recast into a simpler one, without loss of generality.
In Section \ref{main} we state the main result (Theorem \ref{main-theorem} and
the corresponding algorithm). We then prove the main result in Section \ref{mainproof}
assuming the key-lemma (``Renormalization Lemma'') whose proof is deferred to Section \ref{thebiglemma}.
The proof of this lemma relies on a certain number of technical lemmas proved in
Section \ref{lessouslemmes}.
In turn, we postponed to Appendices \ref{proof-of-periodic}
and \ref{auxiliary-inequalities} the proof of a number of statements used in the proofs of those lemmas.
We have made a section with three examples (Section \ref{examples}) to illustrate our result, in particular
the renormalization procedure. Two of these examples seem to be new and we would not have
found them without the renormalization idea. 

\bs \section{Settings and Generalities}\label{setup} \

\ms Let  $A$ be a finite set with at least two elements and $A^{\zz}$
be the set of two-sided infinite sequences of symbols drawn from $A$.
Elements  of $A^{\zz}$ will be denoted by $\a$, $\b$, $\c$, etc.
The shift map $\sigma:A^{\zz}\to A^{\zz}$ is defined by $(\sigma\a)_i=\a_i$ for all $i\in\zz$.
Endowed with the product topology, $A^{\zz}$ is a compact metrizable space.
Given $\b\in A^{\zz}$ and $p, q \in\zz$ such that  $p \leq q$, we denote by $[\b_q^p]$ the 
cylinder set $\{\a \in A^{\zz}:\ \a_p = \b_p,\ldots, \a_q = \b_q\}$. 

\ms \subsection{Equilibrium States and Pressure}\

\ms We recall a few facts on equilibrium states. We refer the reader to, {\em e.g.}, \cite{Bowen} 
and~\cite{Ruelle} for details.

\ms Let $\phi: A^{\zz} \to \rr$ be a continuous function which we call a potential.
This means that $\var_{n}\phi\to 0$ as $n\to +\infty$, where 
$\var_{n}\phi := \sup\{|\phi(\a)-\phi(\b)| : \a_j = \b_j , |j|\leq n\}$ 
is the modulus of continuity of $\phi$.
The equilibrium states $\phi$ are those shift-invariant probability measures
realizing the supremum of $h(\nu) + \int\phi d\nu$  among all shift--invariant probability measures $\nu$
(where $h(\nu)$ is the measure--theoretic entropy of $\nu$).
The supremum equals $P(\phi)$, the topological pressure of $\phi$.
Endowed with the weak-$^*$ topology, the set of shift-invariant probability measures is a compact convex set,
in fact a Choquet simplex~\cite{Ruelle}. The set of equilibrium states of $\phi$ is a face of that simplex.

\ms Let $Y\subset A^{\zz}$ be a subshift of $A^\zz$
({\em i.e.}, a closed shift-invariant subset of $A^\zz$) and $\psi: A^{\zz}\to \rr$ be continuous.
If we restrict $\psi$ to $Y$, we can define its equilibrium states on $Y$. We denote the 
corresponding topological pressure by $P(\psi|Y)$, which we call the $\psi$-pressure on $Y$.
When $\psi\equiv 0$, $P(\psi|Y)=h_{\top}(Y)$, {\em i.e.} the $0$-pressure on $Y$ is
the topological entropy of $Y$.

\ms \subsection{Maximizing Measures and the Maximizing Subshift $\bX$}\

\ms For $\phi: A^{\zz}\to\rr$ continuous, let
\[
\bphi:=\sup\left\{\int\phi d\nu:\ \nu\;\text{{\small shift-invariant probability measure}}\right\}.
\]
By compactness, there always exists an invariant measure 
which realizes this supremum; we call it a {\em maximizing measure for $\phi$}. The
shift-invariant, compact  set containing the support of all $\phi$-maximizing measures is denoted by
$\bX = \bX(\phi)\subset A^{\zz}$ and we call it the {\em $\phi$-maximizing subshift}. 
The set of $\phi$-maximizing measures is a face of the Choquet simplex of shift-invariant probability
measures, hence it is also a Choquet simplex. Its extreme points are precisely the ergodic $\phi$-maximizing measures. We use the usual notations $S_p\phi:=\sum_{i=0}^{p-1}\phi\circ\sigma^i$ and 
$\per_p(A^{\zz}):=\{\a\in A^{\zz}: \sigma^p\a=\a\}$ ($p\in \nn$).
First,
\begin{equation}\label{phibarutile}
\bphi=\sup_{p\in\nn}\ \max_{\sa\in\sper_p(A^{\zz})}\frac{S_p\phi(\a)}{p}\cdot
\end{equation}
Second,
$$
\bX=\bX(\phi):=\overline{\bigcup_{p\in\nn}\left\{\a\in \per_p(X):\ S_p\phi(\a)=p\bphi\right\}}.
$$

\ms \subsection{Basic facts about Zero--temperature Accumulation Points}\

\ms Given a continuous $\phi$, we consider the one-parameter family of equilibrium
states $\{\mu_{\bp}:\, \beta > 0\}$, where $\beta$ is interpreted in statistical physics as the inverse
temperature. For each $\beta$, the potential $\bp$ admits at least one equilibrium state.

\ms We collect general basic facts (that we will not need) relating zero-temperature accumulation points
of $(\mu_{\bp})_{\beta>0}$ to maximizing measures, when $\phi$ is continuous. 

\ms \begin{facts}
Let $\phi:A^{\zz}\to \rr$ be a continuous potential. Then
\begin{enumerate}
\item $(\mu_{\bp})_{\beta>0}$ has a weak${}^*$ accumulation point, as $\beta\to +\infty$.
\item Every accumulation point of $(\mu_{\bp})_{\beta>0}$ is a $\phi$-maximizing measure.
Its support is contained in $\bX=\bX(\phi)$.
\item Every accumulation point of $(\mu_{\bp})_{\beta>0}$ is of maximal entropy among  
$\phi$-maximizing measures.
\end{enumerate}
\end{facts}

\ms Let us make a few comments on these facts. 
One can find the proofs of the above facts in~\cite{ConzeGivarch95}.
Statement (1) results by compactness of the set of shift--invariant 
probability measures. Statement (2) is a straightfoward consequence of the variational principle. 
The last statement is also a consequence of the variational principle and of the convexity of 
$\beta\mapsto P(\beta\phi)$. In the formalism of statistical mechanics of lattice systems, these statements were 
known long before and can be found in~\cite[Appendix B.2]{vanEnter&al93}.
The above proposition allows to conclude the existence of the zero-temperature limit 
in some special cases: when $\bX$ supports a unique shift-invariant measure or when $\bX$ has
a unique measure of maximal entropy ({\em e.g.}, when it is a transitive SFT).

A slight generalization is to consider the case $\bpp$ where $\psi:A^{\zz}\to \rr$
is another continuous potential. The above facts are valid with this new family of potentials
if one replaces `maximal entropy'  by `maximal $\psi$-pressure' in statement (3). Notice
that $\bX$ does not depend on $\psi$ but only on $\phi$.

\bs \section{Preparatory results and reduction of the problem}\  \label{reduc}

\ms \subsection{Locally Constant Potentials and Markov Measures}\

\ms We say that $\phi$ is locally constant if there exists a strictly positive integer $r$ such that 
\[
\phi(\a) = \phi(\b), \, \forall\, \a, \b \in A^{\zz}\, \text{ such that }\, \a\in [\b_0^r].
\]
We say that $\phi$ is a $(r+1)$--symbol potential. 
A $(r+1)$--symbol potential can be identified with a function from
$A^{r+1}$ to $\rr$ that we can still denote by $\phi$ by a slight abuse of notation.
In that case, for each $\beta > 0$, 
the potential $\bp$ admits a unique equilibrium measure, which is also a ($r$--step)
Markov measure. (Notice that the case $r = 0$ corresponds to the case of product measures for
which the zero-temperature limit problem is trivial.)

Without loss of generality, we can reduce our problem to the case of (1-step) Markov measures, {\em i.e.}, to the case of 2--symbol potentials. We make this precise in the next section.

\ms \subsection{Recodification for Locally Constant Potentials, Maximizing SFT
and Heavy Components}\ 

\ms Let $\Phi,\Psi:A^{r+1}\to\rr$ be $(r+1)$-symbol potentials and let $\AA$ be the alphabet of words
of length $r+1$ in the alphabet $A$. Then $A^\zz$ can be recoded as a topological Markov chain (hereafter
TMC) $X\subset \AA^{\zz}$ (that is, a SFT which can be defined by words of length two). 
 
\ms We identify $X$ with the set of all bi--infinite paths on a directed graph 
$\GG_{X}:=(\AA,\EE)$, with vertex set $\AA$ and arrow set
$\EE:=\{(a,a')\in\AA\times\AA:\ [aa']\cap X\neq \emptyset\}$.
In this representation, the potentials 
$\Phi,\Psi$ become `arrow functions' ($2$-symbol potentials)  $\phi,\psi:\EE\to\rr$. 

\ms An {\em elementary circuit} is a cyclic path in $\GG_{X}$ with no repeated vertices.
Any circuit $C$ in $\GG_{X}$ can be represented as a sum of elementary circuits $C=C_1+C_2+\cdots+C_n$, where all the arrows in $C$ appear in one and only one of the elementary circuits $C_i, \ 1\leq i\leq n$, whence
$|C|=\sum_{i=1}^n|C_i|$. 
In general, this representation is not unique.

\ms Given a periodic point $\a\in\per_p(X)$ we define a cyclic path $C(\a):=(\a_0,\a_1,\ldots,\a_{p-1})$ 
in $\GG_{X}$. This path can be represented as the sum of elementary circuits
$C(\a)=C_1+\cdots+C_{n(\sa)}$. Using this representation we can write 
$S_p\phi(\a)\equiv \phi(C(\a)):=\sum_{i=1}^{n(\sa)}\phi(C_i)$, where
$$
\phi(C):=\sum_{i=0}^{|C|-1}\phi(b_i,b_{i+1})
$$
for any circuit $C=(b_0,\ldots,b_{|C|-1})$ in $\GG_{X}$. 
Since $\sum_{i=1}^n|C_i|=|C(\a)|=p$, it follows from \eqref{phibarutile} that 
\begin{eqnarray*}
\bphi & = &  \sup_{p\in\nn}\ \max_{\a\in\sper_p(X)}\frac{\phi(C(\a))}{|C(\a)|}
                      = \sup_{p\in\nn}\ \max_{\a\in\sper_p(A)}\sum_{i=1}^n\frac{|C_i|}{|C(\a)|} 
                                                                               \frac{\phi(C_i)}{|C_i|} \\
                     & = & \max\left\{\frac{\phi(C)}{|C|}:\ C \ \text{is an elementary circuit in } \GG_{X}\right\}\cdot
\end{eqnarray*}
Notice that the set of elementary circuits in $\GG_{X}$ has cardinality bounded by $(\#\AA)!$. 

\ms The maximizing subshift $\bX$ is such that $S_p\phi(\a)=p\bphi$ for all $p\in\nn$ and for all
$\a\in\per_p(\bX)$,  {\em i.e}, $\bX$ is the smallest subshift containing all the periodic points
corresponding to maximizing circuits  in $\GG_{X}$.
Let
$$
\bar{\CC}:=\left\{C\text{ elementary circuit in } \GG_{X}: \ \frac{\phi(C)}{|C|}=\bphi\right\}\cdot
$$
Notice that any maximizing circuit in $\GG_{X}$ is necessarily the sum of circuits 
in $\bar{\CC}$.
On the other hand, if all the arrows in a circuit $C=(b_0,\ldots,b_{|C|-1})$ appear in a 
maximizing circuit, then necessarily $\phi(C)/|C|=\bphi$. From this observation
it follows that 
\[
\bX=\bX(\phi)=\left\{\a\in X:\ \forall \ i\in \zz,\; \exists C\in\bar{\CC} \text{ such that }\, 
(\a_i,\a_{i+1}) \text{ occurs in } C\right\}\cdot
\]
It is clearly a (non-wandering) subshift of finite type defined by a subgraph 
$\GG_{\bX}:=(\bAA,\bEE)\prec\GG_{X}$, where $\bEE$ is the arrow set spanning
the maximizing circuits in $\GG_{X}$ and where $\bAA$ is the set of corresponding
vertices.

\bs The maximizing SFT $\bX$ is a disjoint union of transitive subshifts and we set
$\bX=\bigcup_{J=0}^N\bX_J$. To each transitive component $\bX_J$ it corresponds a transitive subgraph 
$\GG_J:=(\bAA_J,\bEE_J)\prec \GG_{X}$. Let us order these transitive components so that 
$P(\psi|\bX_J)=P(\psi|\bX)$ for each $1\leq J\leq N_\phi$ and $P(\psi|\bX_J)<P(\psi|\bX)$
for each $N_\phi<J\leq N$. (Recall that $P(\psi|\bX)=\max\{P(\psi|\bX_J): 1\leq J\leq N\}$.)

\begin{definition}[Heavy components]\label{defheavy}
We will refer to the SFT's $\bX_J$ in the subcollection $\{\bX_J\subset\bX:\ 1\leq J\leq N_\phi\}$ as 
{\it the heavy components of $\bX$}. 
We define the set
\begin{equation}\label{Ebarphi}
\bEE_{\phi}:=\bigsqcup_{J=1}^{N_\phi}\bEE_J\subseteq\bEE
\end{equation}
of all arrows in digraphs associated to heavy components.
\end{definition}
As we shall see, the zero-temperature limit is concentrated {\em only on some} of the heavy components. 

\ms From now on, and without loss of generality, we assume that we are given a directed graph
$\GG_{X}:=(\AA,\EE)$ and two potentials $\phi,\psi:\EE\to\rr$. 

\ms \subsection{Transition Matrices and Equilibrium States}\

\ms Let $\phi,\psi:\EE\to\rr$.
To each $\beta\in\rr$ it corresponds a unique equilibrium state $\mu_{\bpp}$, which is a ($1$-step) Markov measure completely determined by the irreducible matrix
$\MM_{\bpp}:\AA\times \AA\to [0,\infty)$ 
defined by
\begin{equation}\label{matrix}
\MM_{\bpp}\left(a,a'\right)=\left\{\begin{array}{ll} 
                   e^{\bp(a,a')+\psi(a,a')} & \text{ if } (a,a')\in \EE,\\
                   0 & \text{ otherwise.} \end{array}\right.
\end{equation}
According to Perron--Frobenius Theorem (see Theorem~\ref{perron-frobenius} in 
Appendix~\ref{proof-of-periodic}) there are unique left and right maximal eigenvectors $\v_{\bpp}$ and 
$\w_{\bpp}$ associated to the maximal eigenvalue $\rho_{\bpp}:=\max|\spec(\MM_{\bpp})|$.
We can choose them in such a way that $\w_{\bpp}^{\dag}\v_{\bpp}=1$. 
We then have the formula
\begin{equation}\label{equilibrium-state-mubpp}
\mu_{\bpp}[\b_0^n]= \frac{1}{\rho_{\bpp}^{n+1}}\thinspace\w_{\bpp}(\b_0)\left(\prod_{i=0}^{n-1}
                                                      \MM_{\bpp}(\b_i,\b_{i+1})\right)\v_{\bpp}(\b_n),
\end{equation}
for every $\b_0^n\in X$ ($n\in\nn_0$). 

\ms For each $1\leq J\leq N$, let $\bAA_J\subset \AA$ be the vertex set of the digraph 
$\GG_J:=(\bAA_J,\bEE_J)$ associated to the transitive component $\bX_J\subset \bX$.
The transition matrices $\MM_{\psi,J}:\bAA_J\times \bAA_J\to\rr^+$ defined by
\begin{equation}\label{psiJ-matrix}
\MM_{\psi,J}(a,a')      = \left\{\begin{array}{ll}
                           e^{\psi(a,a')} & \text{ if } (a,a')\in \bEE_J\\
                                        0 & \text{ otherwise} \end{array}
                             \right. 
\end{equation}
are irreducible. Therefore they have associated to their maximal eigenvalue 
$\rho_{\psi,J}:=\max|\spec(\MM_{\psi,J})|$ unique left and right eigenvectors
$\v_{\psi,J},\w_{\psi,J}:\bAA_J\to\rr^+$ such that $\w_{\psi,J}^{\dag}\v_{\psi,J}=1$.
We can associate to each transitive component $\bX_J\subset \bX$  the Markov measure $\nu_{\psi,J}$ 
defined by
\begin{equation}\label{nu-J}
\nu_{\psi,J}[\b_0^n]:=\frac{1}{\rho_{\psi,J}^{n+1}}\thinspace\w_{\psi,J}(\b_0)
\left(\prod_{i=0}^{n-1}\MM_{\psi,J}(\b_i,\b_{i+1})\right)
\v_{\psi,J}(\b_n),
\end{equation}
for each $\b\in X$ such that $[\b_0^n]\cap \bX_J\neq\emptyset$.
This is precisely the equilibrium state on $\bX_J$ associated to the potential $\psi|\bX_J$.
We recall that $P(\psi|\bX_J)=\log \rho_{\psi,J}$. We of course have $\nu_{\psi,J}(\bX_J)=1$.

\medskip

\ms Without loss of generality, our original question can be recast in the following way.
We slightly generalize it and consider potentials of the form $\beta\phi+\psi$ which will appear naturally when we make the renormalization described in the next section.
\begin{question}[Recasting Question \ref{question}]\label{newquestion}
Let $X\subset \AA^\zz$ be a topological Markov chain defined by a finite alphabet $\AA$
and an arrow set $\EE$ and let $\phi,\psi:\EE\to\rr$ be potentials. Denote by
$\mu_{\bpp}$ the unique equilibrium state of $\bpp$ for each $\beta>0$, which is a Markov measure.\\
Does the limit of $(\mu_{\bpp})_{\beta>0}$ exist in weak-$^*$ topology ? Can we compute its barycentric
decomposition on the finite-dimensional simplex generated by the ergodic $\phi$-maximizing measures ?
\end{question}

\ms \subsection{Normalization of $\phi$ and $\psi$}\label{normalisationpot}    \

\ms It will be convenient later on to assume that $\bphi=P(\psi|\bX)=0$.
If it is not the case, then one can redefine $\phi$ as $\phi-\bphi$.
These potentials are cohomologous and yield the same Gibbs state $\mu_{\beta\phi}$.
Moreover, this does not change the $\phi$-maximizing SFT $\bX$. 
One can also redefine $\psi$ as $\psi-P(\psi|\bX)$ without changing 
the equilibrium states $\nu_{\psi,J}$. Finally, for every $\beta$, 
we can assume that the potential $\beta\phi+\psi$ is such
that $\bphi=0$ and $P(\psi|\bX)=0$, otherwise we normalize $\phi$ and $\psi$
as above and obtain the same equilibrium state $\mu_{\bpp}$.
Since heavy components (Definition \ref{defheavy}) maximize
$P(\psi|\bX_J)$, $J=1\ldots N$, and since $P(\psi|\bX)=\max\{P(\psi|\bX_J): 1\leq J\leq N\}$,
then, if $P(\psi|\bX)=0$, we have $P(\psi|\bX_J)<0$, $N_\phi<J\leq N$.
Since $P(\psi|\bX_J)=\log\rho_{\psi,J}$, this means that
$\rho_{\psi,J}\leq 1$ for each $1\leq J\leq N$, and $\rho_{\psi,J}=1$ for the
heavy components $1\leq J\leq N_\phi$.

\ms \subsection{Notations}\

\ms Henceforth, we will use the following convenient short-hand notations:
\begin{eqnarray*}
a=\pm b & \text{for} & -b\leq a\leq b \\
a=e^{\pm b} & \text{for} & e^{-b}\leq a\leq e^b
\end{eqnarray*}
where $a,b$ are positive real numbers.

\bs \section{Main Theorem}\label{main} \

\ms We formulate our main result (Theorem \ref{main-theorem}) which answers Question \ref{newquestion} (hence Question \ref{question}). For convenience, the algorithm describing how to compute the zero-temperature
limit is in Subsection \ref{thealgo}.

\ms \subsection{Convergence when temperature goes to zero}\label{thetheorem} \

\begin{theorem}\label{main-theorem}
Let $X\subset \AA^\zz$ be a topological Markov chain defined by a finite alphabet $\AA$
and an arrow set $\EE$ and let $\phi,\psi:\EE\to\rr$ be potentials.
Denote by $\mu_{\bpp}$ the unique equilibrium state of $\bpp$ for each $\beta>0$.
Let $\bX$ be the $\phi$--maximizing SFT and $\bX_1,\ldots, \bX_{N_\phi}\subset\bX$
its heavy components.
For each heavy component $\bX_J$, $1\leq J\leq N_{\phi}$, let $\nu_{\psi,J}$
be the (unique) equilibrium state on $\bX_J$ associated to the potential $\psi$ (restricted to $\bX_J$).
\ms Then the sequence $(\mu_{\bpp})_{\beta>0}$ converges in the weak-$^*$ topology and
\[
\lim_{\beta\to+\infty}\mu_{\bpp}=\sum_{J=1}^{N_{\phi}} \alpha_J\thinspace\nu_{\psi,J}
\]
where $0\leq \alpha_J\leq 1$ and $\sum_{J=1}^{N_\phi} \alpha_J=1$.

\ms The $\alpha_J$'s can be computed by means of
a recursive algorithm, detailed below, which converges after a finite number of steps.
Notice that some of the $\alpha_J$'s may be zero.

\ms Furthermore, there exist $\beta_0=\beta_0(\AA,\phi,\psi)$ and $C=C(\AA,\phi)>0$ (independent of $\psi$)
such that for all $\beta\geq \beta_0$
\begin{equation}\label{speed}
\mu_{\bpp}[\b_0^n]=\alpha_J\thinspace\nu_{\psi,J}[\b_0^n]\pm e^{-C\beta},
\end{equation}
whenever $[\b_0^n]\cap \bX_J\neq \emptyset$ for some $1\leq J\leq N_\phi$.
\end{theorem}


\ms \subsection{The Algorithm to compute the zero-tempeature limit}\label{thealgo} \

\ms To describe the algorithm for computing the coefficients $\alpha_J$, we need more
notations and definitions. We assume that $\bphi=P(\psi|\bX)=0$ (see Subsection \ref{normalisationpot}).

\bigskip
\begin{definition}[Renormalized SFT]\label{newsft}
To each heavy component $\bX_J$, $1\leq J\leq N_\phi$, it corresponds a sub-alphabet
$\bAA_J\subset\AA$ such that $\bX_J\subset (\bAA_J)^{\zz}$.  
We define the new alphabet $\AA':=\{1,2,\ldots,N_\phi\}$.
The new arrow set $\EE'\subset \AA'\times\AA'$ is
defined as follows: $(J,K)\in\EE'$  if and only if there exists $\a\in X$ and $n\in\nn$ such that
$\a_0\in \bAA_J$, $\a_n\in\bAA_K$ and $\bigcup_{L=1}^{N_\phi}(\{\a_m:\ 1 < m <n\}\cap\bAA_L)=\emptyset$.
Then the ``renormalized'' SFT is defined by
$$
X':=\{\a\in \AA'\,{}^{\zz}:\ (\a_n,\a_{n+1})\in \EE',\ \forall n\in\zz\}.
$$
\end{definition}

\bs Let $\MM_{\psi}:\AA\times \AA\to\rr^+$ be defined in the same way as $\MM_{\bpp}$, with
$\bpp$ instead of $\psi$, and let $\bMM_{\psi}$ be the restriction of $\MM_{\psi}$ to all the 
transitive components of $\bX$, either heavy or not.
It is easy to verify that $\sum_{k=0}^{\infty}\bMM_{\psi}^k(b,b) < \infty$,
for each $b\in \AA\setminus \bigsqcup_{J=1}^{N_\phi}\bAA_J$.

\ms We need some notations for paths in digraphs.

\begin{notations}[Paths in digraphs]
Given two vertices $a,c$ in a digraph $\GG$, we denote by
$a\lds c$ any finite path starting from $a$ and ending at $c$.
When the path is elementary \textup{(}{\em i.e.}, it does not contain any
circuit\textup{)} we write $a\to c$. The symbols `$\lds$' and `$\to$'
will be naturally used as variables in path-depending functions.
We will also use the notation $b \cira$ to denote a finite circuit based on $b$.
\end{notations}

Given $a,c\in \AA$, we denote by $\path[a,c]$ the collection of all elementary
paths in $\GG_{X}$ (the digraph defining $X$) going from $a$ to $c$.
To a given path of the form $a\to c=(a,b_1,\ldots,b_{m},c)\in \path[a,c]$, with $m\in\nn$
and $b_i\in \AA\setminus\bigsqcup_{J=1}^{N_\phi}\bAA_J$ ($1\leq i\leq m$),
we associate the following {\it transition pressure}
\begin{equation}\label{transition-pressure}
P_\psi(a\to c):=\log\left(\sum_{i=1}^{m}\sum_{k=0}^{\infty}\bMM_{\psi}^k(b_i,b_i)\right)<+\infty.
\end{equation}

\ms Let us now fix, for each heavy component $1\leq J\leq N_\phi$, a {\it central vertex} $c_J\in \bAA_J$. 
Since $\phi(C)=0$ for each circuit $C$ in $\GG_J$ (the digraph defining $\bX_J$), it follows that 
$\phi\left({c_J\lds a}\right)$ has the same value for all paths $c_J\lds a$ in $\GG_J$ connecting
$c_J$ to $a$. Therefore we can define, for each $a\in \bAA_J$, the {\it central term}
\begin{equation}\label{central-term}
\lfphi(a):=\phi\left({c_J\lds a}\right) \text{ with } c_J\lds a\text{ a path in } \GG_J \text{ from } 
c_J \text{ to } a.
\end{equation}

\ms For $a\in \bAA_J$ and $c\in \bAA_K$, let $\tpath[a,c]$ be the set of all elementary paths in $\GG_X$, starting 
at $a$ and ending at $c$, with no arrows in $\bEE_{\phi}$. 
We define the {\em transition term}
\begin{equation}\label{transition-term}
\rtphi(a,c):=\max\big\{\phi\left({a\to c}\right):\ a\to c\in \tpath[a,c]\big\},
\end{equation}
where $\phi\left({a\to c}\right):=\phi(a,b_1)+\sum_{i=1}^{m-1}\phi(b_i,b_{i+1})+\phi(b_{m},c)$ for 
$a\to c=(a,b_1,\ldots,b_{m},c)$. 

\begin{definition}[Renormalized potentials]\label{newpotentials}
With the notations just introduced, we define the {\it renormalized potentials} $\phi',\psi':\EE'\to\rr$
by
$$
\phi'(J,K):=\max_{a\in \bAA_J, c\in\bAA_K}\big\{\lfphi(a)+\rtphi(a,c)-\lfphi(c)\big\}
$$
$$
\psi'(J,K):=
$$
$$
\log\left(
\sum_{a,c\in \bAA_{J,K}}  \v_{\psi,J}(a)\w_{\psi,K}(c) \sum_{a\to c\in\bpath[a,c]}
\exp\big(\psi\left({a \to c}\right)+P_\psi(a\to c)\big)
\right), 
$$
where:
\begin{itemize}
\item $\bAA_{J,K}\subset \bAA_J\times\bAA_K$ denotes the set where
$\lfphi(a)+\rtphi(a,c)-\lfphi(c)$ is maximal,
\item $\bpath[a,c]$ is the set of all elementary paths from $a$ to $c$ maximizing $\phi$, and
\item
$P_\psi(a\to c)$ is defined in \eqref{transition-pressure}.
\end{itemize}
\end{definition}

\begin{algo}[Continuation of Theorem \ref{main-theorem}]
\ms \begin{enumerate}
\item 
Let $X'\subset \AA'\,{}^{\zz}$ be the renormalized subshift of finite type, and $\phi',\psi':\EE'\to\rr$ be
the renormalized two--symbol potentials defined above. 

\item Compute $\bX'\subset X'$, the $\phi'$--maximizing SFT,
and normalize $\phi'$, {\em i.e.}, replace $\phi'$ by $\phi'-\bphi'$. 

\item Normalize $\psi'$, {\em i.e.}, replace $\psi'$ by $\psi'-P(\psi'|\bX')$, and 
identify the heavy components $\bX'_1,\ldots, \bX'_{N_{\phi'}}$ in $\bX'$. 

\item Compute, for each $1\leq J\leq N_{\phi'}$, the equilibrium state $\nu_{\psi',J}$ associated to the potential 
$\psi'$ on $\bX'_J$ (which is a Markov measure).

\item If $J\in\AA'$ is such that $[J]\cap (\bigcup_{K=1}^{N_{\phi'}}\bX'_K)=\emptyset$, 
then $\alpha_J=0$. Otherwise $\alpha_J=\nu_{{\psi'},K}[J]\thinspace \alpha'_K$, where $K\leq N_{\phi'}$ 
is such that $[J]\cap \bX'_K\neq\emptyset$, and $\alpha'_K\geq 0$ is computed following steps (1) to (5) using
$X'$ instead of $X$, $\phi'$ instead of $\phi$ and $\psi'$ instead of $\psi$.
\end{enumerate}
\end{algo}

\ms \subsection{Comparison with previous results}\

\ms Let us compare our result with those in \cite{Bremont03,Leplaideur05}.
In \cite{Bremont03}, it is proved that the limit exists but it is not
identified. In \cite{Leplaideur05}, the limit is proved to exist
by a somewhat indirect argument and the coefficients of the
barycentric decomposition of the limiting measure are not
explicitly identified. Here we directly prove the existence of
the limit and identify it precisely. 
Besides the fact that we make explicit calculations, the main idea is renormalization.
The reader can see our method at work in Section \ref{examples}.
Finally, we are able to prove that the Gibbs measure concentrates
exponentially fast in $\beta$ about its limit. 

\bs \section{The Renormalization Lemma: Proof of Theorem \ref{main-theorem}}\label{mainproof} \

\ms From now on, by `for $\beta$ large enough' we mean that `there exists
$\beta_0=\beta_0(\AA,\phi,\psi)$ such that for all $\beta\geq \beta_0$'.

\ms The following lemma is the crux of our approach.

\ms \begin{renorm}
Let $\AA$ be a finite alphabet and $X\subset \AA^{\zz}$  a transitive topological Markov chain
defined by an arrow set $\EE$. 
Let $\phi,\psi:\EE\to\rr$ be $2$-symbol potentials
and $\bX:=\bigsqcup_{J=1}^N\bX_J$ be the $\phi$--maximizing SFT which decomposes into transitive components
$\bX_J$.
Let $\phi',\psi':\EE'\to\rr$ be the renormalized potentials (cf. Definition \ref{newpotentials}),
where $\EE'$ is the arrow set of the renormalized SFT $X'$ (Definition \ref{newsft}),
and $\mu_{\beta\phi'+\psi'}$ the equilibrium state of $\beta\phi'+\psi'$.\\
Then there exists $\delta=\delta(\AA,\phi)$ (independent of $\psi$) such that for all $\beta$ large enough
\[
\mu_{\bpp}[\b_0^n]=\nu_{\psi,J}[\b_0^n]\thinspace\mu_{\beta\phi'+\phi'}[J] \pm e^{-\beta\, \delta},
\]
for all cylinder $[\b_0^n]$ intersecting an heavy component (Definition \ref{defheavy}), {\em i.e.},
whenever $[\b_0^n]\cap \bX_J\neq \emptyset$ for some $1\leq J\leq N_\phi$.
\end{renorm}

\ms The proof of Theorem \ref{main-theorem}
consists in the recursive application of the Renormalization Lemma, which 
we assume true for the moment and which we shall prove in Section~\ref{thebiglemma}. 

\ms \paragraph{{\bf \em Step 1} (First Renormalization)} We compute the $\phi$--maximizing subshift 
$\bX\subset X$, which we decompose into its transitive components $\bX:=\bigsqcup_{J=1}^N\bX_J$.
We order the heavy components so that they are indexed by $J=1,\ldots,N_\phi$.
According to the Renormalization Lemma, there exists $\delta$  such that for $\beta$ large enough
\[
\mu_{\bpp}[\b_0^n]=\nu_{\psi,J}[\b_0^n]\thinspace\mu_{\beta\phi'+\phi'}[J] \pm e^{-\beta\, \delta},
\]
whenever $[\b_0^n]\cap \bX_J\neq \emptyset$ for some $1\leq J\leq N_\phi$.
In particular, for each $b\notin\bigsqcup_{J=1}^{N_\phi}\bAA_J$ and $\beta$ large enough we have
\begin{eqnarray*}
\mu_{\bpp}[b]&\leq& 1-\sum_{J=1}^{N_\phi}\sum_{a\in\bAA_J}\mu_{\bpp}[a]\\
             &\leq&
             1-\sum_{J=1}^{N_\phi}\mu_{\beta\phi'+\phi'}[J]\left(\sum_{a\in\bAA_J}\nu_{\psi,J}[a]\right)
                    +\#\left(\bigsqcup_{J=1}^{N_\phi}\bAA_J\right) e^{-\beta\, \delta}\\
             &\leq& \#\AA\thinspace e^{-\beta\, \delta}.       
\end{eqnarray*}
Hence, it follows, for $\beta$ large enough, that 
\begin{equation}\label{chou}
\mu_{\bpp}[\b_0^n]\leq \thinspace e^{-\beta\, \frac{\delta}{2}}
\end{equation}
whenever $[\b_0^n]\cap \bX_J\neq \emptyset$ for some $1\leq J\leq N_\phi$.
Indeed, such a cylinder contains at least a letter $b\notin\bigsqcup_{J=1}^{N_\phi}\bAA_J$.
By shift-invariance, we can assume that $\b_0=b$ and we have $\mu_{\bpp}[\b_0^n]\leq \mu_{\bpp}[b]$.

\ms \paragraph{{\bf \em Step 2} (Recursion)} Now, in order to compute $\mu_{\beta\phi'+\phi'}([J])$, 
we apply {\it Step 1} to the renormalized system $(X',\beta\phi'+\psi')$ and use \eqref{chou}. 
This yields
\[
\begin{array}{ll}
\mu_{\beta\phi'+\psi'}[J]= 
                                       \nu_{\phi',K}[J]\thinspace\mu_{\beta\phi''+\phi''}[K] 
                                       \pm e^{-\beta\, \delta'} & \text{ if } J\in\bAA'_K,\, K\leq N_{\phi'},\\ \\
\mu_{\beta\phi'+\psi'}[J]\leq e^{-\beta\, \frac{\delta'}{2}} & \text{ if } J\notin\bigsqcup_{K=1}^{N_{\phi'}}\bAA'_J,
\end{array} 
\]
for $\beta$ large enough.\\
In particular, $\nu_{{\psi'},K}[J]$ is the first term of the factorization of $\alpha_J$.\\ 
When we iterate this procedure $i$ times, we obtain a sequence of renormalized systems, the cumulative error
term, and the first $i$ terms of the factorization of $\alpha_J$, for each $1\leq J\leq N_\phi$.

\ms \paragraph{{\bf \em Final Step} (Convergence of the Recursion)}
Suppose we have done $i$ times the renormalization, so that we have
a renormalized system $(X^{(i)},\beta\phi^{(i)}+\psi^{(i)})$ with a corresponding
$\delta^{(i)}$. Then there are two cases:
either $\#\AA^{(i)} < \#\AA^{(i-1)}$,
or $\#\AA^{(i)}=\#\AA^{(i-1)}$.\\
In the second case, we have a maximizing SFT made only of fixed points,
{\em i.e.}, each symbol of the $(i-1)$-th normalization defines a heavy component. 
Then we necessarily have $(J,J)\in\EE^{(i-1)}$ for each $J\in\AA^{(i-1)}$, whence
$\EE^{(i)}=\EE^{(i-1)}\setminus\{(J,J)\in\EE^{(i-1)}\}$.
Therefore, if we are in the second case, then
$\#\AA^{(i+1)}<\#\AA^{(i-1)}$.
The renormalization process ends because there exists some 
$m=m(\AA,\phi,\psi)$ such that $\AA^{(m)}$ is a singleton, say $\{1\}$, and
we necessarily have $\mu_{\beta\phi^{(m)}+\psi^{(m)}}[1]=1$.
Therefore, for $\beta$ large enough, we end up with 
\[
\mu_{\bpp}[\b_0^n]=\alpha_J\nu_{\psi,J}[\b_0^n]\pm \left(\sum_{j=1}^{m-1}e^{-\beta\delta^{(j)}}\right)
\]
whenever $[\b_0^n]\cap \bX_J\neq \emptyset$ for some $1\leq J\leq N_\phi$.\\
The last claim of the theorem follows from this, by taking 
$C:=\delta^{(m-1)}/2$ and $\beta$ large enough.
Theorem \ref{main-theorem} is proved.
\endproof

\bs \section{Examples}\label{examples} \

\ms \subsection{A Basic Example}\

\ms  Let $X=A^{\zz}$ with $A=\{a,b,c\}$ and the following $2$-symbol potential:
\[
\phi:=\left(\begin{matrix} 0 & -1 & -2\\ -1 & 0 & -2\\-2 & -2 & 0 \end{matrix}\right).
\]
One can of course compute $\lim_{\beta\to+\infty}\mu_{\bp}$ directly, as was done in \cite{Jenkinson01}.
For the sake of illustration of our method, let us compute it by following the algorithm described in 
Section \ref{main}. 
In this case $\bX=\per_1(X):=\{\a,\b,\c\}$, with $\a,\b,\c$ such that $\a_n=a, \b_n=b,\c_n=c$, for all $n\in\zz$.
Hence
\[
\lim_{\beta\to+\infty}\mu_{\bp}=\alpha_1\delta_{\sa}+\alpha_2\delta_{\sb}+\alpha_3\delta_{\sc},
\]
where $\delta_{\sx}$ denotes the Dirac measure at $\x$.
Notice that $\bphi=0$. Since $\psi\equiv 0$ and $\bX$ is a finite union of periodic points,
we have $P(\psi|\bX)=h_{\top}(\bX)=0$. Hence potentials are already normalized.
The renormalized alphabet is $A'=\{1,2,3\}$, and the renormalized system is the topological Markov chain $X'\subset A'\,{}^{\zz}$ described by the digraph $(A',\EE')$ shown in the following picture:

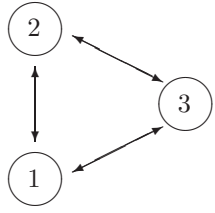
\begin{center}

\unitlength=1truecm
\begin{picture}(5,4)(0,0)
\put(1,1){\circle{0.7}}
\put(0.9,0.9){$1$}
\put(1,3){\circle{0.7}}
\put(0.9,2.9){$2$}
\put(3,2){\circle{0.7}}
\put(2.9,1.9){$3$}

\put(1,1.5){\vector(0,1){1}}
\put(1,2.5){\vector(0,-1){1}}

\put(1.5,1.1){\vector(2,1){1.2}}
\put(2.5,1.6){\vector(-2,-1){1}}

\put(1.5,2.9){\vector(2,-1){1.2}}
\put(2.5,2.4){\vector(-2,1){1}}


\end{picture}

\end{center}

\ms In this case, the renormalized potentials $\phi',\psi':\EE'\to\rr$ are 
given by
\[
\phi':=\left(\begin{matrix} -\infty & -1 & -2\\ -1 & -\infty & -2\\-2 & -2 & -\infty \end{matrix}\right)
\]
and 
\[
\psi':=\left(\begin{matrix} -\infty & 0 & 0\\ 0 & -\infty & 0\\ 0 & 0 & -\infty \end{matrix}\right).
\]
(By `$-\infty$' we mean that there is no arrow.)
The maximizing topological Markov chain $\bX'$ reduces to the periodic orbit of $\x:=(\ldots1212\ldots)$,
which is the only heavy component, and carries only one shift-invariant measure, namely
$\frac12(\delta_{\sx}+\delta_{\sigma\sx})$. Hence $\alpha_1=\alpha_2=\frac12$ and $\alpha_3=0$
and 
\[
\lim_{\beta\to+\infty}\mu_{\bp'}=
\frac12(\delta_{\sx}+\delta_{\sigma\sx}).
\]
Therefore the limit of the original measure is
\[
\lim_{\beta\to+\infty}\mu_{\bp}=
\frac12(\delta_{\sa}+\delta_{\sb}).
\]

\ms \subsection{An Example with an Irrational Barycenter} \

\ms Let $X=A^{\zz}$ with $A=\{a,b,c,d\}$ and the following $2$-symbol potential:
\[
\phi:=\left(\begin{matrix} 0 & -1 & -1 & -2\\ -1 & 0 & -1 & -2\\-1 & -1 & 0 & -1\\ -2 & -2 & -1 & 0\end{matrix}\right).
\]
In this case $\bX=\per_1(X):=\{\a,\b,\c,\d\}$, with $\a,\b,\c$ and $\d$ such that $\a_n=a, \b_n=b,\c_n=c$ 
and $\d_n=d$ for all $n\in\zz$. Hence
\[
\lim_{\beta\to+\infty}\mu_{\bp}=\alpha_1\delta_{\sa}+\alpha_2\delta_{\sb}+\alpha_3\delta_{\sc}+\alpha_4\delta_{\sd}.
\]
Note that potentials are already normalized.
The renormalized alphabet is $A'=\{1,2,3,4\}$, and the renormalized system is the topological Markov chain $X'\subset A'\,{}^{\zz}$ described by the digraph $(A',\EE')$ shown in the following picture:

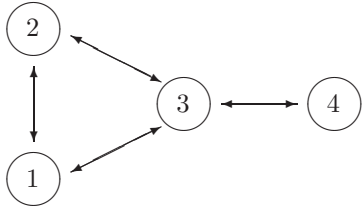
\begin{center}

\unitlength=1truecm
\begin{picture}(6,4)(0,0)
\put(1,1){\circle{0.7}}
\put(0.9,0.9){$1$}
\put(1,3){\circle{0.7}}
\put(0.9,2.9){$2$}
\put(3,2){\circle{0.7}}
\put(2.9,1.9){$3$}
\put(5,2){\circle{0.7}}
\put(4.9,1.9){$4$}

\put(1,1.5){\vector(0,1){1}}
\put(1,2.5){\vector(0,-1){1}}

\put(1.5,1.1){\vector(2,1){1.2}}
\put(2.5,1.6){\vector(-2,-1){1}}

\put(1.5,2.9){\vector(2,-1){1.2}}
\put(2.5,2.4){\vector(-2,1){1}}

\put(3.5,2){\vector(1,0){1}}
\put(4.5,2){\vector(-1,0){1}}

\end{picture}

\end{center}

\ms In this case, the renormalized potentials $\phi',\psi':\EE'\to\rr$ are constant and such that $\phi'(J,K)=-1$
and $\psi'(J,K)=0$, for each arrow $(J,K)\in \EE'$. The maximizing topological Markov chain $\bX'$ coincides with the renormalized topological Markov chain, therefore there is only one heavy component $\bX'=X'$, and the second renormalization is trivial, {\em i.e.}, the resulting topological Markov chain is a fixed point. From this we obtain $\alpha_J=\nu_{X'}[J]$, with $\nu_{X'}$ the measure of maximal entropy
on $X'$. It can be explicitly computed, and we finally obtain
\[
\lim_{\beta\to+\infty}\mu_{\bp}=\frac{\delta_{\sa}+\delta_{\sb}}{2(4-\rho)}+
\frac{(\rho-1)^2\delta_{\sc}}{2(4-\rho)}+\frac{(\rho-1)^2\delta_{\sd}}{2\rho^2(4-\rho)},
\]
with $\rho=\frac13(1+2\sqrt{10}\cos(\frac13\arctan(3\sqrt{111})))$, the largest root of the polynomial
$p(x)=x^4-4x^2-2x+1$, which is an irrational number. In the table below we present the comparison between the limiting
measure and $\mu_{\bp}$ for different values of the inverse temperature.

\ms 
\begin{center}
\scriptsize{
\begin{tabular}{|r|c|c|c|c|c|c|c|c|}
\hline
      $\beta=$  &$\log(2)$  &$2\log(2)$&$3\log(2)$&$4\log(2)$&$5\log(2)$&$6\log(2)$& $\cdots$ &$\infty$\\
\hline
$\mu_{\bp}[a]=$ &  0.253298 & 0.259815 & 0.265413 & 0.269011 & 0.271041 & 0.272118 & $\cdots$ &  0.273237\\
$\mu_{\bp}[b]=$ &  0.253298 & 0.259815 & 0.265413 & 0.269011 & 0.271041 & 0.272118 & $\cdots$ &  0.273237\\
$\mu_{\bp}[c]=$ &  0.316672 & 0.349361 & 0.363356 & 0.369239 & 0.371810 & 0.372988 & $\cdots$ &  0.374089\\
$\mu_{\bp}[d]=$ &  0.176732 & 0.131010 & 0.105818 & 0.092738 & 0.086109 & 0.082777 & $\cdots$ &  0.079437\\
\hline
\end{tabular}
}
\end{center}

\ms \subsection{An Example with a Two--step Renormalization} \

\ms Let us now take $X=A^{\zz}$ with $A=\{a,b,c,d,e\}$, and consider the following $2$-symbol potential:
\[
\phi:=\left(\begin{matrix} 0 & -4 & -1 & -3 & -4\\ 
                          -1 & 0  & -4 & -3 & -3\\
                          -4 & -1 &  0 & -3 & -3\\ 
                          -4 & -4 & -4 &  0 & -1\\
                          -3 & -4 & -4 & -1 & 0 \end{matrix}\right).
\]
We have
$\bX=\per_1(X):=\{\a,\b,\c,\d,\e\}$, with $\a,\b,\c,\d$ and $\d$ such that $\a_n=a, \b_n=b,\c_n=c,\d_n=d$ 
and $\e_n=e$ for all $n\in\zz$, {\em i.e},
\[
\lim_{\beta\to+\infty}\mu_{\bp}=\alpha_1\delta_{\sa}+\alpha_2\delta_{\sb}+\alpha_3\delta_{\sc}+\alpha_4\delta_{\sd}+\alpha_5\delta_{\se}.
\]
As in the previous examples, potentials are already normalized.
The renormalized alphabet is $A'=\{1,2,3,4,5\}$, and the renormalized system is the topological Markov chain $X'\subset A'\,{}^{\zz}$ described by
the digraph $(A',\EE')$ shown in the following picture.

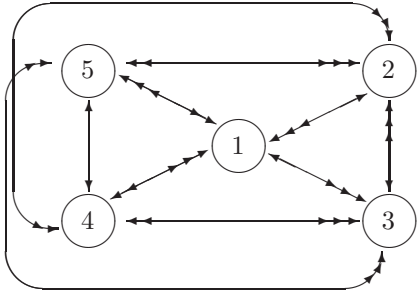
\begin{center}

\unitlength=1truecm
\begin{picture}(6,4)(0,0)
\put(3,2){\circle{0.7}}
\put(2.9,1.9){$1$}
\put(5,1){\circle{0.7}}
\put(4.9,0.9){$3$}
\put(5,3){\circle{0.7}}
\put(4.9,2.9){$2$}

\put(1,1){\circle{0.7}}
\put(0.9,0.9){$4$}
\put(1,3){\circle{0.7}}
\put(0.9,2.9){$5$}

\put(3.5,2.1){\vector(2,1){1.2}}
\put(4.6,1.3){\vector(-2,1){1.2}}
\put(5,2.6){\vector(0,-1){1.2}}
\put(3.9,2.3){\vector(-2,-1){0.1}}
\put(3.7,2.2){\vector(-2,-1){0.1}}
\put(3.5,2.1){\vector(-2,-1){0.1}}
\put(5,2.5){\vector(0,1){0.1}}
\put(5,2.3){\vector(0,1){0.1}}
\put(5,2.1){\vector(0,1){0.1}}
\put(4.6,1.3){\vector(2,-1){0.1}}
\put(4.4,1.4){\vector(2,-1){0.1}}
\put(4.2,1.5){\vector(2,-1){0.1}}

\put(1,2.6){\vector(0,-1){1.2}}
\put(1,1.4){\vector(0,1){1.2}}

\put(1.5,2.9){\vector(2,-1){1.2}}
\put(1.5,2.9){\vector(2,-1){1}}
\put(2.5,1.9){\vector(-2,-1){1.2}}
\put(2.5,1.9){\vector(-2,-1){1}}
\put(1.5,2.9){\vector(-2,1){0.1}}
\put(1.7,2.8){\vector(-2,1){0.1}}
\put(1.9,2.7){\vector(-2,1){0.1}}
\put(2.5,1.9){\vector(2,1){0.1}}
\put(2.3,1.8){\vector(2,1){0.1}}
\put(2.1,1.7){\vector(2,1){0.1}}

\put(4.5,1){\vector(-1,0){3}}
\put(4.5,1){\vector(-1,0){2.8}}
\put(4.5,1){\vector(1,0){0.1}}
\put(4.3,1){\vector(1,0){0.1}}
\put(4.1,1){\vector(1,0){0.1}}

\put(4.5,3.1){\vector(-1,0){3}}
\put(4.5,3.1){\vector(-1,0){2.8}}
\put(4.5,3.1){\vector(1,0){0.1}}
\put(4.3,3.1){\vector(1,0){0.1}}
\put(4.1,3.1){\vector(1,0){0.1}}

\put(2.5,3.4){\oval(5,1)[t]}
\put(0.5,2.4){\oval(1,3)[l]}
\put(0.5,0.9){\vector(1,0){0.1}}
\put(0.3,0.95){\vector(2,-1){0.1}}
\put(4.9,0.5){\vector(0,1){0.1}}
\put(4.8,0.25){\vector(1,2){0.1}}
\put(4.7,0.2){\vector(1,1){0.1}}

\put(2.4,0.6){\oval(5,1)[b]}
\put(0.4,1.6){\oval(1,3)[l]}
\put(0.4,3.1){\vector(1,0){0.1}}
\put(0.2,3.05){\vector(2,1){0.1}}
\put(5,3.5){\vector(0,-1){0.1}}
\put(4.9,3.75){\vector(1,-2){0.1}}
\put(4.8,3.8){\vector(1,-1){0.1}}

\end{picture}

\end{center}

\ms In the previous figure, we use the following convention: an arrowhead from $J$ to $K$ 
means that $\phi'(J,K)=-1$, a double arrowhead corresponds to $\phi'(i,j)=-2$, while a triple arrowhead
means that $\phi'(J,K)=-3$.
Since $P(\psi|\bX_J)=0$ for each heavy component $1\leq J\leq 4$,  
$\psi'(J,K)=0$, for each arrow $(J,K)\in \EE'$. The maximizing SFT of the first renormalization $\bX'$ is composed by two heavy components as indicated in the picture below.

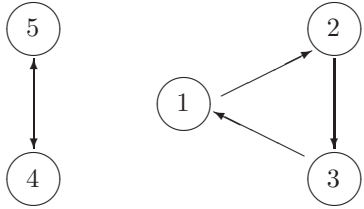
\begin{center}

\unitlength=1truecm
\begin{picture}(6,4)(0,0)
\put(3,2){\circle{0.7}}
\put(2.9,1.9){$1$}
\put(5,1){\circle{0.7}}
\put(4.9,0.9){$3$}
\put(5,3){\circle{0.7}}
\put(4.9,2.9){$2$}

\put(1,1){\circle{0.7}}
\put(0.9,0.9){$4$}
\put(1,3){\circle{0.7}}
\put(0.9,2.9){$5$}

\put(3.5,2.1){\vector(2,1){1.2}}
\put(4.6,1.3){\vector(-2,1){1.2}}
\put(5,2.6){\vector(0,-1){1.2}}

\put(1,2.6){\vector(0,-1){1.2}}
\put(1,1.4){\vector(0,1){1.2}}

\end{picture}

\end{center}
Hence, according to the algorithm, the second renormalization yields
\[
\lim_{\beta\to+\infty}\mu_{\bp'+\psi'}=\alpha'_1\frac{\delta_{\sx}+\delta_{\sigma\sx}}{2}
+\alpha'_2\frac{\delta_{\sy}+\delta_{\sigma\sy}+\delta_{\sigma^2\sy}}{3},
\]
where $\x$ is the only periodic point in $\per_2(X')\cap [45]$, and $\y$ is the only periodic point
in $\per_3(X')\cap[123]$. In order to compute the coefficients $\alpha'_1$ and $\alpha'_2$, we need
a second renormalization. The second renormalization gives the topological Markov chain 
$X''\subset\{1,2\}^{\zz}$, defined by the digraph

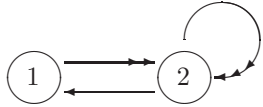
\begin{center}

\unitlength=1truecm
\begin{picture}(5,3)(0,0)
\put(1,1){\circle{0.7}}
\put(0.9,0.9){$1$}
\put(3,1){\circle{0.7}}
\put(2.9,0.9){$2$}

\put(1.4,1.2){\vector(1,0){1.2}}
\put(1.4,1.2){\vector(1,0){1}}
\put(2.6,0.8){\vector(-1,0){1.2}}

\put(3.5,1.5){\oval(1,1)[t]}
\put(3.5,1.5){\oval(1,1)[r]}
\put(3.5,1){\vector(-1,0){0.1}}
\put(3.7,1.05){\vector(-2,-1){0.1}}
\put(3.9,1.2){\vector(-1,-1){0.1}}

\end{picture}

\end{center}
\ms According to the algorithm, we have
\[
\phi''=\left(\begin{matrix} -\infty & -1 \\ -1      &       -2\end{matrix}\right)\hskip 20pt 
\psi''=\left(\begin{matrix} -\infty &  0 \\ \log(5) & \log(3) \end{matrix}\right).\hskip 20pt
\]
The second renormalization gives only a single heavy component,  
$\bX''=\{\z,\sigma\z\}:=\per_2(X'')\cap[12]$, therefore 
\[
\lim_{\beta\to\infty}\mu_{\bp''+\psi''}=\frac{\delta_{\sz}+\delta_{\sigma\sz}}{2},
\]
and we have $\alpha'_1=\alpha'_2=1/2$ for the coefficients in the limit of the first renormalization.
Therefore
\[
\lim_{\beta\to+\infty}\mu_{\bp}=\frac{\delta_{\sa}+\delta_{\sb}+\delta_{\sc}}{6}+\frac{\delta_{\sd}+\delta_{\se}}{4},
\]
for the limit of the original measure. In the table below we present the comparison between the limiting
measure and $\mu_{\bp}$ for different values of the inverse temperature.

\bs
 
\begin{center}
\scriptsize{
\begin{tabular}{|r|c|c|c|c|c|c|c|c|}
\hline
      $\beta=$  &$\log(2)$  &$2\log(2)$&$3\log(2)$&$4\log(2)$&$5\log(2)$&$6\log(2)$& $\cdots$ &$\infty$\\
\hline
$\mu_{\bp}[a]=$ &   0.19273 & 0.18423  & 0.17668  & 0.17200  & 0.16942  & 0.16807  & $\cdots$ &  0.166667\\
\hline
$\mu_{\bp}[b]=$ &  0.18399  & 0.17722  & 0.17395  & 0.17115  & 0.16918  & 0.16800  & $\cdots$ &  0.166677\\
\hline
$\mu_{\bp}[c]=$ &  0.19326  & 0.18343  & 0.17607  & 0.17176  & 0.16935  & 0.16805  & $\cdots$ &  0.166667\\
\hline
$\mu_{\bp}[d]=$ &   0.21312 & 0.22565  & 0.23582  & 0.24227  & 0.24595  & 0.24792  & $\cdots$ &  0.250000\\
\hline
$\mu_{\bp}[e]=$ &  0.21690  & 0.22946  & 0.23748  & 0.24282  & 0.24610  & 0.24796  & $\cdots$ &  0.250000\\
\hline
\end{tabular}
}
\end{center}

\ms \section{Proof of the Renormalization Lemma} \label{thebiglemma} \

\ms We assume that $\phi$ and $\psi$ are normalized as described in Subsection \ref{normalisationpot}.

\ms Let $E_{\phi}=\left\{\phi(C)/|C|:\ C\ \text{is an elementary circuit in } \GG_{X}\right\}$.
Clearly $E_{\phi}$ is finite and $\max E_{\phi}=\bphi=0$. 
We let
\begin{equation}\label{phig}
\phi_g:=\max(E_{\phi}\setminus\{\bphi\}) < 0
\end{equation}
denote the second largest value in $E_{\phi}$.
\ms \paragraph{{\bf \em First Step} (Factorization on Heavy Components)}
For each $1\leq J\leq N_\phi$, let $I_J:=\{\b\in X:\ \b_0\in \bAA_J\}$
be the set of all points whose orbit visits the heavy component $\bX_J$ during an interval 
of time containing the origin.
According to Lemma~\ref{factorization-Lemma} (proved in Section \ref{lessouslemmes}), 
we have for $\beta$ large enough
\begin{equation}\label{factorization-formula}
\mu_{\bpp}[\b_0^n]=\nu_{\psi,J}[\b_0^n]\thinspace \mu_{\bpp} \left(I_J \right)\pm 2e^{\beta\frac{\phi_g}{4}},
\end{equation}
whenever $[\b_0^n]\cap \bX_J\neq \emptyset$ for some $1\leq J\leq N_\phi$.

\ms \paragraph{{\bf \em Second Step} (Excursion Markov Chain)} 
For each $\a\in X$ 
such that $\a_0\in\bigsqcup_{J=1}^{N_\phi}\bAA_J\subset\bAA$, let
$$
i(\a):= \max\{i < 0:\ (\a_{i-1},\a_i)\notin \bEE_{\phi}\}, \; o(\a):=\min\{i > 0:\ (\a_i,\a_{i+1})\notin \bEE_{\phi}\},
$$
$$
i'(\a):=\min\{i > j(\a):\ (\a_{i-1},\a_i)\in\bEE_{\phi}\},\; o'(\a):=\min\{i > i'(\a):\ (\a_i,\a_{i+1})\notin\bEE_{\phi}\},
$$
where $\bEE_{\phi}$ is defined in \eqref{Ebarphi}.
Indices $i(\a), o(\a),i'(\a)$ and $o'(\a)$ are the first and second input/output times to/from heavy
components of the orbit of $\a\in X$.

\ms For $1\leq J\leq N_\phi$ and $a,a'\in \bAA_J$, let
$[a,a']_J:=\{\a\in X:\ \a_{i(\sa)}=a \in \bAA_J,\, \a_{o(\sa)}=a'\}.$
This is the set of all points whose orbit enters the component $\bX_J$ at vertex $a$ and leaves it
at vertex $a'$. Similarly, for $1\leq J,K\leq N_\phi$, $a,a'\in \bX_J$ and $c,c'\in\bX_K$, let
\[
\left[[a,a']_J,[c,c']_K\right]:=
\]
\[
\{\a\in X:\ \a_{i(\sa)}=a\in \bAA_J,\,\a_{o(\sa)}=a'\in \bAA_J,\,
\c_{i'(\sa)}=c\in \bAA_K,\text{ and }\a_{o'(\sa)}=c'\in \bAA_K\},
\]
which is the set of all points entering, under the shift action, the component $\bX_J$ at vertex $a$
and going out at $a'$, and such that the next heavy component they visit is $\bX_K$,
entering at $c$ and going out at $c'$.  
We clearly have
\[
\mu_{\bpp}(I_J)=\sum_{\sd,\se\in \bAA_J}\mu_{\bpp}([a,a']_J).
\]

\ms To estimate $\mu_{\bpp} \left(I_J \right)$, we use the Markov chain on the extended alphabet
\begin{equation}\label{extalphabet}
\AA_{\rm ext}:=\{[a,a']_J:\ 1\leq J\leq N_\phi,\text{ and } a,a'\in \bAA_J\},
\end{equation}
with transition matrix $M_{\bpp}:\AA_{\rm ext}\times\AA_{\rm ext}\to[0,1]$ such that
\[
M_{\bpp}([a,a']_J,[c,c']_{K}):=\frac{\mu_{\bpp}([[a,a']_J,[c,c']_{K}])}{\mu_{\bpp}([a,a']_J)}\cdot
\]
From the shift--invariance of $\mu_{\bpp}$ it follows that $M_{\bpp}$ is a stochastic matrix. Furthermore,
since $\mu_{\bpp}$ is ergodic, then $M_{\bpp}([a,a']_J,[c,c']_{K})$ is irreducible and has a unique invariant 
distribution $\eta_{\bpp}$, which by construction satisfies
\begin{equation}\label{numu}
\eta_{\bpp}([a,a']_J):=
\frac{\mu_{\bpp}([a,a']_J)}{\sum_{K=1}^{N_\phi}\sum_{c,c'\in \bAA_K}\mu_{\bpp}([c,c']_K)}\cdot
\end{equation}

\ms \paragraph{{\bf \em Third Step} (Concentration on Heavy Components)} \

\ms Lemma~\ref{concentration-Lemma} states that 
the measure $\mu_{\bpp}$ concentrates on heavy the components: for $\beta$
large enough
\begin{equation}\label{pouic}
\mu_{\bpp}\left(\bigcup_{K=1}^{N_\phi} I_K\right)\geq 1-e^{\beta\frac{\phi_g}{4}}.
\end{equation}
We defer its (lengthy) proof to Appendix~\ref{concentration}.
Then it follows that the factor $\mu_{\bpp} \left(I_J \right)$ in \eqref{factorization-formula} 
can be approximated by the invariant distribution $\eta_{\bpp}$ of the stochastic matrix $M_{\bpp}$. 
Using the fact that 
$\mu_{\bpp}(I_J)=\sum_{a,a'\in \bAA_J}\mu_{\bpp}([a,a']_J)$ and \eqref{numu}
we get
\[
\mu_{\bpp}(I_J)=\sum_{K=1}^{N_\phi}\mu_{\bpp}(I_K)\thinspace
\thinspace\sum_{a,a'\in \bAA_J}\eta_{\bpp}([a,a']_J)
\]
from which it follows by using \eqref{pouic} that for $\beta$ large enough
\begin{equation}\label{concentration-approximation}
\mu_{\bpp}[\b_0^n]=\nu_{\psi,J}[\b_0^n]\thinspace\sum_{a,a'\in \bAA_J} \eta_{\bpp}([a,a']_J)
\thinspace \pm 3e^{\beta\frac{\phi_g}{4}},
\end{equation}
whenever $[\b_0^n]\cap \bX_J\neq \emptyset$ for some $1\leq J\leq N_\phi$.
Hence, the convergence of $\mu_{\bpp}$ when 
$\beta\to+\infty$, is controlled by the behavior of the invariant distribution $\eta_{\bpp}$, which we 
investigate now. 
\ms \paragraph{{\bf \em Fourth Step} (Excursion potentials)} \

We will replace the stochastic matrix $M_{\bpp}$ 
by a transition matrix $M_{\btp+\tpsi}$, defined by two--symbol potentials 
$\tilde{\phi},\tpsi:\AA_{\rm ext}\to \rr$.
These {\it excursion potentials} are such that the 
one--marginal $\mu_{\btp+\tpsi}^{(1)}$ of the Gibbs measure $\mu_{\btp+\tpsi}$ approaches the invariant 
distribution $\eta_{\bpp}$ of the stochastic matrix $M_{\bpp}$ as $\beta\to+\infty$. 

\ms The excursion potentials are defined as follows.
Let $\v_{\bpp}$ be the right maximal eigenvector 
of $\MM_{\bpp}$, and $\tMM_{\bpp} < \MM_{\bpp}$ be the submatrix of $\MM_{\bpp}$ 
obtained by excluding all the heavy components. 
For each heavy component $1\leq J\leq N_\phi$ and $a,a'\in \bAA_J$, let 
$\tv_{\bpp}(a):=\sum_{b\in\AA}\tMM_{\bpp}(a,b)\v_{\bpp}(b)$.
Then define 
$\tphi,\tpsi:\AA_{\rm ext}\to \rr$ such that
\begin{equation}
\left\{
\begin{array}{l}
\tphi([a,a']_J,[c,c']_{K}) :=\lfphi(a')+\rtphi(a',c)-\lfphi(c),\\ 
\tpsi([a,a']_J,[c,c']_{K})) :=
\log\left(\sum_{a'\to c\in {\bpath}[a',c]}e^{\psi\left({a'\to c}\right)+P_{\psi}(a'\to c)}\right)  \\
\qquad\qquad\qquad\qquad\quad\quad +\log\left(\v_{\psi,J}(a')\w_{\psi,K}(c)\right),
\end{array}\right.
\end{equation}
where ${\bpath}[a',c]\subset \path[a',c]$ is the set of all elementary paths from $a'$ to $c$ 
maximizing $\phi$.
Here we have used the transition pressure, central term, and transition term
as defined in \eqref{transition-pressure},\eqref{central-term} and \eqref{transition-term} respectively. 

\ms To the excursion potential we associate a transition matrix
$M_{\btp+\tpsi}:\AA_{\rm ext}\times\AA_{\rm ext}\to\rr^{+}$ such that 
\[
M_{\btp+\tpsi}\big([a,a']_J,[c,c']_{K}\big):=\exp\left((\btp+\tpsi)([a,a']_J,[c,c']_{K})\right).
\]
The matrices $M_{\btp+\tpsi}$ and $M_{\bpp}$ can be related using 
Lemma~\ref{approximated-excursion-lemma}: for $\beta$ large enough
we have
$$
M_{\bpp}\big([a,a']_J,[c,c']_{K}\big)=
$$
\begin{equation}\label{closeness}
M_{\btp+\tpsi}\big([a,a']_J,[c,c']_{K}\big)\thinspace
         \frac{e^{\beta\lfphi(c')}\v_{\psi,K}(c')\tv_{\bpp}(c')}{e^{\beta\lfphi(a')}\v_{\psi,J}(a')\tv_{\bpp}(a')}
         \thinspace \frac{e^{\pm e^{-\beta\,\delta}}}{\rho_{\bpp}-1},
\end{equation}
for all $1\leq J,K\leq N_\phi$, $a,a'\in \bAA_J$, $c,c'\in \bAA_K$.
Note that $\delta_{\phi} > \phi_g/6$. (By Proposition \ref{propspectralradius}, $\rho_{\bpp}-1\neq 0$.)

\ms The closeness between the one--marginal $\mu_{\btp+\tpsi}^{(1)}$ of the Gibbs measure $\mu_{\btp+\tpsi}$ 
and the invariant distribution $\eta_{\bpp}$ of the stochastic matrix $M_{\bpp}$ follows from 
\eqref{closeness} after the following considerations. 

\ms First notice that the matrix $\bN:\AA_{\rm ext}\times \AA_{\rm ext}\to \rr^+$ 
defined by
\[
\bN([a,a']_J,[c,c']_{K}):= M_{\btp+\tpsi}([a,a']_J,[c,c']_{K})\thinspace
\frac{e^{\beta\lfphi(c')}\v_{\psi,K}(c')\tv_{\bpp}(c')}{e^{\beta\lfphi(a')}\v_{\psi,J}(a')\tv_{\bpp}(a')}
\thinspace \frac{1}{\rho_{\bpp}-1}
\]
is precisely the transition matrix associated to the potential $\Phi:=(\btp+\tpsi)+(h-\sigma\circ h)-\log(\rho_{\bpp}-1)$, 
where $h:\AA_{\rm ext}^\zz\to\rr$ is given by
\[
h([a_0,a_0']_{J_0}[a_1,a_1']_{J_1}\cdots)=-\beta\lfphi(a_0')-\log(\v_{\psi,J_0})(a_0')-\log(\tv_{\bpp})(a_0').
\]
The potentials $\btp+\tpsi$ and $(\btp+\tpsi)+(h-\sigma\circ h)$ are cohomologous, so they define exactly the 
same Gibbs state (see~\cite{Bowen} for details). Furthermore, the potentials $(\btp+\tpsi)+(h-\sigma\circ h)$ and
$\Phi$ differ only by a constant term, therefore they define the same Gibbs state as well.

\ms The one--marginal of the Gibbs measure $\mu_{\Phi}$, which coincides with the one--marginal of $\mu_{\btp+\tpsi}$, 
is completely determined by the maximal eigensystem of the matrix $\bN$. Indeed, the analogous of 
\eqref{equilibrium-state-mubpp} holds, and we have
\[
\nu_{\Phi}([a,a']_J)=\v_{\bN}([a,a']_J)\w_{\bN}([a,a']_J),
\]
for each $[a,a']_J\in \AA_{\rm ext}$. Here $\v_{\bN}$ and $\w_{\bN}$ are respectively the right and left eigenvectors
of $\bN$, associated to the maximal eigenvalue $\rho_{\bN}:=\max(\spec(\bN))$, and normalized such that 
$\w_{\bN}^{\dag}\v_{\bN}=1$. 

\ms Now, we claim that 
$\eta_{\bpp}=\nu_{\Phi}e^{\pm 8(\#\AA_{\rm ext}-1) e^{-\beta\,\delta}}$. By the above cohomological
arguments we also have
$\eta_{\bpp}=\mu_{\btp+\tpsi}^{(1)}\ e^{\pm 8(\#\AA_{\rm ext}-1) e^{-\beta\,\delta}}$.
The claim follows from the application of Proposition~\ref{proposition-stability}
to $M_{\bpp}$ and $\bN$ which are proved to be projectively close
by Lemma~\ref{approximated-excursion-lemma}.
Therefore, taking into account \eqref{concentration-approximation}, it follows that
for $\beta$ large enough
\begin{equation}\label{excursion-approximation}
\mu_{\bpp}[\b_0^n]=\nu_{\psi,J}[\b_0^n]\thinspace\sum_{a,a'\in \bAA_J} \mu_{\btp+\tpsi}^{(1)}([a,a']_J)
 \thinspace\pm 9(\#\AA_{\rm ext}-1)\, e^{-\beta \,\delta},
\end{equation}
whenever $[\b_0^n]\cap \bX_J\neq \emptyset$ for some $1\leq J\leq N_\phi$.
\ms \paragraph{{\bf \em Last Step} (The Projection and the Renormalized Potentials)} \

Let us now simplify the expression for $\mu_{\bpp}[\b_0^n]$ we just proved by first making a
dimensional reductionn (projection of the alphabet), followed by a simplification of the resulting potential,
which will allow us to define the renormalized system.

\ms Let $\AA'=\{1,2,\ldots,N_\phi\}$ be set of indices of heavy components. Define the projection 
$\pi:\AA_{\rm ext}\to \AA'$ such that $\pi([a,a']_J)=J$ for all $1\leq J\leq N_\phi$ and $a,a'\in \bAA_J$,
and extend it coordinatewise to $(\AA_{\rm ext})^{\zz}$. Let $\bmu:=\mu_{\btp+\tpsi}\circ\pi^{-1}$ 
denote the pull back of the measure $\mu_{\btp+\tpsi}$ under the projection $\pi$.
Since
\[
\bmu[J]:=\sum_{a,a'\in\bAA_J} \mu_{\btp+\tpsi}([a,a']_J)\equiv
\sum_{a,a'\in \bAA_J} \mu_{\btp+\tpsi}^{(1)}([a,a']_J), 
\]
for each $J\in\AA'$ and $a,a'\in\bAA_J$, then \eqref{excursion-approximation} can be writen,
for $\beta$ large enough, as
\[
\mu_{\bpp}[\b_0^n]=\nu_{\psi,J}[\b_0^n]\thinspace \bmu[J]\pm 9(\#\AA_{\rm ext}-1)\, 
e^{-\beta\,\delta},
\]
whenever $[\b_0^n]\cap \bX_J\neq \emptyset$ for some $1\leq J\leq N_\phi$.
In Lemma~\ref{projection-Lemma} we prove that  
$\bmu$ equals the Gibbs measure (more properly called Parry measure) $\mu_{\bar{\Phi}}$ 
defined by the $2$-symbol potential  $\bar{\Phi}:\AA_{\rm ext}\times\AA_{\rm ext}\to\rr$ 
such that
\[
\bar{\Phi}(J,K)=\log\thinspace \sum_{a'\in \bAA_J,\,c\in \bAA_K} \exp\big((\btp+\tpsi)([a,a']_J,[c,c']_K)\big).
\]

\ms Hence, in order to compute $\bmu_{\btp+\tpsi}(J)\equiv\mu_{\bar{\Phi}}(J)$, 
we only have to find the left and right positive eigenvectors, $\w_{\bar{\Phi}}$ and 
$\v_{\bar{\Phi}}$, associated to the maximal eigenvalue $\rho_{\bar{\Phi}}$ of the transition
matrix $\MM_{\bar{\Phi}}:\AA'\times\AA'\to\rr^+$ given by
\[
\MM_{\bar{\Phi}}(J,K):=\sum_{a'\in \bAA_J,c\in \bAA_K}\exp \big(\btp+\tpsi)([a,a']_J,[c,c']_K\big).
\]

\ms Instead of computing the left and right positive eigenvectors of $\bMM_{\btp+\tpsi}$, let us first 
approximate this matrix by a more convenient one.
Let us recall the definition of the renormalized potentials
(cf. Definition \ref{newpotentials}).
For $1\leq J,K\leq N_\phi$, let  
\begin{eqnarray*}
\phi'(J,K)&:=&\max\left\{\tphi([a,a']_J,[c,c']_K):\,a'\in\bAA_J,\,c\in\bAA_K\right\},\\
\psi'(J,K)&:=&\log\thinspace \sum_{(a',c)\in\bAA_{J,K}} e^{\tpsi([a,a']_J,[c,c']_K)},\nonumber
\end{eqnarray*}
where $\bAA_{J,K}:=\left\{(a',c)\in\bAA_J\times\bAA_K:\ \tphi([a,a']_J,[c,c']_K)=\phi'(J,K)\right\}$. With this 
renormalized potentials, a rather direct computation allows us to write, for $\beta$ large enough,
\[
\MM_{\bar{\Phi}}=\MM_{\beta\phi'+\psi'} \exp\left(\pm e^{-\beta\tdelta}\right)
\]
where
\[
\tdelta:=\frac{1}{2}\max_{J,K\in\AA'}
\left(\phi'(J,K)-\max\left\{\tphi([a,a']_J,[c,c']_K)\ : \ (a',c)\in \bAA_J\times \bAA_K\setminus \bAA_{J,K}\right\}\right).
\] 
Here $\MM_{\beta\phi'+\psi'}$ is defined in the same way as $\MM_{\bpp}$, using the renormalized 
potentials just defined. 
Now, according to Proposition~\ref{proposition-stability}, 
the invariant distribution $\eta_{\beta\phi'+\psi'}$ is close to the one-marginal of
$\bmu_{\bar{\Phi}}$.
We have 
$$
\nu_{\bar{\Phi}}=\eta_{\beta\phi'+\psi'}\exp\left(\pm 8(\#\AA'-1) \ e^{-\beta\,\tdelta}\right),
$$
from which it follows that for $\beta$ large enough
\begin{eqnarray*}
\mu_{\bpp}[\b_0^n]
&=&
\nu_{\psi,J}[\b_0^n]\thinspace \eta_{\beta\phi'+\psi'}([J])\thinspace \pm 18(\#\AA_{\rm ext}-1)\, e^{-\beta\,\tdelta}\\
&=&\nu_{\psi,J}[\b_0^n]\thinspace \mu_{\beta\phi'+\psi'}([J])\thinspace \pm 18(\#\AA_{\rm ext}-1)\, e^{-\beta\,\tdelta},
\end{eqnarray*}
whenever $[\b_0^n]\cap \bX_J\neq \emptyset$ for some $1\leq J\leq N_\phi$.
The Renormalization Lemma follows by taking 
$\delta':=\min(\frac{\tdelta}{2},\delta)$ and for large enough $\beta$.

\ms \section{Auxiliary Lemmas}\label{lessouslemmes} \

\ms We devote this section to statements and proofs of the auxiliary lemmas used in 
the proof of the Renormalization Lemma (Section \ref{thebiglemma}). We start with the more technical one, 
Lemma~\ref{factorization-Lemma}, which we prove by using periodic approximations of the
Gibbs measure $\mu_{\bpp}$. We are able to give precise estimates of the speed of convergence 
of these periodic approximations, based on a refined version of the Perron--Frobenius Theorem
which we present in Appendix~\ref{proof-of-periodic}.

\ms \subsection{Periodic Approximations and Factorization on Heavy Components}\

\ms \subsubsection{Transition Matrices and Periodic Approximations}\ \label{subsection-transition-matrix}

\ms The matrix $\MM_{\bpp}$ defined in \eqref{matrix} is irreducible and periodic.
If we let $\barp$ be its period, this means that there exists a partition $\{\AA_i, 0\leq i\leq \barp\}$ of $\AA$
($\AA=\bigsqcup_{i=0}^{\barp-1}\AA_i$), and nonnegative rectangular matrices $\QQ_{\bpp,i}:\AA_{i+1}\times \AA_{i}\to\rr^+$, where indices are taken mod $\barp$, such that 
\[
\MM_{\bpp}=\left(\begin{matrix}
                0     & \QQ_{\bpp,1}    &          0     &  \cdots   &            0            \\
                0     &         0       &\QQ_{\bpp,2}    &  \cdots   &            0            \\
                0     &         0       &          0     &  \cdots   &            0            \\
              \vdots  &    \vdots       &       \vdots   &  \ddots   &       \vdots            \\
                0     &         0       &          0     &  \cdots   &   \QQ_{\bpp,\barp-1}    \\
         \QQ_{\bpp,0} &         0       &          0     &  \cdots   &             0
     \end{matrix}\right).
\]
\bs
We also have
\[
\MM_{\bpp}^{\barp}=\left(\begin{matrix}
              \MM_{\bpp,0}  &         0      &   0            &  \cdots  &      0              \\
                       0    &\MM_{\bpp,1}    &   0            &  \cdots  &      0              \\
                       0    &        0       &  \MM_{\bpp,2}  &  \cdots  &      0              \\
                \vdots      &      \vdots    &  \vdots        &  \ddots  & \vdots              \\
                      0     &         0      &   0            &  \cdots  & \MM_{\bpp,\barp-1}
     \end{matrix}\right)
\]

\bs

\noindent
where $\MM_{\bpp,i}:\AA_{i}\times \AA_{i}\to \rr^+$ is primitive for each $0\leq i < \barp$. 
Therefore, according to Perron--Frobenius Theorem (Theorem~\ref{perron-frobenius} in 
Appendix~\ref{proof-of-periodic}) there are unique left and right maximal eigenvectors $\v_{\bpp,i}>0$ 
and $\w_{\bpp,i}>0$ associated to the maximal eigenvalue 
$\max|\spec(\MM_{\bpp,i})|=\rho_{\bpp}^{\barp}$, and normalized such that $\w_{\bpp,i}^{\dag}\v_{\bpp,i}=1$. 
With this we define $\v_{\bpp}:=1/\sqrt{\barp}\otimes_{i=0}^{\barp-1}\v_{\bpp,i}$, and similarly for 
$\w_{\bpp}$. The vectors $\v_{\bpp}$ and $\w_{\bpp}$ so defined are the unique 
left and right eigenvectors associated to the maximal eigenvalue $\rho_{\bpp}:=\max|\spec(\MM_{\bpp})|$, 
normalized such that $\w_{\bpp}^{\dag}\v_{\bpp}=1$. 

\ms For each $\b_0^n\in \AA^{n+1}$ and $p=k\barp >n $ we define the period--$p$ approximation of
$\mu_{\bpp}[\b_0^n]$ by
\begin{eqnarray}\label{periodic-measure}
\PP^{(p)}_{\bpp}[\b_0^n]& :=  &\frac{\sum_{\sa\in\sper_p(X)\cap [\sb_0^n]}e^{S_p(\bpp)(\sa)} 
                                        }{\sum_{\sa\in\sper_p(X)} e^{\beta S_p(\bpp)(\sa)}} \nonumber\\ 
                                      &\equiv&\frac{\prod_{i=0}^{n-1}\MM_{\bpp}(\b_i,\b_{i+1})
                                  \MM_{\bpp}^{p-n}(\b_n,\b_0) }{\trace\left( \MM_{\bpp}^p \right)},
\end{eqnarray}
with $S_p(\bpp)(\a):=\beta\sum_{i=0}^{p-1}\phi(\a_i,\a_{i+1})+\sum_{i=0}^{p-1}\psi(\a_i,\a_{i+1})$, 
$\a\in X$. 

\ms The following result provides an estimate of the convergence rate of $\PP^{(p)}_{\bpp}$ towards 
$\mu_{\bpp}$ when $p\to\infty$, as a function of $\beta$ and the potential $\phi$. 
The proof is deferred to Appendix~\ref{proof-of-periodic}. 
\ms
\begin{proposition}[Periodic Approximation]~\label{periodic-proposition}
Let $\mu_{\bpp}$ be the unique equilibrium state of $\bpp$ defined in \eqref{equilibrium-state-mubpp}
and let $\ell$ be an upper bound for the primitivity indices of the matrices $\MM_{\bpp,i}$, $0\leq i\leq \barp$.
Then, for $\gamma > s_{\phi}:=2\barp\ell \|\phi\|_{\infty}$ and $\beta$ large enough, we have
\[
\mu_{\bpp}[\b_0^n]=\PP^{(k\barp)}_{\bpp}[\b_0^n]\exp\left(\pm e^{-\beta(\gamma-s_{\phi})}\right),
\] 
for all $[\b_0^n]$ such that $[\b_0^n]\cap X\neq \emptyset$.
such that $k\barp > e^{\beta\gamma}+n$. 
\end{proposition}

\ms Let us recall that for each $1\leq J\leq N$, $\bAA_J\subset \AA$ denote the vertex set of the digraph 
$\GG_J:=(\bAA_J,\bEE_J)$ associated to the transitive component $\bX_J\subset \bX$. Let us define the matrices 
$\MM_{\bpp,J}:\bAA_J\times \bAA_J\to\rr^+$ by
\[
\MM_{\bpp,J}(a,')
= \left\{\begin{array}{ll}
e^{(\bpp)(a,a')} & \text{ if } (a,a')\in \bEE_J,\\
0 & \text{ otherwise.} 
\end{array}\right. 
\]

\ms Now, for each $1\leq J\leq N$ and $a,a'\in \bAA_J$, let $a\lds a'$ be any path in the digraph 
$\GG_J$, going from $a$ to $a'$.
Since any circuit $C$ in $\GG_J$ is such that $\phi(C):=|C|\bphi=0$, 
the sum $\phi\left({a\lds a'}\right)$ does not depend on the chosen path $a\lds a'$. 
Therefore for each $1\leq J\leq N$ and $a,a'\in \bAA_J$, we can define {\it the $(a,a')$--compensation term}
by
\begin{equation}\label{compensation-term}
\lfphi(a,a'):=\phi\left({a\lds a'}\right),\, \text{ with }  a\lds a' \text{a path in } 
\GG_J \text{ from } a' \text{ to } a. 
\end{equation}
which relates products of the matrix $\MM_{\psi,J}$ to products of the matrix $\MM_{\bpp,J}$ defined in 
\eqref{psiJ-matrix}. Indeed, since we have fixed $\bphi=0$, it readily follows that
\begin{equation}\label{compensation-effect}
\prod_{i=0}^n\MM_{\bpp,J}^k(\a_i,\a_{i+1})=e^{\beta \lfphi(\sa_0,\sa_n)}\,\prod_{i=0}^n
\MM_{\psi,J}^k(\a_i,\a_{i+1}),
\end{equation}
for each $1\leq J\leq N$ and $\a\in X$ such that $[\a_0^n]\cap\bX_J\neq \emptyset$. 

\ms Now, for each $1\leq J\leq N$ the transition matrix $\MM_{\psi,J}$ is irreducible, therefore there
exists a partition $\bAA_J:=\bigsqcup_{i=0}^{p_J-1}\bAA_{J,i}$, such that
\[
\MM_{\psi,J}=\left(\begin{matrix}
                0       & \QQ_{J,1}    &   0     &\cdots& 0\\
                0       &   0          &\QQ_{J,2}&\cdots& 0\\
                0       &   0          &   0     &\cdots& 0\\
              \vdots    & \vdots       &\vdots   &\ddots& 0\\
                0       &   0          &   0     &\cdots&\QQ_{J,p_J-1}\\
               \QQ_{J,0}&   0          &   0     &\cdots& 0
     \end{matrix}\right).
\]
The rectangular matrices $\QQ_{J,i}:\bAA_{J,i+1}\times \bAA_{J,i}\to\rr^+$, (indices are taken mod $p_J$), are 
non--negative, and such that
\[
\MM_{\psi,J}^{p_J}=\left(\begin{matrix}
              \MM_{J,0}   &   0             &   0        & \cdots &    0        \\
                0         & \MM_{J,1}       &   0        & \cdots &    0        \\
                0         &   0             & \MM_{J,2}  & \cdots &    0        \\
            \vdots        &  \vdots         & \vdots     & \ddots &    0        \\
                0         &    0            &   0        & \cdots & \MM_{J,p_J-1}
     \end{matrix}\right)
\]
with $\MM_{J,i}:\bAA_{J,i}\times \bAA_{J,i}\to \rr^+$ primitive for each $0\leq i < p_J$. 

\ms Once again, Perron--Frobenius Theorem ensures that, for each $0\leq j<p_J$, there are unique left 
and right maximal eigenvectors $\v_{J,i}>0$ and $\w_{J,i}>0$, associated to the maximal eigenvalue 
$\rho_{\psi,J}^{p_J}:=\max|\spec(\MM_{J,i})|$, satisfying $\w_{J,i}^{\dag}\v_{J,i}=1$. 
With this we define $\v_{\psi,J},\w_{\psi,J}:\bAA_J\to\rr^+$ such that $\v_{\psi,J}(\f)=(\sqrt{p_J})\v_{J,i}(a)$ 
for $a\in \bAA_{J,i}$  and similarly for $\w_{\psi,J}$, which are the left and right eigenvectors associated
to $\rho_{\psi,J}:=\max|\spec(\MM_{\psi,J})|$, normalized such that $\w_{\psi,J}^{\dag}\v_{\psi,J}=1$.

\ms Corollary~\ref{periodic-perron-frobenius} in Appendix~\ref{proof-of-periodic} ensures the existence of
constants $C_{\psi}>0$ and $0\leq \tau_{\psi} < 1$ such that
\begin{equation}~\label{periodic-perron-frobenius-J}
\MM_{\psi,J}^{kp_J+r}(a,a')=\frac{\v_{\psi,J}(a)\w_{\psi,J}(a')}{\rho_{\psi,J}^{kp_J+r}}\thinspace
e^{\pm C_{\psi} \tau_{\psi}^{k}}
\end{equation}
for each transitive component $1\leq J\leq N$, and each $a,a'\in \bAA_J$ such that $a\in \bAA_{J,i}$ and 
$a'\in \bAA_{J,(i+r)}$ (indices taken mod $p_J$). 

\ms As mentioned above, to each transitive component $\bX_J\subset \bX$ we can associate the Markovian 
measure $\nu_{\psi,J}$ defined by \eqref{nu-J}, which is precisely the equilibrium 
state on $\bX_J$ associated to the potential $\psi|\bX_J$. Because of the normalization $P(\psi|\bX)=0$, we have
$\rho_{\psi,J}\leq 1$ for each $1\leq J\leq N$, and $\rho_{\psi,J}=1$ for the heavy components 
$1\leq J\leq N_\phi$.

\ms Notice that the transition matrices $\MM_{\psi,J}$ depend only on the potential $\psi$, so that they do
not change with the inverse temperature $\beta$. 

\ms \subsubsection{Approximating the Measure of Incursions}\

\ms The Renormalization Lemma involves the system of incursions into heavy components.
It turns out that the measure $\mu_{\bpp}[\b_0^n]$ of cylinders $[\b_0^n]$ intersecting a heavy
component $\bX_J$ is almost proportional to its $\nu_{\psi,J}$ measure
(Lemma~\ref{factorization-Lemma} below). This is what we call the approximated measure for 
incursions.
 
\ms Fix $q< p=k\barp$, and let
\[
X^{(p,q)}:=\left\{\a\in\per_{p}(X):\ [\a_{p-q}^{q-1}]\cap\bX\neq \emptyset\right\}.
\]
A periodic point $\a\in X^{(p,q)}$ is such that its suffix--prefix factor 
$\a_{p-q}\a_{p-q+1}\cdots\a_0\a_1\cdots\a_{q-1}$ defines a path in $\GG_{X}$ composed by maximizing 
elementary circuits, therefore
\[S_{2q}\phi(\sigma^{-q}(\a))\geq - \#\AA\ \|\phi\|_{\infty},\]
for each $\a \in X^{(p,q)}$. 

\ms Recall that on each heavy component $\bX_J\subset \bX$ we have a 
Markovian measure $\nu_{\psi,J}$  given by
\[
\nu_{\psi,J}[\b_0^n]:=\w_{J}(\b_0)\prod_{i=0}^{n-1}\MM_{\psi,J}(\b_i,\b_{i+1})\v_{\psi,J}(\b_n),
\]
for each $\b\in X$ such that $[\b_0^n]\cap \bX_J\neq\emptyset$. 

\ms Now, for each $1\leq J\leq N_\phi$ and $p=k\barp$ with $k\in \nn$, let  
\[
I_J^{(p)}:=\{\b\in \per_p(X):\ \b_0\in \bAA_J\},
\]
where $\bAA_J$ is the alphabet associate to the heavy component $\bX_J$ as defined above. The set 
$I_J^{(p)}$ can also be obtained as $I_J\cap\per_p(X)$, where 
\[
I_J:=\{\b\in X:\ \b_0\in \bAA_J\},
\]
which is the set of all points whose orbit visits the heavy component $\bX_J$ during an interval 
of time containing the origin. 

\ms We have the following important lemma.

\ms \begin{Lemma}[Factorization on Heavy Components]~\label{factorization-Lemma}
For $\beta$ large enough one has
\[
\mu_{\bpp}[\b_0^n]=\nu_{\psi,J}[\b_0^n]\thinspace \mu_{\bpp} \left(I_J \right)\pm 2 e^{\beta\frac{\phi_g}{4}},
\] 
whenever $[\b_0^n]\cap \bX_J\neq \emptyset$ for some $1\leq J\leq N_\phi$,
where $\phi_g<0$ is defined in \eqref{phig}.
\end{Lemma}

\ms \begin{proof}
Let $q=\lfloor e^{\eta\,\beta}\rfloor$ with  $-\phi_g/4 < \eta < -\phi_g/3$, and
$p=k\barp\geq n+e^{\beta\,\gamma}$ with $\gamma = s_\phi-\phi_g/4$.
In Appendix~\ref{proof-incursion-length}  (Proposition \ref{incursiontime})
we prove that for $\beta$ large enough
\begin{equation}\label{ineq-incursion-length}
\sum_{\sa\in\sper_{p}(X)\cap[\sb_0^n]}e^{S_{p}(\bpp)(\sa)} = \exp 
 \left(\pm 2e^{\beta \frac{\phi_g}{4}}\right) \sum_{\sa\in X^{(p,q)}\cap[\sb_0^n]}
  e^{S_{p}(\bpp)(\sa)},
\end{equation}
for all $\b\in X$ such that $[\b_0^n]\cap \cup_{J=1}^N\bX_J\neq\emptyset$.
Using this and Proposition~\ref{periodic-proposition} we obtain
\begin{eqnarray}\label{first-incursion}
\mu_{\bpp}[\b_0^n]&=&
\frac{\sum_{\sa\in X^{(p,q)}\cap[\sb_0^n]}e^{S_{p}(\bpp)(\sa)}
    }{\sum_{\sa\in \sper_p(X)}e^{S_{p}(\bpp)(\sa)}}
\exp\left(\pm \left(
    2e^{\beta \frac{\phi_g}{4}}+
    e^{-\beta(\gamma-s_{\phi})} \right)\right)\\
                    &=&
\frac{\sum_{\sa\in X^{(p,q)}\cap[\sb_0^n]}e^{S_{p}(\bpp)(\sa)}
    }{\sum_{\sa\in \sper_p(X)}e^{S_{p}(\bpp)(\sa)}}
\exp\left(\pm 3e^{\beta \frac{\phi_g}{4}}\right). \nonumber
\end{eqnarray}

\ms We organize the periodic orbits in $X^{(p,q)}\cap[\b_0^n]$ according to input--output vertices 
as follows. For each $\a\in X^{(p,q)}$ such that $\a_0\in \bAA_J$, let 
\[i(\a):=\max\{0\leq i < p:\ (\a_{i-1},\a_i)\neq \bEE_J\} \text{ and } 
                                  o(\a):=\min\{0\leq i < p:\ (\a_i,\a_{i+1})\neq \bEE_J\}.\] 
Notice that $p-i(\a) > q$ and $o(\a)\geq q-1$. Now, for $a,a'\in \bAA_J$ and $0\leq j < i < p$, 
let
\[
\left[a,i;a',j\right]_J^{(p)}:= \{\a\in X^{(p,q)}:\a_0\in \bAA_J,\, i(\a)=i,\, o(\a)=j,\, 
\a_{i}=a \text{ and } \a_{j}=a'\}
\] 
be the set of all $p$--periodic orbits in $[\b_0^n]$ intersecting the heavy component $\bX_J$, 
leaving this component at vertex $a'$ and time $j$, and entering the last time at vertex $a$ 
and time $i$.

\ms It follows from \eqref{compensation-effect} that for each $a,a'\in \bAA_J$ and $0\leq j < i < p$ 
we have
\begin{eqnarray*}
\lefteqn{\sum_{\sa\in [\sb_0^n]\cap [a,i;a',j]_J^{(p)}}e^{S_{p}(\bpp)(\sa)}}\\
&&=\MM_{\bpp,J}^{p-i}(a,\b_0)\left(\prod_{i=0}^{n-1}\MM_{\bpp,J}(\b_i,\b_{i+1})\right)
           \MM_{\bpp,J}^{j-n}(\b_n,a')\MM_{\bpp}^{i-j}(a',a) \\
 &&=\MM_{\psi,J}^{p-i}(a,\b_0)\left(\prod_{i=0}^{n-1}\MM_{\psi,J}(\b_i,\b_{i+1})\right)
           \MM_{\psi,J}^{j-n}(\b_n,a') e^{\lfphi(a,a')} \MM_{\bpp}^{i-j}(a',a), 
\end{eqnarray*}
where $\lfphi(a,a')$ is the compensation term, defined in~\eqref{compensation-term}, which relates 
$\MM_{\psi,J}$ to $\MM_{\bpp,J}$. Now, since $\min(p-i, j-n, p-i+j)\geq q-n$, then using three 
times \eqref{periodic-perron-frobenius-J}, and the definition \eqref{nu-J} of $\nu_{J,\psi}$, we obtain
\begin{eqnarray*}
\lefteqn{\sum_{\sa\in [\sb_0^n]\cap [a,i;a',j]_J^{(p)}}e^{S_{p}(\bpp)(\sa)}}\\
&& \hskip 30pt =\v_{\psi,J}(a)\w_{\psi,J}(\b_0)\left(\prod_{i=0}^{n-1}\MM_{\psi,J}(\b_i,\b_{i+1})\right)\v_{\psi,J}(\b_n)
       \w_{\psi,J}(a') e^{\lfphi(a,a')}\\ \\
& &\hskip 90pt
     \times \, \MM_{\bpp}^{i-j}(a',a) \,\,e^{\pm 2C_{\psi} \tau_{\psi}^{\lfloor (q-n)/p_J\rfloor}}\\ \\
&& \hskip 30pt =\nu_{\psi,J}[\b_0^n]\v_{J}(a)\w_{\psi,J}(a') \, 
      \MM_{\bpp}^{i-j}(a',a)\, e^{\pm 2C_{\psi} \tau_{\psi}^{\lfloor (q-n)/p_J\rfloor}}\, e^{\lfphi(a,a')} \\ \\
&& \hskip 30pt =\nu_{\psi,J}[\b_0^n]\,\MM_{\psi,J}(a,a')^{p-i+j}\,\MM_{\bpp}^{i-j}(a',a)\,
       e^{\lfphi(a,a')}\, e^{\pm 3C_{\psi} \tau_{\psi}^{\lfloor (q-n)/p_J\rfloor}}\\ \\
&& \hskip 30pt =\nu_{\psi,J}[\b_0^n]\,\MM_{\bpp,J}(a,a')^{p-i+j}\, \MM_{\bpp}^{i-j}(a',a)\,
      e^{\pm 3C_{\psi} \tau_{\psi}^{\lfloor (q-n)/p_J\rfloor}}\\ \\
&& \hskip 30pt =e^{\pm 3C_{\psi} \tau_{\psi}^{\lfloor (q-n)/p_J\rfloor}}\,  \nu_{\psi,J}[\b_0^n]\,
   \sum_{\sa\in [a,i;a',j]_J^{(p)}}e^{S_{p}(\bpp)(\sa)}.
\end{eqnarray*}
With this, \eqref{first-incursion}, and taking into account that
\[
\{\a\in X^{(p,q)}:\ \a_0\in \bAA_J\}=\bigcup_{a,a'\in \bAA_J}\bigcup_{q\leq j\leq i\leq p-q}[a,i;a',j]_J^{(p)},
\]
it follows that
\begin{eqnarray*}
\mu_{\bpp}[\b_0^n]&=&\nu_{\psi,J}[\b_0^n]\thinspace\PP_{\bpp}^{(p)}\left(\{\a\in X^{(p,q)}:\ \a_0\in \bAA_J\}\right)\\
                              &   & \hskip 80pt \times
\exp\left(\pm 3\left(e^{\beta \frac{\phi_g}{4}}+C_{\psi} \tau_{\psi}^{\lfloor (q-n)/p_J\rfloor}\right)\right).
\end{eqnarray*}

\ms Now, \eqref{ineq-incursion-length} 
implies that
\[
\PP_{\bpp}^{(p)}\left(\{\a\in X^{(p,q)}:\ \a_0\in \bAA_J\}\right)=\PP_{\bpp}^{(p)}\left(I_J^{(p)}\right)
\thinspace \exp\left(\pm 2e^{-\beta \frac{\phi_g}{2}}\right),
\] 
therefore
\[
\mu_{\bpp}[\b_0^n]=\nu_{\psi,J}[\b_0^n]\thinspace\PP_{\bpp}^{(p)}\left(I_J^{(p)}\right)
\thinspace  \exp\left(\pm  \left(5e^{\beta \frac{\phi_g}{4}}+3C_{\psi} \tau_{\psi}^{\lfloor (q-n)/p_J\rfloor}\right)\right),
\]
for each $q = \lfloor e^{\eta\,\beta}\rfloor$, with $-\phi_g/4<\eta < -\phi_g/3$.

\ms Since $\tau_\phi < 1$, then $C_{\psi} \tau_{\psi}^{\lfloor (q-n)/p_J\rfloor} 
\leq e^{\beta \phi_g/4}$ if we take $q\geq \beta^2$ and $\beta$ large enough. 
Therefore
\[
\mu_{\bpp}[\b_0^n]=
\nu_{\psi,J}[\b_0^n]\thinspace\PP_{\bpp}^{(p)}\left(I_J^{(p)}\right)
\thinspace \exp\left(\pm 6 e^{\beta \frac{\phi_g}{4}}\right),
\]
and since $\PP_{\bpp}^{(p)}\left(I_J^{(p)}\right)=\PP_{\bpp}^{(p)}\left(I_J\right)$, then 
we obtain
\begin{equation}\label{incursion-equation}
\mu_{\bpp}[\b_0^n]=\nu_{\psi,J}[\b_0^n]\thinspace \PP_{\bpp}^{(p)}
\left(I_J\right)\thinspace
\exp\left(\pm 6e^{\beta\frac{\phi_g}{4}}\right)
\end{equation}
for each $\b\in X$ such that $[\b_0]\cap\bX_J\neq\emptyset$.

\ms Notice now that $X\setminus I_J:=\{\b\in X:\ \b_0\not\in \bAA_J\}$ is a union of cylinder sets, therefore 
Proposition~\ref{periodic-proposition} implies that for $\beta$ large enough
\[
1-\PP_{\bpp}^{(p)}(I_J)=(1-\mu_{\bpp}(I_J))\exp\left(\pm e^{-\beta(\gamma-s_{\phi})}\right).
\]
Hence,
\begin{equation}\label{approximationPP} 
\PP_{\bpp}^{(p)}(I_J)=\mu_{\bpp}(I_J)+(1-\mu_{\bpp}(I_J))\left(1-\exp\left(\pm e^{-\beta(\gamma-s_{\phi})}\right)\right).
\end{equation}
Since
\[
\exp\left(\pm  e^{-\beta(\gamma-s_{\phi})}\right)  = 1 \pm \frac32\, e^{-\beta(\gamma-s_{\phi})}\text{ and } 
\exp\left(\pm 6e^{\beta \frac{\phi_g}{4}}\right)     = 1 \pm 9\, e^{\beta \frac{\phi_g}{4}},
\]
for $\beta$ large enough,  then, substituting in 
\eqref{incursion-equation} and~\eqref{approximationPP} we obtain
\[
\mu_{\bpp}[\b_0^n]= \nu_{\psi,J}[\b_0^n]\thinspace
\mu_{\bpp}(I_J)\pm \frac32\, e^{\beta\frac{\phi_g}{4}}\thinspace\left(1 + 9e^{\beta \frac{\phi_g}{4}}\right),
\]
from which the result follows by taking $\beta$ big enough.
\end{proof}

\ms \subsection{Excursion Potentials}\

\ms We will now replace $M_{\bpp}$ by another transition matrix closed to it, whose entries can be 
explicitly computed in terms of the renormalized protentials defined above. Let us recall
their definitions.

\ms Let  $\w_{\bpp}$ and $\v_{\bpp}$ be the left and right maximal eigenvectors 
of our original transition matrix $\MM_{\bpp}$ on the alphabet $\AA$. Let 
$\tAA:=\AA\setminus \bigsqcup_{J=1}^{N_\phi}\bAA_J$, and $\tMM_{\bpp} < \MM_{\bpp}$ the submatrix of $\MM_{\bpp}$ obtained 
by excluding all the heavy components. 

\ms Now, for each heavy component $1\leq J\leq N_\phi$ and $a,a'\in \bAA_J$ let
\[
\tw_{\bpp}(a):=\sum_{b\in\AA}\w_{\bpp}(b)\tMM_{\bpp}(b,a),
\]
\[
\tv_{\bpp}(a):=\sum_{b\in\AA}\tMM_{\bpp}(a',b)\v_{\bpp}(b).
\]

\ms  Let us recall the definition of the transition term (see \eqref{transition-term} above).
For $a'\in \bAA_J$ and $c\in \bAA_K$, let $\path[a',c]$ be the set of all elementary paths, starting 
at $a'$ and ending at $c$, with no arrows in $\bEE_{\phi}$. The $(a',c)$--transition 
term is the maximum
\[
\rtphi(a',c):=\max\{\phi\left({a'\to c}\right):\ a'\to c\in \path[a',c]\},
\]
where $\phi\left({a'\to c}\right):=\phi(a',b_1)+\sum_{i=1}^{m-1}\phi(b_i,b_{i+1})+\phi(b_{m},c)$ for 
$a\to c=(a',b_1,\ldots,b_{m},c)$. 

\ms Let $\MM_{\psi}:\AA\times \AA\to\rr^+$ be defined in the same way as $\MM_{\bpp}$, with $\psi$ 
replacing $\bpp$, and let $\bMM_{\psi}$ be the restriction of $\MM_{\psi}$ to all the 
transitive components of $\bX$, either heavy or not. Because of the normalization $P(\psi|\bX)=0$, we 
have $\sum_{n=1}^{\infty}\bMM_{\psi}^n(b,b) < \infty$,
for each $b\not\in \tAA$. This is due to the fact that 
non--heavy transitive components of $\bX$ have maximal eigenvalue strictly smaller than 1.

 
\ms Let us now fix, for each heavy component $1\leq J\leq N_\phi$, a {\it central vertex} $c_J\in \bAA_J$. 
Using this vertex, define, for each $a'\in \bAA_J$, the {\it central term}, whose definition we remind here:
\[
\lfphi(a'):=\phi\left({c_J\lds a'}\right) \text{ with } c_J\lds a'\text{ a path in } \GG_J \text{ from } 
c_J \text{ to } a.
\]
It is easy to see that the central term satisfies the relation
\begin{equation}\label{relation-terms}
\lfphi(a')=\lfphi(a)+\lfphi(a,a').
\end{equation}

\ms Recall that the excursion potentials $\tphi,\tpsi:\AA_{\rm ext}\times\AA_{\rm ext}\to\rr$
are defined as 
\begin{eqnarray*}
\tphi([a,a']_J,[c,c']_{K})&:=& \lfphi(a')+\rtphi(a',c)-\lfphi(c),\\ 
\tpsi([a,a']_J,[c,c']_{K}))&:=&
\log\left(\sum_{a'\to c\in \bpath [a',c]}e^{\psi\left({a'\to c}\right)+P_{\psi}(a'\to c)}\right)  \\
&&+\log\left(\v_{\psi,J}(a')\w_{\psi,K}(c)\right),\nonumber
\end{eqnarray*}
where $P_{\psi}(a'\to c)$ is the transition pressure defined in \eqref{transition-pressure} and
where $\bpath[a',c]\subset \path[a',c]$ is the set of all elementary paths from $a'$ to $c$ 
maximizing $\phi$. As usual, we will denote by $M_{\btp+\tpsi}$ the transition matrix defined by the
potential $\btp+\tpsi$. 

\ms Recall the notation $\phi(C):=\sum_{i=0}^{|C|-1}\phi(b_i,b_{i+1})$ for the circuit 
$C=(b_0,\ldots,b_{|C|-1})$. 

\ms\begin{Lemma}[Approximated Cohomology]\label{approximated-excursion-lemma}
For $\beta$ large enough, and for all  $a,a'\in \bAA_J$,  
$c,c'\in \bAA_K$, $1\leq J,K\leq N_\phi$, we have
\[
M_{\bpp}([a,a']_J,[c,c']_{K})=
\]
\[
M_{\btp+\tpsi}([a,a']_J,[c,c']_{K})\thinspace
         \frac{e^{\beta\lfphi(c')}\v_{\psi,K}(c')\tv_{\bpp}(c')}{e^{\beta\lfphi(a')}\v_{\psi,J}(a')\tv_{\bpp}(a')}
         \thinspace \frac{e^{\pm e^{-\beta\,\delta}}}{\rho_{\bpp}-1}
\]
where
\[
\delta:=\frac{1}{6}\min\left\{|\phi(C)-\phi(C')|:\  C,\, C'\in \CC, \, \phi(C)\neq\phi(C'), 
\text{ and }  |C|,|C'| \leq 2\#\AA\right\}.
\]
\end{Lemma}
\ms \begin{proof} First of all notice that
\begin{eqnarray*}
\mu_{\bpp}([[a,a']_J,[c,c']_{K}])&=&\tw_{\bpp}(a)
\left(\sum_{k=0}^\infty\frac{\MM_{\bpp,J}^k(a,a')}{\rho_{\bpp}^{k}}\right) 
                   \left(
                   \sum_{k=0}^\infty\frac{\tMM_{\bpp}^k(a',c)}{\rho_{\bpp}^{k}}\right)\\
                    & &\hskip 100pt\times 
                       \left(\sum_{k=0}^\infty\frac{\MM_{\bpp,K}^k(c,c')}{\rho_{\bpp}^{k}}\right)\tv_{\bpp}(c'),\\
\mu_{\bpp}([a,a']_J) &=&\tw_{\bpp}(a)\left(\sum_{k=0}^\infty\frac{\MM_{\bpp,J}^k(a,a')}{\rho_{\bpp}^{k}}\right)
                           \tv_{\bpp}(a').
\end{eqnarray*}
Hence, taking into account~\eqref{compensation-effect} and~\eqref{relation-terms}
we have
\begin{eqnarray*}
\lefteqn{M_{\bpp}([[a,a']_J,[c,c']_{K}])=}\\
&& \left(\sum_{k=0}^\infty\frac{\tMM_{\bpp}^k(a',c)}{\rho_{\bpp}^{k}}\right)
                               \thinspace
\left(\sum_{k=1}^\infty\frac{\MM_{\psi,K}^k(c,c')}{\rho_{\bpp}^{k}}\right)
       \thinspace\frac{\tv_{\bpp}(c')e^{\beta(\lfphi(c')-\lfphi(c))}}{\tv_{\bpp}(a')}.
\end{eqnarray*}

\ms We now deal with 
\[
Q(a',c):= \sum_{k=0}^\infty\frac{\tMM_{\bpp}^k(a',c)}{\rho_{\bpp}^{k}} , \quad
R(c,c'):= \sum_{k=1}^\infty\frac{\MM_{\psi,K}^k(c,c')}{\rho_{\bpp}^{k}}.
\]

\ms Notice that $Q(a',c)=\sum_{a'\lds c \in \widehat\path[a',c]}
e^{\beta\phi\left({a'\lds c}\right)+\psi\left({a'\lds c}\right)-|a'\lds c|\log(\rho_{\bpp})}$,
where $\widehat\path[a',c]$ denotes the collection of all paths in $\GG_{X}$ going from $a'$ to $c$, with no arrows in 
$\bEE_{\phi}$ (defined in \eqref{Ebarphi}).

\ms Any path $a'\lds c\in\widehat\path[a',c]$ can be decomposed into an elementary path 
$a'\to c=(a,b_1,\ldots,b_{m-1},c)\in \path[a',c]$ and a sum of circuits 
$b_{i_1}\cira + b_{i_2}\cira +\cdots + b_{i_\ell}\cira$, 
with no arrows in $\bEE_{\phi}$.
Let us denote by $\tCC_b$ the collection of all the circuits in $\GG_{X}$ 
with base point $b$, and with no arrows in $\bEE_{\phi}$.
Taking this into account, we can rewrite $Q(a',c)$ as 
\[
Q(a',c)=\sum_{a'\to c\in \path[a',c]}e^{\beta\phi\left({a'\to c}\right)+\psi\left({a'\to c}\right)
-|a'\to c|\log(\rho_{\bpp})}
\thinspace\sum_{b \in a'\to c}\thinspace \sum_{k=0}^\infty\frac{\tMM_{\bpp}^k(b,b)}{\rho^k_{\bpp}},
\] 
where $\path[a',c]$ denotes the collection of all elementary paths in $\GG_{X}$ going from $a'$ to $c$, 
and $b\in a'\to c$ means that the path $a'\to c$ passes through the vertex $b$.
By Proposition \ref{ex} (Appendix~\ref{excursions}) we have
\[
\sum_{k=0}^{\infty}\frac{\tMM^k_{\bpp}(b,b)}{\rho_{\bpp}^k}=
\sum_{k=0}^\infty\bMM_{\psi}^k(b,b)\pm D\thinspace e^{\beta \phi_g},
\]
for $\beta$ large enough, for all $b\not\in\bigsqcup_{J=1}^{N_\phi}\bAA_J:=\bigsqcup_{J=1}^{N_\phi}\bAA_J$,  
and for some constant $D>0$.
Therefore, 
\begin{eqnarray}~\label{convergent-series}
\sum_{b \in a'\to c}\sum_{k=0}^\infty \frac{\tMM_{\bpp}^k(b,b)}{\rho^k_{\bpp}}&=&
e^{P_{\psi}(a'\to b)}\left(1\pm \frac{|a'\to c|D
                                    }{e^{P_{\psi}(a'\to c)}}e^{\beta \phi_g}\right)\\
                                                                               &=&
e^{P_{\psi}(a'\to c)}\exp\left(\pm e^{\beta\, \frac{3\phi_g}{4}}\right),\nonumber
\end{eqnarray}
for all $\beta$ large enough.

\ms Now, if $a'\to c\in \path[a',c]$ is not maximal, then we necessarily have 
\[
\phi\left({a'\to c}\right)\leq 
\rtphi(a',c) - \min\left\{|\phi(C)-\phi(C')|:\ \phi(C)\neq \phi(C'),\, C,C'\in \CC[a',c] \right\},
\] 
where $\CC_0[a',c]$ denote the set of circuits formed by an elementary path from $a'$ to $c$ and followed
by an elementary path from $c$ to $a'$. 
Since
$$
6\delta < \min\left\{|\phi(C)-\phi(C')|:\ \phi(C)\neq \phi(C'),\, C,C'\in \CC[a',c] \right\}
$$ 
then
\begin{eqnarray*}
\lefteqn{\sum_{a'\to c\in \path[a',c]} e^{\beta\phi\left({a'\to c}\right)+\psi\left({a'\to c}\right)+
         P_{\psi}(a'\to c)} }\\
      &=& e^{\beta\rtphi(a',c)}\sum_{a'\to c \in \bpath[a',c]} e^{\psi\left({a'\to c}\right)+
         P_{\psi}(a'\to c)} \times \\
      & &\hskip 20pt \left(1\pm e^{\beta\, 6\delta}
            \frac{\sum_{a'\to c \in \path[a',c]} e^{\psi\left({a'\to c}\right)+ P_{\psi}(a'\to c)}
                }{\sum_{a'\to c\in \bpath[a',c]} e^{\psi\left({a'\to c}\right)+ P_{\psi}(a'\to c)}}\right)\\ \\
      &=&  \exp\left(\pm e^{-\beta\,5\delta}\right) e^{\beta\rtphi(a',c)}
          \sum_{a'\to c \in \bpath[a',c]}  e^{\psi\left({a'\to c}\right)+P_{\psi}(a'\to c)}
\end{eqnarray*}
for $\beta$ large enough.
Proposition \ref{propspectralradius} (Appendix~\ref{spectralstuff}) implies that
$\rho_{\bpp}=1\pm e^{\beta \phi_g/2}$, for $\beta$ large enough.
On the other hand, since $a'\to c\in\path[a',c]$ is  an elementary, then $|a'\to c| \leq  \#\AA$, therefore,
 for $\beta$ large enough, $|a'\to c|\log(\rho_{\bpp})=\pm e^{\beta\frac{2\phi_g}{5}}$.
Taking this into account, and using \eqref{convergent-series}, we obtain
\begin{eqnarray}\label{bound-Q}
Q(a',c)&=&\exp\left(\pm\left(e^{-\beta\,5\delta}+e^{\beta\, 2\phi_g/5}\right)\right)
            e^{\beta\rtphi(a',c)} \sum_{a'\to c\in \bpath[a',c]}e^{\psi\left({a'\to c}\right)+
            P_{\psi}(a'\to c)} \nonumber \\
        &=&\exp\left(\pm e^{-\beta\,2\delta}\right)  e^{\beta\rtphi(a',c)}
           \sum_{a'\to c\in \bpath[a',c]} e^{\psi\left({a'\to c}\right)+P_{\psi}(a'\to c)}
\end{eqnarray}
for $\beta$ large enough. we used the fact that $-2\phi_g/5 \geq  12\delta/5>2\delta$.

\ms The factor $R(c,c'):=\sum_{k=1}^\infty\MM_{\psi,K}^k(c,c')/\rho_{\bpp}^{k}$ can be treated as 
follows. Fix $k_0=k_0(\beta)$ so that $C_\psi\tau_{\psi}^{k_0} \leq e^{\beta \phi_g}$, with $C_\psi$ 
and $\tau_\psi$ as in \eqref{periodic-perron-frobenius-J}.
Since $\tau_{\psi}<1$, then we can choose $k_0$ proportional to $\beta$.
Thus, from the cited lemma we obtain
\[
\sum_{k=k_0}^{\infty}\frac{\MM_{\psi,K}^k(c,c')}{\rho_{\bpp}^{k}}=\exp\left(\pm e^{\beta \phi_g}\right)
\w_{\psi,K}(c)\v_{\psi,K}(c')\thinspace \sum_{k=k_0}^{\infty} \rho_{\bpp}^{-(p_Kk + r)},
\]
where $r=r(c,c')$ is the smallest integer such that $\MM_{\psi,K}^r(c,c') >0$. Now, for $k < k_0$ we have
\[
\sum_{k=0}^{k_0-1}\frac{\MM_{\psi,K}^k(c,c')}{\rho_{\bpp}^{k}} < \sum_{k=0}^{k_0-1}\MM_{\psi,K}^{p_Kk+r}(c,c')
\leq \frac{k_0\,\max_{k\in\nn}\MM_{\psi,K}^k(c,c')}{p_K},
\]
where $\max_{k\in\nn}\MM_{\psi,K}^k(c,c')<\infty$ follows from the fact that $\max|\spec(\MM_{\psi,K})|=1$. 

\ms By Proposition \ref{propspectralradius} one has $\rho_{\bpp} \leq 1+e^{\beta\frac{\phi_g}{2}}$.
From this, after a few computations, it follows that
$\sum_{k=0}^{\infty} \rho_{\bpp}^{-(p_Kk + r)}\geq 2^{-r}e^{-\beta\,\frac{\phi_g}{2}}/p_K$. 
Then, since $k_0$ is proportional to $\beta$, by taking $\beta$ large enough
we obtain
\[
\sum_{k=0}^{k_0-1}\frac{\MM_K^k(c,c')}{\rho_{\bpp}^{k}} < e^{\beta\frac{\phi_g}{3}}
\w_{\psi,K}(c)\v_{\psi,K}(c') \thinspace \sum_{k=0}^{\infty} \rho_{\bpp}^{-(p_Kk + r)},
\]
and from this
\begin{eqnarray*}
\left(e^{e^{\beta\phi_g}}-e^{\beta\frac{\phi_g}{3}}\right)\w_{\psi,K}(c)\v_{\psi,K}(c')
\thinspace \frac{\rho_{\bpp}^{p_K-r}}{\rho_{\bpp}-1} &\leq&\sum_{k=0}^{\infty}\frac{\MM_K^k(c,c')}{\rho_{\bpp}^{k}},\\
\left(e^{e^{\beta\phi_g}}+e^{\beta\frac{\phi_g}{3}}\right) \w_{\psi,K}(c)\v_{\psi,K}(c') 
\thinspace  \frac{\rho_{\bpp}^{p_K-r}}{\rho_{\bpp}-1} &\geq &\sum_{k=0}^{\infty}\frac{\MM_K^k(c,c')}{\rho_{\bpp}^{k}}.
\end{eqnarray*}
From these two inequalities, and taking into account that $\rho_{\bpp}=\exp\left(\pm e^{\beta\phi_g/2}\right)$, 
it follows that
\begin{eqnarray}\label{bound-R}
R(c,c')&=&\exp\left(\pm 3e^{-\beta\frac{4\delta}{3}}\right)\w_{\psi,K}(c)\v_{\psi,K}(c')
           \thinspace  \frac{\rho_{\bpp}^{p_K-r}}{\rho_{\bpp}-1}\\
       &=&\frac{\exp\left(\pm \left(3e^{-\beta\frac{4\delta}{3}}+(p_K-r)e^{\beta\phi_g/2}\right)\right)
              }{\rho_{\bpp}-1} \thinspace \w_{\psi,K}(c)\v_{\psi,K}(c')\nonumber\\
        &=&\frac{\exp\left(\pm 4e^{-\beta\frac{4\delta}{3}}\right)
              }{\rho_{\bpp}-1} \thinspace \w_{\psi,K}(c)\v_{\psi,K}(c')\nonumber
\end{eqnarray}
for $\beta$ large enough. The result follows from the 
bounds~\eqref{bound-Q} and~\eqref{bound-R}.

\end{proof}

\bs \subsection{Projecting the Excursion System}\

\ms The aim here is to ``compress" the excursion system defined on $\AA_{\rm ext}$ 
\eqref {extalphabet}. 

\ms Recall that $\AA'=\{1,2,\ldots,N_\phi\}$ is the index set of the heavy components and that the projection 
$\pi:\AA_{\rm ext}\to \AA'$ is such that $\pi([a,a']_J)=J$ for all $1\leq J\leq N_\phi$ and $a,a'\in \bAA_J$, and 
extend it coordinatewise to $(\AA_{\rm ext})^{\zz}$. Let
\[
\bmu_{\btp+\tpsi}:=\mu_{\btp+\tpsi}\circ\pi^{-1}
\]
denote the pull back of the measure $\mu_{\btp+\tpsi}$ under the projection $\pi$. Since
\[
\bmu_{\btp+\tpsi}[J]=\sum_{a,a'\in\bAA_J} \mu_{\btp+\tpsi}(\left[[a,a']_J\right])\equiv
\sum_{\sd,\se\in \bAA_J} \mu_{\btp+\tpsi}^{(1)}([a,a']_J), 
\]
for each $J\in\AA'$ and $a,a'\in\bAA_J$, then \eqref{excursion-approximation} (fourth step in
the proof of the Renormalization Lemma) can be written
as
\[
\mu_{\bpp}[\b_0^n]=\nu_{\psi,J}[\b_0^n]\thinspace \bmu_{\btp+\tpsi}[J]\pm 9(\#\AA_{\rm ext}-1)\, 
\]
for $\beta$ large enough and whenever $[\b_0^n]\cap \bX_J\neq \emptyset$ for some $1\leq J\leq N_\phi$.

\ms We have the following.

\ms \begin{Lemma}[The Projection is Markovian]~\label{projection-Lemma} 
The pull back $\bmu_{\btp+\tpsi}$ of the measure $\mu_{\btp+\tpsi}$
under the projection $\pi:(\AA_{\rm ext})^{\zz}\to(\AA')^\zz$, coincides with the Gibbs state
defined by the $2$-symbol potential $\btp+\tpsi:\AA_{\rm ext}\times\AA_{\rm ext}\to\rr$. 
\end{Lemma}

\ms \begin{proof}
\ms The Gibbs state (Markov measure) associated to $\btp+\tpsi$, where $\tphi$ and $\tpsi$ are the approximate excursion potentials, is defined by
$$
\mu_{\btp+\tpsi}\left[[a_0,a'_0]_{J_0}\cdots[a_n,a'_n]_{J_n}\right]:=
$$
\begin{equation}
\label{approximated-Parry-measure}
\w_{\btp+\tpsi}([a_0,a'_0]_{J_0})  
\frac{\prod_{i=0}^{n-1}M_{\btp+\tpsi}([a_i,a'_i]_{J_i},[a_{i+1},a'_{i+1}]_{J_{i+1}})}{\rho_{\btp+\tpsi}^{n}}
\v_{\btp+\tpsi}([a_n,a'_n]_{J_n}),
\end{equation}
with $\rho_{\btp+\tpsi}$ the maximal eigenvalue of associated transition matrix $M_{\btp+\tpsi}$. 

\ms Since the matrix element $M_{\btp+\tpsi}([a,a']_J,[c,c']_K)$ depends only on the internal
symbols $a'$ and $c$, we can write $M_{\btp+\tpsi}(a'_J,c_K)$ instead of $M_{\btp+\tpsi}([a,a']_J,[c,c']_K)$. 
Now, for each $1\leq J\leq N_\phi$ and $a'\in \bAA_J$ fixed, let 
$\hat{\w}_{\btp+\psi}(a_J):=\sum_{a\in \bAA_J}\w_{\btp+\tpsi}([a,a']_J)$. Then, since $\w_{\btp+\tpsi}$ 
is the left invariant vector associated to the maximal eigenvalue of $M_{\btp+\tpsi}$, we have
\begin{eqnarray*}
\w_{\btp+\tpsi}([a,a']_K)&=&\frac{1}{\rho_{\btp+\tpsi}}\sum_{J;\, a,a'\in\bAA_J}\w_{\btp+\tpsi}
                                                      ([a,a']_J)M_{\btp+\tpsi}([a,a']_J,[c,c']_K)\\ 
                          &=&\frac{1}{\rho_{\btp+\tpsi}}\sum_{J;\, a,a'\in\bAA_J}\w_{\btp+\tpsi}
                                                      ([a,a']_J)M_{\btp+\tpsi}(a'_J,c_K)\\
                          &=&\frac{1}{\rho_{\btp+\tpsi}}\sum_{J,\, \se\in\bAA_J}\hat{\w}_{\btp+\tpsi}(a'_J)
                                                                          M_{\btp+\tpsi}(a'_J,c_K),
\end{eqnarray*}
for any $1\leq K\leq N_\phi$ and $c,c'\in \bAA_K$. In this way we show that $\w_{\btp+\tpsi}([c,c']_K)$ does 
not depend on $c'$. A similar computation shows that $\v_{\btp+\tpsi}([c,c']_K)$ does not depend on $c$,
and we can write $\w_{\btp+\tpsi}(c_K)$ instead of $\w_{\btp+\tpsi}([c,c']_K)$ and 
$\v_{\btp+\tpsi}(c'_K)$ instead of $\v_{\btp+\tpsi}([c,c']_K)$. 

\ms For each $J\in\AA'$ let $\bw_{\btp+\tpsi}(J):=\sum_{a\in\bAA_J}\w_{\btp+\tpsi}(a_J)$ and 
$\bv_{\btp+\tpsi}(J):=\sum_{a'\in\bAA_J}\v_{\btp+\tpsi}(a'_J)$. Then we have
\begin{eqnarray*}
\bw_{\btp+\tpsi}(K)&=&\frac{1}{\rho_{\btp+\tpsi}}
   \sum_{J;\,a,a'\in\bAA_J;\,c\in\bAA_K}\w_{\btp+\tpsi}([a,a']_J) M_{\btp+\tpsi}([a,a']_J,[c,c']_K)\\
                       &=&\frac{1}{\rho_{\btp+\tpsi}}
   \sum_{J;\,a,a'\in\bAA_J;\,c\in\bAA_K}\w_{\btp+\tpsi}(a_J) M_{\btp+\tpsi}(a'_J,c_K)\\
                       &=&\frac{1}{\rho_{\btp+\tpsi}}
   \sum_{J;\,a,a'\in\bAA_J;\,c\in \bAA_K}\bar{\w}_{\btp+\tpsi}(J)\bMM_{\btp+\tpsi}(J,K),
\end{eqnarray*}
from which it follows $\rho_{\btp+\tpsi}\in\spec(\bMM_{\btp+\tpsi})$, with a left positive eigenvector
$\bw_{\btp+\tpsi}$. A similar computation shows that 
$\bMM_{\btp+\tpsi}\bv_{\btp+\psi}=\rho_{\btp+\tpsi}\bv_{\btp+\psi}$. Corollary~\ref{periodic-perron-frobenius} 
in Appendix~\ref{proof-of-periodic} ensures that 
$\rho_{\btp+\tpsi}=\max(\spec(\bMM_{\btp+\tpsi}))$. 

\ms Now, using the definition of the pull back measure $\bmu_{\btp+\tpsi}:=\mu_{\btp+\tpsi}\circ\pi^{-1}$, 
and taking into account \eqref{approximated-Parry-measure} above, we have 
\begin{eqnarray*}
\lefteqn{\bmu_{\btp+\tpsi}[J_0\cdots J_n]}\\
&&
=\sum_{a_i,a'_i\in\bAA_{J_i}}\w_{\btp+\tpsi}((a_0)_{J_0}) 
\frac{\prod_{i=0}^{n-1} M_{\btp+\tpsi}((a'_i)_{J_i},(a_{i+1})_{J_{i+1}})}{\rho_{\btp+\tpsi}^n}     
\v_{\btp+\tpsi}((a'_n)_{J_n})\\
&&                                   
=\left(\sum_{a_0\in\bAA_{J_0}}\w_{\btp+\tpsi}((a_0)_{J_0})\right)\frac{\prod_{i=0}^{n-1}
\left(\sum_{a'_i\in\bAA_{J_i},a_{i+1}\in\bAA_{J_{i+1}}}
M_{\btp+\tpsi}((a'_i)_{J_i},(a_{i+1})_{J_{i+1}})\right)}{\bar{\rho}_{\btp+\tpsi}^n}\\
&&
\qquad\times \left(\sum_{a'_n\in\bAA_{J_n}}\v_{\btp+\tpsi}((a'_n)_{J_n}) \right)\\
&&
=\bw_{\btp+\tpsi}(J_0)\frac{\prod_{i=0}^{n-1}
\bMM_{\btp+\tpsi}(J_i,J_{i+1})}{\bar{\rho}_{\btp+\tpsi}^n}\bv_{\btp+\tpsi}(J_n),
\end{eqnarray*}
for each $n\in \nn$ and $J_0\cdots J_{n}\in (\AA')^{n+1}$, and the result follows.

\end{proof}


\appendix

\section{Proof of Proposition~\ref{periodic-proposition} (Periodic Approximation)}\label{proof-of-periodic} \

\ms The proof of Proposition~\ref{periodic-proposition} is based upon the following Theorem, which is a slight adaptation of Corollary 6.2 in~\cite{ChRU2005}. 

\ms \begin{theorem}[Perron--Frobenius: Primitive Case]~\label{perron-frobenius}
Let $B$ be a finite alphabet, and $\MM:B\times B\to \rr^+$ a primitive matrix.
Then there exists a unique $\rho\in\spec(\MM)$ such that $\rho:=\max|\spec(\MM)|$.
Associated to $\rho$ there are right and left eigenvectors $\v,\w$, such that $\w^{\dag}\v=1$.
Furthermore, for every probability vector $\x\in(0,1)^B$, and for each $m\in\nn$ we have
\[
\MM^m\x=\rho^m (\w^{\dag}\x)\v \thinspace \exp\left( \pm \frac{\tau^{\lfloor m/\ell\rfloor} \,
\ell d(\x,F\x)}{1-\tau}\right),
\]
where 
\begin{itemize}
\item[{\it i)}] $\ell$ is the primitivity index of $\MM$, {\em i.e.}, the smallest integer such that 
$\MM^\ell >0$, 
\item[{\it ii)}] $d$ is the projective distance in the simplex 
$\Delta_B:=\left\{\x\in(0,1)^B:\ |\x|_1=1\right\}$ \textup{(}$|\cdot|_1$ stands for the $\ell_1$ norm\textup{)} of 
probability vectors, 
\[
d(\x,\y)=\log\left(\max_{b\in B}\{\x(b)/\y(b)\}\right)-
\log\left(\min_{b\in B}\{\x(b)/\y(b)\}\right),
\]
\item[{\it iii)}] $F\x:=\MM\x/|\MM\x|_1$ is the action of the matrix $\MM$
on the simplex $\Delta_B$, and
\item[{\it v)}] $\tau=(1-\Gamma)/(1+\Gamma)$ with
\[
\Gamma:=\sqrt{\min_{a,b,c,d\in B}
\frac{\MM^{\ell}(a,b)\MM^{\ell}(c,d)}
{\MM^{\ell}(a,d)\MM^{\ell}(c,b)}},
\]
is the Birkhoff coefficient of $\MM^\ell$.
\end{itemize}
\end{theorem} 

\ms A rather direct consequence of the previous theorem is the following result.

\ms \begin{corollary}[Perron--Frobenius: Periodic Case]~\label{periodic-perron-frobenius}
Let $B$ be a finite alphabet, and $\MM:B\times B\to \rr^+$ an irreducible matrix of period $p$. 
Let $B:=\bigsqcup_{i=0}^p B_i$ be the partition such that 
\[
\MM_{\psi}^{p}=\left(\begin{matrix}
              \MM_0  &   0  &   0  &\cdots& 0\\
                0    &\MM_1 &   0  &\cdots& 0\\
                0    &   0  &\MM_2 &\cdots& 0\\
               \vdots&\vdots&\vdots&\ddots& \vdots\\
                0    &    0 &   0  &\cdots&\MM_{p-1}
     \end{matrix}\right)
\]
with $\MM_{i}:B_i\times B_i\to \rr^+$ primitive for each $0\leq i < p$. Let $\v_i$ and $\w_i$ denote 
the unique left and right eigenvectors associated to the maximal eigenvalue 
$\rho_i:=\max|\spec(\MM_i)|$, normalized such that $\w_i^{\dag}\v_i=1$.
Then,
$\rho=\sqrt[p]{\rho_i}$ is the maximal eigenvalue of $\MM$, with left and right eigenvectors
$\v:=(1/\sqrt{p})\otimes_{i=0}^{p-1}\v_i$ and $\w:=(1/\sqrt{p})\otimes_{i=0}^{p-1}\w_i$ respectively. 
These are the unique left and right positive eigenvectors satisfying $\w^{\dag}\v=1$.
Furthermore, if we fix $0\leq i <p$, and a probability vector $\x:=\otimes_{j=0}^{p-1}\x_j$, 
with $\x_j\in [0,1]^{B_j}$ such that $\x_j=0$ for $j\neq i$ and $\x_i>0$, then 
$\MM^{kp+r}\x=\otimes_{j=0}^{p-1}\y_j$, with $\y_j\in [0,1]^{B_j}$ such that $\y_j=0$ for $j\neq i-r$ and
\[
\y_{i-r}=
\frac{\rho^m}{p}
(\w_{i}^{\dag}\x_{i})\v_{i-r} \thinspace \exp\left( \pm\frac{\tau^{\lfloor k/\ell\rfloor} \,
\ell\, d_i(\x_i,F_i\x_i)}{1-\tau}\right) 
\]
where 
\begin{itemize}
\item[{\it i)}] $\ell$ is an upper bound for $\ell_i$, the primitivity index of the matrix $\MM_i$,
\item[{\it ii)}] $d_i$ is the projective distance in 
$\Delta_{i}:=\left\{\x_{i}\in(0,1)^{B_i}:\ |\x_{i}|_1=1\right\}$, 
\item[{\it iii)}] $F_i\x_i:=\MM_i\x_{B_i}/|\MM_i\x_{i}|_1$ is the action of the matrix 
$\MM^p$ on $\Delta_{i}$, and
\item[{\it iv)}] $\tau$ is an upper bound for the Birkhoff coefficient $\tau_i:=(1-\Gamma_i)/(1+\Gamma_i)$, 
with
\[
\Gamma_i:=\sqrt{\min_{a,b,c,d\in B_i}
\frac{\MM_i^{\ell}(a,b)\MM_i^{\ell}(c,d)}{\MM_i^{\ell}(a,d)\MM^{\ell}(c,b)}}\; .
\]
\end{itemize}
\end{corollary} 
\noindent{\bf {\em Proof of Proposition~\ref{periodic-proposition}}}

\noindent Let $\MM:=\MM_{\bpp}$, $\rho:=\rho_{\bpp}$, and for each $0\leq i<\barp$, let 
$\MM_{i}:=\MM_{\bpp,i}$, $\v_i:=\v_{\bpp,i}$, and $\w_i:=\w_{\bpp,i}$, as
defined in Subsection~\ref{subsection-transition-matrix}. Let $\ell$ be the maximum of the
primitivity indices of the matrices $\MM_i$, $0\leq i <\barp$, and denote by $i(a)$ the index of the set  
$\AA_i$ containing $a$. Applying Corollary~\ref{periodic-perron-frobenius} and using
\eqref{periodic-measure}, we obtain
\begin{eqnarray}
\nonumber
\lefteqn{\PP^{(k\barp)}_{\bpp}[\b_0^n]}\\
\nonumber
&=&\prod_{i=0}^{n-1}\MM(\b_i,\b_{i+1})\thinspace 
\frac{\delta_{\sb_n}^{\dag}\MM^{k\barp-n}\delta_{\sb_0}
    }{\sum_{b\in B}\delta_{b}^{\dag}\MM^{k\barp}\delta_{b}} \\
                             &=&
\frac{\prod_{i=0}^{n-1}\MM(\b_i,\b_{i+1})}{\rho^n}\thinspace
\frac{\left(\delta_{\sb_n}^{\dag}\MM_{i(\sb_n)}^\ell\right)^{\dag}\, 
                                              \v_{i(\sb_n)}\w_{i(\sb_0)}^{\dag}\,(\MM^\ell\delta_{\sb_0})
    }{\sum_{b\in B}\,(\delta_{b}^{\dag}\MM^\ell)\v_{i(b)}\w_{i(b)}^{\dag}\,(\MM^\ell\delta_{b})} \nonumber\\
                             & &
\times\exp\left(\pm \frac{C\,\tau^{\lfloor \frac{k-n/\barp}{\ell}\rfloor-2}}{1-\tau}\right) \nonumber \\
                             &=&
\frac{\w_{i(\sb_0)}(\b_0) \left(\prod_{i=0}^{n-1}\MM(\b_i,\b_{i+1})\right)\,\v_{i(\sb_n)}(\b_n)
    }{\rho^n \sum_{b\in B}\, \v_{i(b)}(\b)\w_{i(b)}(b) } \thinspace
    \exp\left(\pm\frac{C\,\tau^{\lfloor \frac{k-n/\barp}{\ell}\rfloor-2}}{1-\tau}\right) \nonumber\\
                             &=&
\frac{\w_{i(\sb_0)}(\b_0) \left(\prod_{i=0}^{n-1}\MM(\b_i,\b_{i+1})\right)\,\v_{i(\sb_n)}(\b_n)
    }{p \ \rho^n} \thinspace
    \exp\left(\pm\frac{C\,\tau^{\lfloor \frac{k-n/\barp}{\ell}\rfloor-2}}{1-\tau}\right) \nonumber\\
                                 &=&
\frac{\w(\b_0) \left(\prod_{i=0}^{n-1}\MM(\b_i,\b_{i+1})\right)\,\v(\b_n)
    }{\rho^n} \exp\left(\pm\frac{C\,\tau^{\lfloor \frac{k-n/\barp}{\ell}\rfloor-2}}{1-\tau}\right) \nonumber\\
                             &=&
\mu_{\bpp}[\b_0^n]\thinspace
     \exp\left(\pm\frac{C\,\tau^{\lfloor \frac{k-n/\barp}{\ell}\rfloor-2}}{1-\tau}\right), \nonumber
\end{eqnarray}
where $C:=2\ell\,\max_{0\leq i<\barp}\max_{b\in B_i} d_i(\MM_i^{\ell+1}\delta_{b},\MM_i^{\ell}\delta_{b})$, 
where $\delta_{b}:B_i\to\rr$ is the unit vector in the direction of $b$, and where $\tau=\max_{0\leq i<\barp} \tau_i$, where $\tau_i$ is the Birkhoff coefficient of $\MM_i^\ell$.

\ms Therefore we have 
\begin{equation}\label{periodic-approximation}
\mu_{\bpp}[\b_0^n]=\PP^{(k\barp)}_{\bpp}[\b_0^n]
\thinspace
\exp\left(\pm\frac{C\,\tau^{\lfloor \frac{k-n/\barp}{\ell}\rfloor-2}}{1-\tau}\right).
\end{equation}

\ms It remains to bound $\tau$ and $C$ to conclude the proof of the proposition.
Using the fact that 
$\tau \leq 1-\min_{0\leq i < \barp}\Gamma_i$, with $\Gamma_i$ as in the statement of 
Corollary~\ref{periodic-perron-frobenius}, we obtain
\begin{eqnarray*}
\tau &\leq& 1-\min_{0\leq i<\barp}\min_{a,b,d\in \AA_i}\frac{\MM_i^{\ell}(a,b)}{\MM_i^{\ell}(a,d)}
        =   1-\min_{a,b,d\in \AA}\frac{\MM_i^{\ell}(a,b)}{\MM_i^{\ell}(a,d)} \\ \\
     &\leq& 1-
 \frac{\sum_{\sa\in\AA^{\barp\ell}:\ \sa_0=a,\sa_{\barp\ell}=b} 
                                                    e^{\beta S_{\barp\ell}\phi(\sa)+S_{\barp\ell}\psi(\sa)}
     }{\sum_{\sa\in\AA^{\barp\ell}:\ \sa_0=a,\sa_{\barp\ell}=d} 
                                                    e^{\beta S_{\barp\ell}\phi(\sa)+S_{\barp\ell}\psi(\sa)}}, \\ \\
     &\leq&1-\frac{e^{-\barp\ell(\beta\|\phi\|_{\infty}+\|\psi\|_{\infty})+\log(\#\AA))}}{e^{\barp\ell(\beta\|\phi\|_{\infty}+\|\psi\|_{\infty}+\log(\#\AA))}}\\ \\
     &\leq&1-\exp(-2\barp\ell(\beta\|\phi\|_{\infty}+\|\psi\|_{\infty}+\log(\#\AA))).
\end{eqnarray*}
Now, for each $0\leq i < \barp$ and $b\in \AA_i$ we have
\begin{eqnarray*}
\lefteqn{d_i(\MM_i^{\ell+1}\delta_b,\MM_i^{\ell}\delta_b)}\\
                                                &=  &
  \max_{a,a'\in \AA_i}\log\left(\frac{\di
  \sum_{\sa\in\AA^{\barp(\ell+1)}\atop \sa_0=a,\sa_{\barp(\ell+1)}=b}e^{S_{\barp(\ell+1)}(\beta\phi+\psi)(\sa)}
  \sum_{\sa\in\AA^{\barp\ell}\atop \sa_0=a',\sa_{\barp\ell}=b}e^{S_{\barp\ell}(\beta\phi+\psi)(\sa)}                            
                                     }{\di
  \sum_{\sa\in\AA^{\barp(\ell+1)}\atop \sa_0=a',\sa_{\barp(\ell+1)}=b}e^{S_{\barp(\ell+1)}(\beta\phi+\psi)(\sa)}
  \sum_{\sa\in\AA^{\barp\ell}\atop \sa_0=a,\sa_{\barp\ell}=b}e^{S_{\barp\ell}(\beta\phi+\psi)(\sa)} }\right)\\ \\
                                                &\leq &
 \log\left( \frac{e^{\barp(2\ell+1)(\beta\|\phi\|_{\infty}+\|\psi\|_{\infty}+\log(\#\AA))}
      }{e^{-\barp(2\ell+1)(\beta\|\phi\|_{\infty}+\|\psi\|_{\infty})\log(\#\AA))}}\right) \\ \\
      &\leq& 2\barp(2\ell+1)(\beta\|\phi\|_{\infty}+\|\psi\|_{\infty}+\log(\#\AA)),                                                                                               
\end{eqnarray*}
where $\|\cdot\|$ denotes the supremum norm. With these two bounds, and taking into account 
\eqref{periodic-approximation}, we obtain
\[
\mu_{\bpp}[\a_0^n]=\PP^{(k\barp)}_{\bpp}[\a_0^n]\,\exp\left(\pm (\beta\, n_{\phi}+n_{\psi})
e^{\beta s_{\phi}+s_{\psi}}\left(1-e^{-(\beta s_{\phi}+s_{\psi})}\right)^{\frac{k}{\ell}-\frac{n}{\barp\ell}-3}\right),
\]
for some positive constants $s_\phi$, $s_\psi$, $n_\phi$ and $n_\psi$.
Since $k\barp > e^{\gamma\beta}+n$, and since by assumption $\gamma>s_\phi$, we have 
$$
\frac{k}{\ell}-\frac{n}{\barp\ell}-3 \geq 2\beta(\gamma-s_{\phi})\,e^{\beta s_{\phi}+s_{\psi}}\;\text{for}\; \beta\; 
\text{large enough}.
$$
Taking into account that $1-e^{-(\beta s_{\phi}+s_{\psi})}\leq \exp\left(-e^{-(\beta s_{\phi}+s_{\psi})}\right)$,
we finally obtain, for $\beta$ large enough,
\begin{eqnarray*}
\mu_{\bpp}[\b_0^n] & = & \PP^{(k\barp)}_{\bpp}[\b_0^n]\,
  \exp\left(\pm (\beta\, n_{\phi}+n_{\psi})e^{-2\beta(\gamma-s_{\phi})}\right) \\
  &= &
  \PP^{(k\barp)}_{\bpp}[\b_0^n]\, \exp\left(\pm e^{-\beta(\gamma-s_{\phi})}\right)
\end{eqnarray*} 
for all $\b_0^n$. The proof of Proposition \ref{periodic-proposition} is now finished.
\endproof

\begin{remark}
We have $s_{\phi}=2\barp\ell\|\phi\|_{\infty}$, $s_{\psi}=2\barp\ell(\|\psi\|_{\infty}+\log(\#\AA))$, 
and $n_{*}:=2(2\ell+1)s_{*}$ where $*=\phi,\psi$. 
\end{remark}

\section{Auxiliary Inequalities of Lemma~\ref{factorization-Lemma}}~\label{auxiliary-inequalities}\

\ms \subsection{Incursion Length}\label{proof-incursion-length}\

\ms \begin{proposition}[The Incursion Time is Exponential]\label{incursiontime}
Let us suppose $[\b_0^n]\cap\bX_J\neq\emptyset$ for some $1\leq J\leq N_\phi$.
If $q=\lfloor e^{\eta\,\beta}\rfloor$, with $-\phi_g/4 < \eta < -\phi_g/3$, and if $p\geq 2q^2+1$, then we have
\[
\sum_{\sc\in\sper_{p}(X)\cap[\sb_0^n]}e^{S_{p}(\bpp)(\sc)}=\exp\left(\pm 2e^{\beta \frac{\phi_g}{4}}\right)
\sum_{\sc\in X^{(p,q)}\cap[\sb_0^n]}e^{S_{p}(\bpp)(\sc)},
\]
for $\beta$ large enough.
\end{proposition}

\ms \begin{proof} We obviously have 
\[
\sum_{\sc\in\sper_{p}(X)\cap[\sb_0^n]}e^{S_{p}(\bpp)(\sc)}\geq\exp \left(- 2e^{\beta \frac{\phi_g}{4}}\right)
\sum_{\sc\in X^{(p,q)}\cap[\sb_0^n]}e^{S_{p}(\bpp)(\sc)}.
\] 

\ms For each $a$ and $a'\in \AA$, let us denote by $a\cira a'$ a fixed circuit in $\GG_X$
of period $p_0=p_0(a,a')$, containing both $a$ and $a'$. By choosing $p_0$ the minimal integer 
for which such circuit exists, we ensure that $p_0(a,a')\leq 2\#\AA$.
Let us denote by $a\to a'$ and $a'\to a$ the path segments composing $a\cira a'$, and
by $w(a\to a')$ and $w(a'\to a)$ the corresponding $X$--admissible word. Let us also denote by $p(a\cira a')$ the periodic point in $\per_{p_0}(X)\cap [a]$ defined by the circuit $a\cira a'$.

\ms To each periodic orbit in $\a\in\per_p(X)\cap[\b_0^n]$ such that $\a_{p-q^2-1}=a'$ and $\a_{q^2}=a$, 
we associate the periodic points 
\begin{eqnarray*}
\a_{\rm int}&:=&\left(\a_0^{q^2-1}\,w(a\to a')\,\a_{p-q^2}^{p-1}\right)^{\infty}\in\per_{p_1}(X)
                 \cap\sigma^{-q^2}\left[w(a\to a')\right],\\
\a_{\rm ext}&:=&\left({a'\to a}\,\a_{q^2+1}^{p-q^2-2}\right)^{\infty}\in\per_{p_2}(X)
                \cap\left[w(a'\to a)\right],
\end{eqnarray*}
with $p_1=p_1(a,a'):=2q^2+|a\to a'|$ and $p_2:=p_2(a,a')=p-2q^2-2+|a'\to a|$. Using this notation we can write,
\[
\sum_{\sa\in\sper_{p}(X)\cap[\sb_0^n]}e^{S_{p}(\bpp)(\sa)}=
        \sum_{a,a'\in \AA}e^{-S_{p_0}(\beta\phi+\psi)\left(p(a\cira a')\right)}
         e^{S_{p_1}(\beta\phi+\psi)(\sa_{\rm int})}e^{S_{p_2}(\beta\phi+\psi)(\sa_{\rm ext})},                                                           
\]
which yields 
\begin{equation}\label{ub-incursion-length}
\sum_{\sa\in\sper_{p}(X)\cap[\sb_0^n]}e^{S_{p}(\bpp)(\sa)}=
\end{equation}
$$
\sum_{a,a'\in \AA}e^{-S_{p_0}(\beta\phi+\psi)\left(p(a\cira a')\right)}
\times
$$
$$
\sum_{\sa\in \sper_{p_1}(X)\cap [\sb_0^n]\cap\sigma^{-q^2}[w(a\to a')]} e^{S_{p_1}(\beta\phi+\psi)(\sa)}
\times\sum_{\sa\in \sper_{p_2}(X)\cap\left[w(a'\to a)\right]}e^{S_{p_2}(\beta\phi+\psi)(\sa)}.                                                 
$$
Let us now study the {\it interior sums} 
$\sum_{\sa\in \sper_{p_1}(X)\cap [\sb_0^n]\cap\sigma^{-q^2}\left[w(a\to a')\right]}e^{S_{p_1}(\beta\phi+\psi)(\sa)}$.

\ms Each periodic point $\a\in \per_{p_1}(X)\cap [\b_0^n]\cap\sigma^{-q^2}\left[w(a\to a')\right]$ defines 
a circuit $C(\a)$ in $\GG_{X}$. We decompose this circuit into its {\it incursion--excursion path segments},
\[
C(\a):=a_1\lds a'_1\lds a_2\lds a'_2\lds \cdots \lds a_\kappa\lds a'_\kappa\lds  a_1,
\] 
which are defined as follows:
\begin{itemize} 
\item[a)] the segment $a_1\lds a'_1$ lies on $\GG_J\subset \GG_{X}$, the digraph associated to the 
heavy component $\bX_J$ such that $[\b_0^n]\cap\bX_J\neq \emptyset$;
\item[b)] none of the paths $a'_i\lds a_{i+1}$, for $1\leq i\leq\kappa=\kappa(\a)$, includes arrows from 
the digraph $\GG_{\bX}$ defining $\bX$;
\item[c)] each path $a_i\lds a'_i$, for $1\leq i\leq\kappa$, lies on $\GG_{\bX}$.
\end{itemize}

\ms We extend $C(\a)$ by adding, for each $1\leq i\leq\kappa$, a circuit $a_i\cira a'_i$ lying on 
the same transitive component of $\bX$ as $a_i$, containing both $a_i$ and $a'_i$, and having minimal length. 
Since all the added circuits lie in $\GG_{\bX}$, then the extended circuit 
\[
\begin{array}{ccccccc}
C_{\rm ext}(\a)
  :=a_1\lds&a'_1 &\lds  a_2\lds & a'_2&\lds \cdots \lds a_\kappa\lds&a'_\kappa&\lds a_1 \\
           &\cira&             &\cira&                             & \cira &         \\
           & a_1 &             & a_2  &                            & a_\kappa&             
\end{array}  
\]
is such that $\phi(C(\a))=\phi(C_{\rm ext}(\a))$. Since all the added circuits have minimal length, we also
have $\psi(C(\a))\leq \psi(C_{\rm ext}(\a))+2\,\kappa\,\#\AA\,\|\psi\|_{\infty}$. 

\ms We reorganize the path segment in $C_{\rm ext}(\a)$ in order to obtain
\[
C_{\rm ext}(\a)=\sum_{i=1}^\kappa a_i\cirb a'_i + a_1\to a'_1\lds a_2\to a'_2\lds 
\cdots\lds a_i\to a'_i\lds\cdots\lds a_\kappa\to a'_\kappa\lds a_1,
\]
where, for each $1\leq i\leq\kappa$, the circuit $a\cirb a'$ is obtained by concatenation of $a\lds a'$ and the 
path segment $a'\to a$ of the added circuit $a\cira a'$. For the complementary circuit 
\[
C'(\a):=a_1\to a'_1\lds a_2\to a'_2\lds \cdots\lds a_i\to a'_i\lds\cdots\lds a_\kappa\to a'_\kappa\lds a_1,
\]
we replace the segments $a_i\lds a'_i$ in $C(\a)$, by the paths of minimal length $a_i\to a'_i$
appearing in $a_i\cira a'_i$. All the circuits $a_i\cirb a'_i$ lie in $\GG_{\bX}$, therefore they 
maximize $\phi$. Notice also that the complementary circuit $C'(\a)$ does 
not include any circuit maximizing $\phi$, therefore $\phi(C'(\a))\leq |C'(\a)|\ \phi_g$. 

\ms We can bound from above the sum of the potentials $\phi$ and $\psi$ on the circuit $C(\a)$ by the same sums
over the extended circuit $C_{\rm ext}(\a)$ as follows: 
\begin{eqnarray*}
\phi(C(\a))   &  =  &\phi(C_{\rm ext}(\a))=\phi(C'(\a))\leq |C'(\a)|\thinspace \phi_g \\
\psi(C(\a)) &\leq  &\psi(C_{\rm ext}(\a))+2\,\kappa\,\#\AA\|\psi\|_{\infty}=
     \sum_{i=1}^\kappa\psi\left({a_i\cirb a'_i}\right)+\psi(C'(\a))+2\,\kappa\,\#\AA\|\psi\|_{\infty}\nonumber\\
              & \leq &\sum_{i=1}^\kappa\psi\left({a_i\cirb a'_i}\right)
                      +(2\#\AA+1)\|\psi\|_{\infty}\thinspace |C'(\a)|, \nonumber
\end{eqnarray*}
where we use the fact that $\kappa\leq |C'(\a)|$.

\ms We can group the periodic points in $\PP_1:=\per_{p_1}(X)\cap [\b_0^n]\cap\sigma^{-q^2}
\left[w(a\to a')\right]$ 
according to the number of excursion paths they contain. By doing so we have 
\begin{eqnarray*}
\sum_{\sa\in\PP_1} e^{S_{p_1}(\beta\phi+\psi)(\sa)}
                   &\leq & 
\sum_{\kappa=0}^{\lfloor p_1/2\rfloor}\sum_{\sa\in\PP_1\atop\kappa(\sa)=\kappa}e^{S_{p_1}(\beta\phi+\psi)(\sa)}\\
                   &\leq &
\sum_{\kappa=0}^{\lfloor p_1/2\rfloor}\sum_{\sa\in\PP_1\atop\kappa(\sa)=\kappa}
                                         e^{(\beta\phi_g+(2\#\AA+1)\|\psi\|_{\infty})|C'(\sa)|}
                                        \prod_{i=1}^\kappa e^{\psi\left({a_i\cirb a'_i}\right)}.
\end{eqnarray*}  
Now we group all the periodic points in $\{\a\in\PP_1:\ \kappa(\a)=\kappa\}$ in classes 
defined by the total total length of the complementary circuit, $L(\a):=|C'(\a)|$, the
lengths of the incursion and the excursion segments, $m_i:=|a_i\lds a'_i|$ and $n_i:=|b_i\lds b'_i|$
respectively, and the location of the origin inside the first incursion $a_1\lds a'_1$. 
Taking into account Corollary~\ref{periodic-perron-frobenius} and the fact that $P(\psi|\bX)=P(\psi|\bX_J)$, 
for $J=1\ldots N_\phi$, we obtain
\begin{eqnarray*}
\sum_{\sa\in\PP_1} 
e^{S_{p_1}(\beta\phi+\psi)(\sa)}&\leq&\sum_{\sa\in\PP_1\atop L(\sa)<L_0} e^{S_{p_1}(\beta\phi+\psi)(\sa)}+
        \sum_{\kappa=0}^{\lfloor L/2 \rfloor}\sum_{L=L_0}^{p_1}e^{(\beta\phi_g+(2\#\AA+1)\|\psi\|_{\infty}+\log(\#\AA))L} 
        \\
                                &    & \times\thinspace \sum_{\sum m_i=p_1-L\atop\sum n_i=L}
                 m_1 K_{\psi} e^{S_n\psi(\sb_0^n)}\MM_{\psi,J}^{m_1-n}(\b_n,\b_0)
                                          \prod_{i=2}^\kappa K_{\psi}\trace\left(\MM_{\psi,J}^{m_i}\right),
\end{eqnarray*}
for all $\beta\geq ((2\#\AA+1)\|\psi\|_{\infty}+\log(\#\AA))/|\phi_g|$. The integer $L_0\geq 2$ will be fixed later on.
Here we bound the sums of factors $e^{\psi(a _i\cirb a'_i)}$ by a constant multiple  
of $\trace(\MM_{\psi,J}^{m_i})$, except for the term with $i=1$, which we bound by a constant factor of
$e^{S_n\psi(\sb_0^n)}\MM_{\psi,J}^{m_1-n}(\b_n,\b_0)$. This last bound follows from the fact that 
$a_1\cirb a'_1$ includes the path segment $(\b_0,\b_1,\ldots,\b_n)$. 
The constant $K_{\psi}\geq \#\AA$ is taken large enough to include the counting of transitive components
in $\bX$, and to compensate the differences in trace among transitive components and the difference 
between $|a_i\cirb a'_i|$ and $m_i$ for each 
$1\leq i\leq \kappa$. The factor $m_1$ takes into account the all the possible locations of the origin with
inside the first incursion. We are also using the fact that $2\kappa(\a)\leq L(\a)$. 

\ms Now, the normalization $P(\psi|\bX)=0$ ensures that
\begin{eqnarray*}
\sum_{\sa\in\PP_1\atop L(\sa)\geq L_0} e^{S_{p_1}(\beta\phi+\psi)(\sa)}&\leq& e^{S_n\psi(\sb_0^n)}
  \sum_{L=L_0}^{p_1}e^{(\beta\phi_g+(2\#\AA+1)\|\psi\|_{\infty}+\log(\#\AA))L}\ p_1\ \sum_{\kappa=0}^{\lfloor L/2\rfloor}
                                                    \binom{p_1}{2\kappa-1} D_{\psi}^{\kappa}\\
                                                                  &\leq&
  e^{S_n\psi(\sb_0^n)}\sum_{L=L_0}^{p_1}
           e^{(\beta\phi_g+(2\#\AA+1)\|\psi\|_{\infty}+\log(\#\AA)+\log(p_1)+\log(D_\psi)+\log(L)/L)L}.\nonumber
\end{eqnarray*}
The constant $D_{\psi}\geq K_\psi$ includes an upper bound for the factors
$\trace(\MM_{\phi,J}^{m})$. Since $p_1\leq 2 e^{2\beta\,\eta} +|a\to a'|$, with 
$\eta < -\phi_g/3$, it follows from the previous inequality that
\begin{eqnarray}\label{ub-large-L}
\sum_{\sa\in\PP_1\atop L(\sa)\geq L_0} e^{S_{p_1}(\beta\phi+\psi)(\sa)}&\leq& e^{S_n\psi(\sb_0^n)}
         \sum_{L=L_0}^{p_1} e^{(\beta\phi_g+(2\#\AA+1)\|\psi\|_{\infty}+\log(\#\AA)+\log(p_1)+\log(D_\psi)+e^{-1})L}\\ 
                                                                       &\leq& 
           2 \ e^{S_n\psi(\sb_0^n)} e^{\beta L_0\frac{\phi_g}{3}},\nonumber
\end{eqnarray}
for all $\beta$ greater than a convenient $\beta(\phi,\psi,\AA)$.

\ms Let 
$$
k_1:=\min\left\{k\in\nn:\ \per_{k p_J}(\bX_J)\cap[\a_{p-q^2}^{q^2-1}]\neq \emptyset, \, \,\forall \,\, 
\a_{p-q^2}^{q^2-1} \, \,\bX_J\text{--admissible with }\, \a_0^n=\b_0^n\right\}\cdot
$$
It is not hard to verify that such a minimum exists.
To each $\a\in\per_{p_1}(\bX_J)\cap[\b_0^n]$ such that $\a_{p-q^2}^{q^2-1}$ is $\bX_J$--admissible,
we associate a fixed periodic point 
$\a_{\rm int}\in\per_{k p_J}(\bX_J)\cap[\b_0^n]\cap \sigma^{q^2}[\a_{p-q^2}^{q^2-1}]$, and the 
periodic point $p(a\cira a')\in \per_{p_0}(X)$.
Clearly $S_{p_1}\phi(\a)=S_{p_0}\phi\left(p(a\cira a')\right)\geq -2\,\#\AA\,\|\phi\|_{\infty}$ and 
$S_{p_1}\psi(\a)\geq S_{k_1 p_J}\psi(\a_{\rm int})-2\,\#\AA\,\|\psi\|_{\infty}$. From this, and taking
into account Corollary~\ref{periodic-perron-frobenius}, it follows that
\begin{eqnarray*}
\lefteqn{\sum_{\sa\in\PP_1\atop[\sa_{p-q^2}^{q^2-1}]\cap\bX_J\neq \emptyset} e^{S_{p_1}(\beta\phi+\psi)(\sa)}
                 }\\ & \geq & e^{-2\,\#\AA\,(\beta\|\phi\|_{\infty}+\|\psi\|_{\infty})}\#\AA^{2q^2-k_1p_J}\; 
\sum_{\sa\in\sper_{k_1 p_J}(\bX_J)\cap[\sb_0^n]\neq\emptyset} e^{S_{k_1p_J}\psi(\sa)}\\ 
                 &\geq& e^{-2\,\#\AA\,(\beta\|\phi\|_{\infty}+\|\psi\|_{\infty})}\#\AA^{2q^2-k_1p_J}e^{S_n\phi(\sb_0^n)}
                        \MM_{\psi,J}^{k_1p_J-n}(\b_n,\b_0)\\
                 &\geq& B\thinspace e^{S_n\psi(\sb_0^n)} e^{-2\,\#\AA\,(\beta\|\phi\|_{\infty}+\|\psi\|_{\infty})}
\end{eqnarray*}
with $B=B(\psi):=\#\AA^{2q^2-k_1p_J}\min_{\b_0,\b_n}\v_{\psi,J}(\b_n)\w_{\psi,J}(\b_0)e^{-C_\psi}$. This and 
\eqref{ub-large-L} imply 
\begin{equation}\label{ineq-L_0}
\sum_{\sa\in\PP_1\atop L(\sc)\geq L_0} e^{S_{p_1}(\beta\phi+\psi)(\sa)} \leq
e^{\beta\phi_g/3}\sum_{\sa\in\PP_1\atop[\sa_{p-q^2}^{q^2-1}]\cap\bX_J\neq \emptyset} e^{S_{p_1}(\beta\phi+\psi)(\sa)},
\end{equation}
for all 
$\beta\geq \max(\beta(\phi,\psi,\AA),(6\,\#\AA\,\|\psi\|_{\infty}+3\log(2)-3\log(B))/(|\phi_g|(L_0-1)-6\,\#\AA\,\|\phi\|_{\infty}))$,
as long as $L_0 > L_{\phi}:=\lceil 6\,\#\AA\,\|\phi\|_{\infty}/|\phi_g| \rceil+ 1$.

\ms Let us now consider the set 
$\PP_2:=\{\a\in \PP_1:\ L(\a) \leq  L_{\phi}:=\lceil 6\,\#\AA\,\|\phi\|_{\infty}/|\phi_g|\rceil + 1\}$. 
Consider a periodic points in $\a\in \PP_2$ whose associated circuit contains a segment 
$a'_{i-1}\lds a_i\lds a'_i\lds a_{i+1}$ with a incursion $a_i\lds a'_i$ of length $m$ into a 
non--heavy component. We can replace $a'_{i-1}\lds a_i\lds a'_i\lds a_{i+1}$ by a segment
$a'_{i-1}\lds c_i\lds c'_i\lds a_{i+1}$ such that $c_i\lds c'_i$ lies in $\GG_J$ and 
$|c_i\lds c'_i|\geq m-2\#\AA$, obtaining a new periodic point $\a'\in \PP_2$ such that
\begin{eqnarray*}
S_{p_1}\phi(\a) &\leq& S_{p_1}\phi(\a') + 2(L_{\phi}+\#\AA)\|\phi\|_{\infty}, \\
S_{p_1}\psi(\a) &\leq& S_{p_1}\psi(\a') + 2(L_{\phi}+\#\AA)\|\psi\|_{\infty}+
\psi\left({a_i\cirb a_i'}\right)-\psi\left({c_i\cirb c_i'}\right).                  
\end{eqnarray*}
Following this prescription we can replace all the occurrences of incursions in non--heavy components of length 
$m\geq \epsilon p_1$ by incursions into $\bX_J$ of about the same length, and taking into account 
and the cominatorics of the replacements and the fact that $P(\psi|\bX_J)=P(\psi|\bX)=0$, we obtain
\begin{eqnarray*}
\lefteqn{\sum_{\sa\in\PP_2} e^{S_{p_1}(\beta\phi+\psi)(\sa)}}\\
 &\leq & \!\!\!\!\!\!
\sum_{\sa\in\PP_2\atop |a_i\lds a'_i| > \epsilon p_1\Rightarrow a_i\lds a'_i \text{ in } \GG_J} \!\!\!\!\!\!\!\!\!\!
                 e^{S_{p_1}(\beta\phi+\psi)(\sa)} \sum_{k=0}^{L_{\phi}}\binom{L_{\phi}}{k} 
                                \left(e^{(\epsilon p_1-2\#\AA)(P(\psi|\bX')-P(\psi|\bX))+
                                r_{\psi}+\beta\,r_{\phi}}\right)^{k}\\
                                                     &\leq & 
              \left(1+e^{(\epsilon p_1-2\#\AA)P(\psi|\bX')+r_{\psi}+\beta\,r_{\phi}}\right)^{L_{\phi}}
     \hskip -20pt              
     \sum_{\sa\in\PP_2\atop |a_i\lds a'_i| > \epsilon p_1\Rightarrow a_i\lds a'_i \text{ in } \GG_J} 
     \hskip -20pt e^{S_{p_1}(\beta\phi+\psi)(\sa)} ,\nonumber
\end{eqnarray*}  
where $r_{\psi}:=\log(\#\AA)L_{\phi}+\log(D_{\psi})+2(L_{\phi}+\#\AA)\|\psi\|_{\infty}$, 
$r_{\phi}:=2(L_{\phi}+\#\AA)\|\phi\|_{\infty}$, and $\bX'=\bX\setminus\cup_{K=1}^N\bX_K$. 
Since $P(\psi|\bX')<0$, then we have
follows that
\begin{eqnarray}\nonumber
\sum_{\sa\in\PP_2} e^{S_{p_1}(\beta\phi+\psi)(\sa)}&\leq &
 \left(1+2L_{\phi}e^{(\epsilon p_1-2\#\AA)P(\psi|\bX')+r_{\psi}+\beta\,r_{\phi}}\right)\hskip -20pt
\sum_{\sa\in\PP_2\atop |a_i\lds a'_i| > \epsilon p_1\Rightarrow a_i\lds a'_i \text{ in } \GG_J} 
           \hskip -20pt e^{S_{p_1}(\beta\phi+\psi)(\sa)} \\
                                                    &\leq &
                     \left(1+e^{\beta\frac{\phi_g}{3}}\right)\hskip -20pt
\sum_{\sa\in\PP_2\atop |a_i\lds a'_i| > \epsilon p_1\Rightarrow a_i\lds a'_i \text{ in } \GG_J} 
           \hskip -20pt e^{S_{p_1}(\beta\phi+\psi)(\sa)} \label{ub-short-length}
\end{eqnarray}
as long as 
$\epsilon p_1 > (\beta(|\phi_g|/3+r_{\phi})+\log(2L_{\phi})+2\#\AA|P(\psi|\bX')|+r_{\phi})/|P(\psi|\bX')|$.
We ensure this by taking $\epsilon=q^{-1}\geq e^{-\beta\,\eta}$ and $\beta$ larger than a convenient 
$\beta(\eta)$.

\ms Let us group the periodic points in
\[
\PP_3:=\left\{\a\in\PP_2:\ |a_i\lds a'_i|\geq \epsilon p_1=2q+|a\to a'|/q\, \Rightarrow\, a_i\lds a'_i 
\text{ in } \cup_{K\leq N_\phi}\GG_K\right\},
\]
by classes $\QQ:=\QQ(a'_{i-1}\lds a_i,\,\, a_i\lds a_i'\in\bigcup_{K=N_\phi+1}^{N_\phi}\GG_K,\,\,a_i,\,a'_i\in
\bigcup_{K=N_\phi+1}^{N_\phi}\bAA_{K},\,\, k_0)$
defined by the excursion segments $a'_{i-1}\lds a_i$, the incursion into non--heavy components
$a_i\lds a'_i$, the incursion input--output vertices $a_i$ and $a'_i$ into heavy components, and
the location $k_0$ of the origin in the first incursion segment.  
By definition, two points $\a,\a'$ in the same class are such that $S_{p_1}\phi(\a)=S_{p_1}\phi(\a')$. 
Now, for each $\a\in\PP_3$ we have 
\[
\sum_{a_i\lds a'_i \text{ in } \cup_{K\leq N_\phi}\GG_K}|a_i\lds a'_i|\geq 
(2q^2+|a\lds a'|)\left(1-\frac{2L_{\phi}}{q}\right).
\]
Let $p*={\rm lcm}(p_K:\ 1\leq K\leq N)$ and consider a refinement $\QQ=\bigsqcup\QQ_{\{m_i:i=1,\ldots, \kappa\}}$ 
of a particular class 
$\QQ:=\QQ(a'_{i-1}\lds a_i,\,\, a_i\lds a_i'\in\bigcup_{K=N_\phi+1}^{N_\phi}\GG_{K},\,\,a_i,\,a'_i\in
\bigcup_{K=N_\phi+1}^{N_\phi}\bAA_{K},\,\, k_0)$ 
in subclasses 
$\QQ_{\{m_i:i=1,\ldots, \kappa\}}\subset\QQ$ defined by the data of the lengths $m_i=|a_i\lds a'_i|$ 
of the incursions into heavy components. 

\ms A particular subclass $\QQ_{\{m_i:i=1,\ldots, \kappa\}}\subset\QQ$ contributes with   
\begin{eqnarray}\label{same-amount}
\sum_{\sa\in\QQ_{\{m_i:i=1,\ldots, \kappa\}}}
e^{S_{p_1}(\beta\phi+\psi)(\sa)}&=&e^{\beta\phi(C')+\psi(\tC)}
                                   \prod_{i=1}^{\kappa'}\MM^{m_i}_{\psi,J_i}(a_i,a'_i)\\
                           &=&e^{\pm \kappa'\, C_{\psi}}e^{\beta\phi(C')}e^{\psi(\tC)}
                                \prod_{i=1}^{\kappa'}\v_{J_i}(a_i)\w_{J_i}(a'_i),\nonumber
\end{eqnarray}  
to the sum $\sum_{\sa\in\PP_2} e^{S_{p_1}(\beta\phi+\psi)(\sa)}$.
Here $C':=C'(\QQ)$ is the circuit defined by the concatenation of the excursions, the incursion 
into non--heavy components, and path segments $a_i\to a'_i$ in $\cup_{K\leq N_\phi}\GG_K$ of minimal length. 
The complement $\tC:=\tC(\QQ)$ is a disjoint union of path segments, formed by the concatenation 
of all the excursions and the incursion into non--heavy components.
The constant
$\kappa':=\kappa'(\QQ) < \kappa:=\kappa(\QQ)$ is the number of times a periodic point in the chosen collection visits a heavy component.
Hence, all the subclasses $\QQ_{\{m_i:i=1,\ldots, \kappa\}}\subset \QQ$ contribute with about 
the same amount to the sum $\sum_{\sa\in\PP_2} e^{S_{p_1}(\beta\phi+\psi)(\sa)}$.
Their contributions differ at most by a factor in the range $\exp(\pm 2\kappa' C_\psi)=\exp(\pm 2L_\phi C_\psi)$.

\ms Now, to each subclass $\QQ_{\{m_i:i=1,\ldots, \kappa\}}\subset \QQ$ 
such that $\min(k_0,m_1-k_0) < q$ 
(remember that $k_0$ is the location of the origin with respect to the first incursion), there are at least 
\begin{equation}\label{proportion}
\frac{1}{2qp*}\left(\frac{(2q^2+|a\lds a'|)(1-2L_\phi/q)}{L_{\phi}}-2q\right)\geq \frac{q}{2L_{\phi}p^*}
\end{equation}
subclasses $\QQ_{\{m_i:i=1,\ldots, \kappa\}}\subset \QQ$ with $\min(k_0,m_1-k_0)\geq q$, for all $q$ 
sufficiently large.
Indeed, we can increase the length of the first incursion at both sides of the origin 
by decreasing the length of another incursions.
The size of this length change has to be a common multiple of the periods of the
heavy components involved.
The decrease of length can be done for length sufficiently large.
Each length increase can be associated to at most $2q$ subclasses with $m_1\leq 2q$.
All these lengths changes can be done for all large values of $q$, corresponding to values of
$\beta$ greater than a convenient $\beta(p^*)$.

\ms Taking into account that all subclasses $\QQ_{\{m_i:i=1,\ldots, \kappa\}}\subset \QQ$ contribute with about
the same amount to the sum $\sum_{\sa\in\PP_2} e^{S_{p_1}(\beta\phi+\psi)(\sa)}$ \eqref{same-amount}, 
and since the majority of those subclasses are such that  
$\min(k_0,m_1-k_0)\geq q$ \eqref{proportion}, then, since $\eta > -\phi_g/4$ we have
\begin{eqnarray}\label{last-inequality}
\sum_{\sa\in\PP_3}e^{S_{p_1}(\beta\phi+\psi)(\sa)} 
          &  = & \sum_{b_i\lds b'_i,\,\, a_i\lds a_i'\in\bigcup_{K=N_\phi+1}^{N_\phi}\GG_{K},\atop 
                  a_i,\,a'_i\in \bigcup_{K=N_\phi+1}^{N_\phi}\bAA_{K},\,\, k_0} 
                  \sum_{m_i=|a_i\lds a'_i|} \sum_{\sa\in \QQ_{\{m_i\}}}e^{S_{p_1}(\beta\phi+\psi)(\sa)} \nonumber\\  \nonumber\\     \nonumber            
          &\leq& \left(1+e^{2L_{\phi}\, C_{\psi}}\frac{2L_{\phi}p^*}{q}\right) \thinspace\times \\ \nonumber\\
                 && \sum_{b_i\lds b'_i,\,\, a_i\lds a_i'\in\bigcup_{K=N_\phi+1}^{N_\phi}\GG_{K},\atop 
                  a_i,\,a'_i\in \bigcup_{K=N_\phi+1}^{N_\phi}\bAA_{K},\,\, k_0}
                  \sum_{m_i=|a_i\lds a'_i|\atop \min(k_0,m_1-k_0)\geq q} 
                  \sum_{\sa\in \QQ_{\{m_i\}}}e^{S_{p_1}(\beta\phi+\psi)(\sa)} \nonumber\\ \nonumber\\
          &\leq&\left(1+e^{\beta\phi_g/4}\right) 
                  \sum_{\sa\in\PP_3\atop \min(k_0,m_1-k_0)\geq q}e^{S_{p_1}(\beta\phi+\psi)(\sa)},
                  \end{eqnarray}
for $\beta$ larger than $\max(\beta(\eta),\beta(p^*))$. 

\ms We can conclude now: \eqref{ineq-L_0}, \eqref{ub-short-length} and \eqref{last-inequality}
imply that 
\[
\sum_{\sa\in\PP_1}e^{S_{p_1}(\beta\phi+\psi)(\sa)}
\leq
\]
\[
\left(e^{\beta\phi_g/3}+\left(1+e^{\beta\phi_g/3}\right)\left(1+e^{\beta\phi_g/3}\right)
\left(1+e^{\beta\phi_g/4}\right)\right)
\sum_{\sa\in\PP_1\atop [\sa_{p-q}^{q-1}]\cap\bX\neq\emptyset}e^{S_{p_1}(\beta\phi+\psi)(\sa)},
\]
for all 
$\beta\geq\max(\beta(\phi,\psi,\AA),(6\,\#\AA\,\|\psi\|_{\infty}+3\log(2)-3\log(B))/(|\phi_g|(L_0-1)-6\,\#\AA\,\|\phi\|_{\infty}),\beta(\eta),\beta(p^*))$ 
whence
\[\sum_{\sa\in\PP_1}e^{S_{p_1}(\beta\phi+\psi)(\sa)}
\leq \left(1+2e^{\beta\phi_g/4}\right)
\sum_{\sa\in\PP_1\atop [\sc_{p-q}^{q-1}]\cap\bX\neq\emptyset}e^{S_{p_1}(\beta\phi+\psi)(\sa)},
\]
which, together with \eqref{ub-incursion-length} gives the desired result.

\end{proof}

\ms\subsection{The Spectral Radius}\label{spectralstuff} \

\ms \begin{proposition}[Convergence of the Spectral Radius] \label{propspectralradius}
The normalization $\bphi=P(\psi|\bX)=0$ implies that 
$1<  \rho_{\bpp} \leq 1+e^{\beta\frac{\phi_g}{2}}$
for $\beta$ large enough.
\textup{(}Recall that $\phi_g<0$ and it is defined in \eqref{phig}.\textup{)}
\end{proposition}

\ms \begin{proof} Since $\MM_{\bpp}$ is irreducible, Perron--Frobenius Theorem ensures that
\[
\rho_{\bpp}=\limsup_{p\to\infty}\sqrt[p]{\trace(\MM_{\bpp}^{p})}=\limsup_{p\to\infty}\sqrt[p]{\MM_{\bpp}^{p}(a,a)},
\]
for all $a\in\AA$.
Taking $a\in \bAA_1$ and using \eqref{periodic-perron-frobenius-J} (which we
derived from Corollary~\ref{periodic-perron-frobenius}), we deduce that
\begin{eqnarray*}
\rho_{\bpp}  = \limsup_{p\to\infty}\sqrt[p]{\MM_{\bpp}^{p}(a,a)}
               &>&\lim_{k\to\infty}\left(\MM_{\psi,1}^{kp_1}(a,a)\right)^{\frac{1}{kp_1}}\\
               &\geq&\lim_{k\to\infty}\left(\v_1(a)\w_1(a)\right)^{\frac{1}{kp_1}} 
               e^{-\frac{C_{\psi}\tau_{\psi}^k}{kp_1}} =1.
\end{eqnarray*}

\ms For the upper bound we will follow the same technique as in Subsection~\ref{proof-incursion-length}
above. Let $p=k \barp$ for some $k\in \nn$. To each $\a\in\per_p(X)$ we associate a circuit 
$C(\a)$ in $\GG_{X}$ which we decompose into its incursion--excursion path segments,
$C(\a):=a_1\lds a'_1\cdot a_2\lds a'_é\lds\cdots\lds a_\kappa\lds a'_\kappa\lds a_1$, 
as previously. We extend $C(\c)$ to 
\[
\begin{array}{ccccccc}
C_{\rm ext}(\a)
  :=a_1\lds&a'_1 &\lds  a_2\lds & a'_2&\lds \cdots \lds a_\kappa\lds&a'_\kappa&\lds a_1 \\
           &\cira&             &\cira&                             & \cira &         \\
           & a_1 &             & a_2  &                            & a_\kappa&             
\end{array} 
\]
The complementary circuit 
$C'(\a):=a_1\to a'_1\lds a_2\lds a'_2\lds \cdots\lds a_i\to a'_i\lds
\cdots\lds a_\kappa\to a'_\kappa\lds a_1$,
does not include any circuit maximizing $\phi$. As shown in Subsection~\ref{proof-incursion-length}, 
we have the upper bounds
\[
\phi(C(\a))\leq |C'(\a)| \phi_g,\,\,\, \psi(C(\a))\leq\sum_{i=1}^\kappa
\psi\left({a_i\cirb a'_i}\right)+(2\#\AA+1) |C'(\a)| \|\psi\|_{\infty}.
\]

\ms We first group the periodic points in $\per_{p}(X)$ according to the number of its incursion--excursion 
path segments, then we refine the groups so obtained by considering the end sites 
$n_1<n'_1<\cdots<n_\kappa<n'_\kappa$ of the incursion segments. By doing so we obtain
\[
\sum_{\sa\in\sper_p(X)} e^{S_{p}(\beta\phi+\psi)(\sa)}\leq
                           \sum_{\kappa=0}^{\lfloor p/2\rfloor}\sum_{\sa\in\sper(X)\atop\kappa(\sa)=\kappa} 
                           e^{(\beta\phi_g+(2\#\AA+1)\|\psi\|_{\infty})|C'(\sa)|}
                           \prod_{i=1}^\kappa e^{\psi\left({a_i\cirb a'_i}\right)}.
\]
Now, taking into account $|C'(\a)|\geq 2\kappa$ and that all the matrices $\MM_{\psi,j}$ 
associated to heavy components have the same spectral radius $e^{P(\bX|\psi)}$, we obtain
\[
\sum_{\sa\in\sper_p(X)} e^{S_{p}(\beta\phi+\psi)(\sa)}
\]
\[
\leq                           \sum_{\kappa=0}^{\lfloor p/2\rfloor}\sum_{L=2\kappa}^p 
                           e^{(\beta\phi_g+(2\#\AA+1)\|\psi\|_{\infty}+\log(\#\AA))L}
                           \sum_{n_1<n'_1<\cdots<n_\kappa<n'_\kappa}
                          \prod_{i=1}^\kappa K_{\psi}e^{(n_i'-n_i)P(\bX|\psi)},
\]
where $K_{\psi}\geq \#\AA$ is taken so large to include the counting of transitive components in $\bX$ and 
to compensate the differences in trace among transitive components and the difference between lengths 
$|a_i\cirb a'_i|$ and $n'_i-n_i$ for each $1\leq i\leq \kappa$. Since $\phi_g<0$ and
$P(\psi|\bX)=0$, the combinatoics in the distribution of the heavy components gives us 
the upper bound
\begin{eqnarray*}
\sum_{\sa\in\sper_p(X)} e^{S_{p}(\beta\phi+\psi)(\sa)}) &\leq&
                          \sum_{\kappa=0}^{\lfloor p/2\rfloor}(p-2\kappa) \binom{p}{2\kappa} 
                           e^{(\beta\phi_g+(\#\AA+1)\|\psi\|_{\infty}+\log(\#\AA))2\kappa}  K_{\psi}^{\kappa}\\
                                                       &\leq& p
                         \sum_{\kappa=0}^{p}\binom{p}{\kappa} 
                         e^{(\beta\phi_g+(2\#\AA+1)\|\psi\|_{\infty}+\log(\#\AA)+\log(K_{\psi})/2)\kappa},                                                      
\end{eqnarray*} 
for all $\beta\geq (2\#\AA+1)\|\psi\|_{\infty}+\log(\#\AA)/|\phi_g|$. 
Finally, by taking $\beta> (  \log(K_{\psi})+2(2\#\AA+1)\|\psi\|_{\infty}+\log(\#\AA))/|\phi_g|$,  
we obtain
\[
\sum_{\sa\in\sper_p(X)} e^{S_{p}(\beta\phi+\psi)(\sa)}\leq p\, \left(1+e^{\beta\phi_g/2}\right)^{p},
\]
therefore
\[ 
\rho_{\bpp}=\limsup_{p\to\infty}\sqrt[p]{\trace\left(\MM_{\bpp}^{p}\right)}
           =\limsup_{p\to\infty}\sqrt[p]{\sum_{\sa\in\sper_p(X)} e^{S_{p}(\beta\phi+\psi)(\sa)}}
           \leq 1+e^{\beta\phi_g/2},
\]
and the result follows.
\end{proof}

\ms\subsection{Excursion Series} \label{excursions} \

\ms \begin{proposition}[Convergence of the Excursion Series] \label{ex}
There exists a constant $D=D(\psi)>0$ such that for each $a\not\in\bAA:=\bigsqcup_{J=1}^{N_\phi}\bAA_J$ 
we have
\[
\sum_{j=0}^{\infty}\frac{\tMM^j_{\bpp}(a,a)}{\rho_{\bpp}^j}=
\sum_{j=0}^\infty\bMM_{\psi}^j(a,a)\pm D\thinspace e^{\beta\phi_g/3},
\]
for $\beta$ large enough.
\end{proposition}

\ms \begin{proof}
For the lower bound notice that
\begin{eqnarray*}
\sum_{j=0}^{\infty}\frac{\tMM^j_{\bpp}(a,a)}{\rho_{\bpp}^j}
    &\geq&\sum_{j=0}^{\infty}\frac{\tMM^j_{\bpp}(a,a)}{(1+e^{\beta\phi_g/2})^j}
     \geq \sum_{j=0}^{\infty}\left(1-e^{\beta\phi_g/2}\right)^j\bMM^j_{\psi}(a,a) \\
    & \geq &\sum_{j=0}^{\infty}\left(1-je^{\beta\phi_g/2}\right)\bMM^j_{\psi}(a,a) \\
    &\geq&\sum_{j=0}^{\infty}\bMM^j_{\psi}(a,a)-e^{\beta\phi_g/2}\sum_{j=0}^{\infty}j\bMM^j_{\psi}(a,a)\\
   & \geq &\sum_{j=0}^{\infty}\bMM^j_{\psi}(a,a)-e^{\beta\phi_g/2}\sum_{j=0}^{\infty}j\,
          \trace\left(\bMM^j_{\psi}\right).                                              
\end{eqnarray*}
The convergence of the series $\sum_{j=0}^{\infty}j\bMM^j_{\psi}(a,a)$ is ensured by the fact that 
$a\not\in\bAA$. The upper bound $\tMM^j_{\bpp}(a,a)\geq \bMM^j_{\psi}(a,a)$ is obtained by restricting the 
sum 
\[
\tMM^j_{\bpp}(a,a):=\sum_{\sb\in\sper_j(X)\cap\tAA^j\cap[a]}e^{S_j(\bpp)(\sb)}
\]
to periodic points maximizing $S_j\phi$.

\ms Now, since $\rho_{\bpp}> 1$ (by Proposition \ref{propspectralradius}) we obviously have 
\[
\sum_{j=0}^{\infty}\frac{\tMM^j_{\bpp}(a,a)}{\rho_{\bpp}^j}\leq
\sum_{j=0}^{\infty}\tMM^j_{\bpp}(a,a)=\sum_{j=0}^\infty \sum_{\sb\in\PP_j}e^{S_j(\bpp)(\sb)},
\]
where $\PP_j:=\{\b\in \per_j(X):\ (\b_i,\b_{i+1})\notin\bEE_{\phi}\,\,\forall\,i\in\zz\}$. 
To upper bound the terms $\sum_{\sb\in\PP_j}e^{S_j(\bpp)(\sb)}$ we use an incursion--excursion 
decomposition of the circuits associated to periodic points, similar to those employed in the proof of 
the two previous results.
For this, let $\GG'_{\bX}:=\GG_{\bX}\setminus \bigsqcup_{J=1}^{N_\phi}\GG_J$ be the 
digraph associated to the subshift $\bX'\subset \bX$ obtained by excluding all the heavy components, 
and let $\tGG_{X}$ be the subgraph spanned by the arrows in the 
complement of $\GG_{\bX}$.
Then, each $\b\in\PP_j$ defines a circuit $C(\b)$ in $\GG'_{\bX}+\tGG_{X}$
which we decompose into its incursion--excursion path segments,
\[
C(\b):=a_1\lds a'_1\lds a_2\lds a'_2\lds\cdots\lds a_\kappa\lds a'_\kappa\lds a_1,
\] 
where for each $1\leq i\leq \kappa$, $a_i\lds a'_i$ is a path in $\GG'_{\bX}$ and $a'_i\lds a_{i+1}$
is a path in $\tGG_X$. We include the two extreme cases $C(\b)$ in $\GG'_{\bX}$ or in $\tGG_X$,
by taking $\kappa=0$ and specifying which of these two possibilities holds. As previously, we extend $C(\b)$ to 
$[C_{\rm ext}(\b):=\sum_{i=1}^\kappa a_i\cirb a'_i + C'(\b)$, by adding convenient circuits of minimal length 
in $\GG'_{\bX}$. The complementary circuit 
$C'(\b):=a_1\to a'_1\lds a_2\lds a'_2\lds\cdots\lds a_\kappa\to a'_\kappa\cdot \lds a_1$,
does not include any circuit maximizing $\phi$. As in Subsection~\ref{proof-incursion-length}, 
we have the upper bounds: 
\[
\phi(C(\b))\leq |C'(\b)|\thinspace \phi_g,\,\,\, \psi(C(\b))\leq\sum_{i=1}^\kappa
\psi\left({a_i\cirb a'_i}\right)+(2\#\AA+1)\|\psi\|_{\infty}\thinspace |C'(\b)|.
\]

\ms Let $\tilde{P}=P(\psi|\bX')<0$ be the topological pressure of $\psi$ restricted to $\bX'\subset\bX$, 
the collection of all the non--heavy transitive components of $\bX$ (which is of course a union of subshifts 
of finite type). Let us suppose that $[a]\cap\bX\neq\emptyset$, then we have
\begin{eqnarray*}
\lefteqn{\sum_{\sb\in\PP_j}e^{S_j(\bpp)(\sb)}}\\
&& \leq\sum_{\sb\in\PP_j\atop\kappa(\sb)=0}e^{S_j(\bpp)(\sb)}
+\sum_{\kappa=1}^{\lfloor j/2 \rfloor}
\sum_{\sb\in\PP_j\atop\kappa(\sb)=\kappa} e^{(\beta\phi_g+(2\#\AA+1)\|\psi\|_{\infty})|C'(\sb)|}
\prod_{i=1}^\kappa e^{\psi\left({a_i\cirb a'_i}\right)}\\
&&                                             
\leq\bMM_{\psi}^j(a,a)+\sum_{\kappa=1}^{\lfloor j/2 \rfloor} \binom{j}{2\kappa}
e^{(\beta\phi_g+(2\#\AA+1)\|\psi\|_{\infty}+\log(\#\AA))(j-m)} K_{\psi}^{\kappa}e^{\tilde{P}m}\\
&&\leq \bMM_{\psi}^j(a,a)+\sum_{\kappa=1}^{\lfloor j/2 \rfloor} \binom{j}{2\kappa}
e^{(\beta\phi_g+(\#\AA+1)\|\psi\|_{\infty}+\log(\#A)+\log(K_{\psi}))(j-m)}e^{\tilde{P}m},
\end{eqnarray*}
for each $\beta\geq (2\#\AA+1)\|\psi\|_{\infty}+\log(\#\AA)/|\phi_g|$ and $j\geq 1$. Here 
$m=\sum_{i=1}^\kappa|a_i\lds a'_i|\leq j-\kappa$ and $K_\psi$ is a constant which includes the 
count and compensates the differences among transitive components of $\bX'$. Now, since $j-m\geq \kappa$, 
we have
\begin{eqnarray*}
\lefteqn{\sum_{\sb\in\PP_j}e^{S_j(\bpp)(\sb)}}\\
&\leq&
\bMM_{\psi}^j(a,a)+\sum_{\kappa=1}^{\lfloor j/2\rfloor}\binom{j}{2\kappa}
             e^{(\beta\phi_g+(\#\AA+1)\|\psi\|_{\infty}+\log(\#A)+\log(K_{\psi}))\kappa}e^{\tilde{P}(j-\kappa)}\\
                                        &\leq&
\bMM_{\psi}^j(a,a)+\sum_{\kappa=1}^{j}\binom{j}{\kappa}
\left(e^{\beta \phi_g/3}\right)^{(2\kappa)}\left(e^{\tilde{P}/2}\right)^{(j-2\kappa)}\\
                                        &\leq&
\bMM_{\psi}^j(a,a)+\left(e^{\beta \phi_g/3}+e^{\tilde{P}/2}\right)^{k}-e^{j\tilde{P}/2}
\end{eqnarray*}
for all $\beta$ greater than $(|\tilde{P}|+ (2\#\AA+1)\|\psi\|_{\infty}+\log(\#\AA)+\log(D_{\psi}))/|\phi_g|$ and all $j\geq 1$. 
With this we finally obtain
\begin{eqnarray*}
\sum_{j=0}^{\infty}\frac{\tMM^j_{\bpp}(a,a)}{\rho_{\bpp}^j}&\leq &
\sum_{j=0}^{\infty}\bMM_{\psi}^j(a,a)+
\frac{e^{\beta \phi_g/3}+e^{\tilde{P}/2}}{1-e^{\beta \phi_g/3}-e^{\tilde{P}/2}}
                                                     -\frac{e^{\tilde{P}/2}}{1-e^{\tilde{P}/2}}\\
                                                           &\leq&
\sum_{j=0}^{\infty}\bMM_{\psi}^j(a,a)+\frac{e^{\beta \phi_g/3}}{\left(1-e^{\beta \phi_g/3}
                                        -e^{\tilde{P}/2}\right)\left(1-e^{\tilde{P}/2}\right)}\\
                                                           &\leq &
\sum_{j=0}^{\infty}\bMM_{\psi}^j(a,a)+\frac{2e^{\beta \phi_g/3}}{\left(1-e^{\tilde{P}/2}\right)^2},
\end{eqnarray*}
for $\beta$ large enough.
A similar computation, for the case $[a]\cap\bX=\emptyset$ leads to 
\[
\sum_{j=0}^{\infty}\frac{\tMM^j_{\bpp}(a,a)}{\rho_{\bpp}^j}\geq
\frac{2e^{\beta\phi_g/3}}{\left(1-e^{\tilde{P}/2}\right)^2}
\]
for the same values of $\beta$. Since in this last case $\bMM_{\psi}^j(a,a)=0$ for all $j\in\nn$, 
the result follows by taking
\[
D=2\max\left(\frac{1}{\left(1-e^{\tilde{P}/2}\right)^2}\ , \
\sum_{j=0}^{\infty}j\,\trace\left(\bMM^j_{\psi}\right) \right).
\]
\end{proof}

\ms \section{Projective Stability of the Eigensystems}\

\ms \begin{proposition}[Projective Stability of the Eigensystem]~\label{proposition-stability}
Let $E$ be a finite set (with at least two elements) and let $M,N:E\times E\to\rr^+$ be irreducible matrices such that $M=e^{\pm \eta}N$, for some $\eta > 0$. 
Let $\rho_M$ the maximal eigenvalue of $M$, and $\w_M$, $\v_M$ the associated left and right positive 
eigenvectors, normalized such that $\w_M^{\dag}\v_M=1$. Let $\rho_N,\w_N$ and $\v_N$ the corresponding
quantities for $N$.  Then we have
\[\rho_M=e^{\pm\eta}\rho_N,\, \w_M=e^{\pm 2(\#E-1)\eta}\w_N\, \text{ and }
\v_M=e^{\pm 2(\#E-1)\eta}\v_N.
\]
\end{proposition}

\ms \begin{proof} First notice that the matrices $N$ and $M$ necessarily have the same period. Using 
Corollary~\ref{periodic-perron-frobenius} we readily obtain
\[
\rho_M=\limsup_{p\to\infty}\sqrt[p]{\trace(M^{p})}=\limsup_{p\to\infty}\sqrt[p]{\trace(e^{\pm\,p\eta}N^{p})}
=e^{\pm\eta}\limsup_{p\to\infty}\sqrt[p]{\trace(N^p)}=e^{\pm\eta}\rho_N.
\]

\ms For the left and right eigenvectors, let $E=\{e_1,e_2,\ldots,e_{\#E}\}$ and consider the reduced
matrices $M',N':\{e_2,\ldots,e_{\#E}\}\times\{e_2,\ldots,e_{\#E}\}\to\rr^+$ such that $M'(e,e')=M(e,e')$
for all $e,e'\in\{e_2,\ldots,e_{\#E}\}$, and similarly for $N'$. 
The right eigenvector $\v_M$ associated to $\rho_M$, normalized such that $\v_M(e_1)=1$, corresponds to 
the unique solution to the system $(M'-\rho_M\1)\x_M=\y_M$, where $\1$ is the $(\#E-1)$--dimensional identity
matrix and where $\y_M:\{e_2,\ldots,e_{\#E}\}\to\rr^+$ is such that $\y_M(e_k)=-M(e_k,e_1)$. This solution can be
obtained by using the Cramer's method, so that
\[
\x_M(e_k)=\frac{{\rm det}(M_k)}{{\rm det}(M'-\rho_M\1)}\, \, \text{ for each } 2\leq k\leq \#E.
\] 
Here $M_k$ is obtained from $M'-\rho_M\1$ by replacing its $(k-1)$--th column by the vector $\y_M$. 
The same procedure can be employed to obtain the right eigenvector $\v_N$ associated to $\rho_N$, 
and normalized such that $\v_N(e_1)=1$, by solving the equation $(N'-\rho_N\1)\x_N=\y_N$ by the Cramer's
method. Now, since $M'-\rho_M\1=e^{\pm\eta}(N'-\rho_N\1)$, $M_k=e^{\pm\eta}N_k$, and
the determinant is a $(\#E-1)$--homogeneous function, then 
\[
\x_M(e_k)=\frac{{\rm det}(M_k)}{{\rm det}(M'-\rho_M\1)}=e^{\pm\,2(\#E-1)\eta}
\frac{{\rm det}(M_k)}{{\rm det}(M'-\rho_M\1)}=e^{\pm\,2(\#E-1)\eta}\x_N(e_k),
\]
for each $2\leq k\leq\#E$, which implies that $\v_M=e^{\pm\,2(\#E-1)\eta}\v_N$ as long as they
are normalized such that $\v_M(e_1)=\v_N(e_1)$.
The argument goes the same for the left eigenvectors 
$\w_M$ and $\w_N$ associated to $\rho_M$ and $\rho_N$ respectively. 
The proposition is proved.
\end{proof}

\ms \section{Concentration of the Measure on the Heavy Components}\label{concentration}\

\ms \begin{Lemma}[Concentration on the Heavy Components]~\label{concentration-Lemma}
Let $I_K:=\{\a\in X:\ [\a_0]\cap \bX_K\neq \emptyset\}$.
Then for $\beta$ large enough we have
$\mu_{\bpp}(\cup_{K=1}^{N_\phi} I_K)\geq 1-e^{\phi_g/4}$ for $\beta$ large enough.
\end{Lemma}

\ms \begin{proof} According to Proposition~\ref{periodic-proposition}, for $\beta$ large enough and for
$\gamma = s_{\phi}-\phi_g/4$ one has
\begin{eqnarray*}
\mu_{\bpp}[a]&=&\PP^{(k\barp)}_{\bpp}[a]\exp\left(\pm e^{-\beta(\gamma-s_{\phi})}\right)\\
                   &=&\frac{\sum_{\sa\in \sper_{k\barp}(X)\cap[a]}e^{S_{k\barp}(\bpp)(\sa)}
                          }{\sum_{\sa\in \sper_{k\barp}(X)}e^{S_{k\barp}(\bpp)(\sa)}}
                                           \exp\left(\pm e^{-\beta(\gamma-s_{\phi})}\right)
\end{eqnarray*}
for each $a\in \AA$, and every $k\barp > e^{\beta\gamma}+1$.
Following the arguments developed in Subsection~\ref{proof-incursion-length}, we will find bounds for the numerator 
$\sum_{\sa\in \sper_{k\barp}(X)\cap[a]}e^{S_{k\barp}(\bpp)(\sa)}$ when $a\notin\bigsqcup_{J=1}^{N_\phi}\bAA_J$.

\ms Let $p=k\barp$, $q=\lfloor e^{\eta\,\beta}\rfloor$, with $-\phi_g/4 < \eta < -\phi_g/3$, and for 
each $b,b'\in \AA$ let $p_1=2q^2+|b\to b'|$, $p_2=p-2q^2+|b'\to b|$ and $p_0=|b\cira b'|$, where $b\cira b'$ 
is a circuit connecting $b$ and $b'$, formed by the concatenation of the minimal length paths $b\to b'$ and 
$b'\to b$. By the same argument as in Subsection~\ref{proof-incursion-length}, it follows that
\begin{equation}\label{this-is-required-before-first}
\sum_{\sa\in\sper_{p}(X)\cap[a]}e^{S_{p}(\bpp)(\sa)}=
\sum_{b,b'\in \AA}e^{-S_{p_0}(\beta\phi+\psi)\left({b\cira b'}\right)} \; \times
\end{equation}
$$  
\left(\sum_{\sa\in \sper_{p_2}(X)\cap\left[{b'\to b}\right]}e^{S_{p_2}(\beta\phi+\psi)(\sa)}\right)\; \times
\left(\sum_{\sa\in \sper_{p_1}(X)\cap[a]\cap\sigma^{-q^2}\left[{b\to b'}\right]} e^{S_{p_1}(\beta\phi+\psi)(\sa)}\right).                                                       
$$

\medskip

\ms We focus on the {\it interior sums} 
$\sum_{\sa\in \sper_{p_1}(X)\cap [a]\cap\sigma^{-q^2}\left[w(a\to a')\right]}e^{S_{p_1}(\beta\phi+\psi)(\sa)}$.

\ms Each periodic point $\a\in \per_{p_1}(X)\cap \sigma^{-q^2}\left[w(a\to a')\right]$ defines a circuit 
$C(\a)$ in $\GG_{X}$. As previously, we decompose this circuit into its {\it incursion--excursion path segments},
\[
C(\a):=a_1\lds a_1\cdot b_1\lds b'_1\cdot \ldots \cdot a_\kappa\lds a'_\kappa\cdot b_\kappa\lds b'_\kappa,
\] 
and we extend it by adding, for each $1\leq i\leq\kappa$, a circuit $a_i\cira a'_i$ in $\bX$. 
The extended circuit 
\[
C_{\rm ext}(\a)=\sum_{i=1}^\kappa a_i\cirb a'_i + 
a_1\to a_1\cdot b_1\lds b'_1\cdot \ldots \cdot a_\kappa\to a'_\kappa\cdot b_\kappa\lds b'_\kappa,          
\]
is such that $\phi(C(\a))=\phi(C_{\rm ext}(\a))$ and $\psi(C(\a))\leq 
\psi(C_{\rm ext}(\a))+2\kappa\#\AA\|\psi\|_{\infty}$. 

\ms The complementary circuit 
$C'(\a):=a_1\to a_1\cdot b_1\lds b'_1\cdot \ldots \cdot a_\kappa\to a'_\kappa\cdot b_\kappa\lds b'_\kappa$,
does not include any circuit maximizing $\phi$, therefore $\phi(C'(\a))\leq |C'(\a)|\ \phi_g$, and once 
again, we have the upper bounds
\begin{eqnarray*}
\phi(C(\a))&\leq&|C'(\a)|\times \phi_g \\  
\psi(C(\a))&\leq&\sum_{i=1}^\kappa\psi\left({a_i\cirb a'_i}\right) +(2\#\AA+1)\|\psi\|_{\infty}\ |C'(\a)|.
\end{eqnarray*}

\ms We group the periodic points in $\PP_1:=\per_{p_1}(X)\cap\sigma^{-q^2}\left[{b\to b'}\right]$ 
according to the number of its incursion--excursion path segments.
By doing so we have 
\begin{eqnarray*}
\sum_{\sa\in\PP_1} e^{S_{p_1}(\beta\phi+\psi)(\sa)}
                   &\leq & 
\sum_{\kappa=0}^{\lfloor p_1/2\rfloor}\sum_{\sa\in\PP_1\atop\kappa(\sa)=\kappa}e^{S_{p_1}(\beta\phi+\psi)(\sa)}\\
                   &\leq &
\sum_{\kappa=0}^{\lfloor p_1/2\rfloor}\sum_{\sa\in\PP_1\atop\kappa(\sc)=\kappa}
                                      e^{(\beta\phi_g+(2\#\AA+1)\|\psi\|_{\infty})|C'(\sa)|}
                                        \prod_{i=1}^\kappa e^{\psi\left({a_i\cirb a'_i}\right)}.
\end{eqnarray*}  
Now we group all the periodic points in $\a\in\PP_1$ with $\kappa(\a)=\kappa$ in classes defined by the
total length of the complementary circuit, $L(\a):=|C'(\a)|$, the relative position of the incursion 
and excursion segments, and the location of the origin inside the first incursion $a_1\lds a'_1$. 
We distinguish periodic orbits for which $L(\a)\leq L_\phi:=\lceil 6\#\AA\|\phi\|_{\infty}/|\phi_g|\rceil+1$. 
Taking into account Theorem~\ref{periodic-perron-frobenius}, the fact that $P(\psi|\bX)=0$, and the fact that
$p_1\leq 2\lfloor e^{2\beta\,\eta}\rfloor + \#\AA$, with $-\phi_g/4 < \eta < -\phi_g/3$, we obtain
$$
\sum_{\sa\in\PP_1} e^{S_{p_1}(\beta\phi+\psi)(\sa)}
$$
\begin{eqnarray}
\nonumber
&& \leq \sum_{\sa\in\PP_1\atop L(\sa)< L_\phi} e^{S_{p_1}(\beta\phi+\psi)(\sa)}+ 
\sum_{\kappa=0}^{\lfloor L/2 \rfloor}\sum_{L=L_\phi}^{p_1}
e^{(\beta\phi_g+(2\#\AA+1)\|\psi\|_{\infty}+\log(\#\AA))L}\, 
p_1 \binom{p_1}{2\kappa-1} D_{\psi}^{\kappa} \\
\nonumber
&& \leq \sum_{\sa\in\PP_1\atop L(\sa)<L_\phi} e^{S_{p_1}(\beta\phi+\psi)(\sa)}+ 
\sum_{L=L_\phi}^{p_1}e^{(\beta\phi_g+(2\#\AA+1)\|\psi\|_{\infty}+\log(\#\AA)+\log(p_1)+\log(D_\psi)+e^{-1})L}\\
\nonumber
&& \leq\sum_{\sa\in\PP_1\atop L(\sc)<L_\phi} e^{S_{p_1}(\beta\phi+\psi)(\sa)}+
2e^{\beta\frac{\phi_g}{3}L_\phi} \\
\nonumber
&& \label{this-is-required-first}
\leq\sum_{\sa\in\PP_1\atop L(\sa)<L_\phi} e^{S_{p_1}(\beta\phi+\psi)(\sa)}+
2e^{-\eta\,\beta}\sum_{\sa\in\PP_1\atop \left[\sa_{p-q^2}^{q^2-1}\right]\cap\bX_1\neq\emptyset} 
e^{S_{p_1}(\beta\phi+\psi)(\sa)},
\end{eqnarray}
for all $\beta$ larger than a convenient $\beta(L,\eta)$.
The constant $D_{\psi}$ is taken large enough to include the counting of transitive components
in $\bX$, and to compensate the differences in trace among transitive components, and the difference 
between lengths $|a_i\cirb a'_i|$ and $m_i$ for each $1\leq i\leq \kappa$. The factor $p_1$ takes into account 
the all the possible locations of the origin with inside the first incursion. We are also considering the fact 
that $2\kappa(\a)\leq L(\a)$. 

\ms Let us now consider the periodic points 
$\PP_2:=\{\a\in \per_{p_1}(X)\cap\sigma^{-q^2}\left[w(a\to a')\right]:\ L(\a) <  L_\phi\}$.
Using exactly the same argument as in Subsection~\ref{proof-incursion-length}, we deduce that  
\begin{equation}\label{this-is-required}
\sum_{\sa\in\PP_2} e^{S_{p_1}(\beta\phi+\psi)(\sa)}\leq
                     \left(1+e^{\beta\frac{\phi_g}{3}}\right)
\sum_{\sa\in\PP_2\atop m_i > p_1/q\Rightarrow a_i\lds a'_i \text{ in } \bigsqcup_{K\leq N_\phi}\GG_K} 
                              e^{S_{p_1}(\beta\phi+\psi)(\sa)} 
\end{equation}
for $\beta$ large enough.

\ms We organize the periodic points in 
$\PP_3:=\{\a\in\PP_2:\,|a_i\lds a'_i|\geq p_1/q\,\Rightarrow\,a_i\lds a'_i\text{ in }\bigsqcup_{K\leq N_\phi}\GG_K\}$,
by classes defined by fixing the excursion segments $b_i\lds b'_i$, the incursion into non--heavy components
$a_i\lds a'_i$, and the input--output vertices and lengths of the incursions into heavy components.  
By definition, two points $\a,\a'$ in the same class are such that $S_{p_1}\phi(\a)=S_{p_1}\phi(\a')$. 
Now, for each $\a\in\PP_3$ we have 
\[
\sum_{a_i\lds a'_i \text{ in } \cup_{K}\GG_K}|a_i\lds a'_i|\geq (2q^2+|b\to b'|)\left(1-\frac{2L_\phi}{q}\right).
\]
Let $p*={\rm lcm}(p_K:\ 1\leq K\leq N)$ and consider a refinement of a particular class, in 
subclasses defined by location of the origin. Now, for each subclass such that the segment containing the
origin is not included in $\bigsqcup_{K=1}^{N_\phi}\GG_K$, there are at least 
\[
\frac{1}{2qp*}\left(\frac{(2q^2+|b\to b'|)(1-2L_\phi/q)}{L_\phi}-2q\right)\geq \frac{q}{2p^*L_\phi}
\]
subclasses such that the segment containing the origin is included in $\cup_{K=1}^{N_\phi}\GG_K$, as long as 
$q$ is large enough.
The argument is exactly the same as the one developed in 
Subsection~\ref{proof-incursion-length} of Appendix~\ref{auxiliary-inequalities}. 
From~\eqref{this-is-required-first}, \eqref{this-is-required} and the previous inequality it follows that
\[
\sum_{\sa\in\PP_1\atop \sa_0\notin\cup_{K=1}^{N_\phi}\bAA_K} e^{S_{p_1}(\beta\phi+\psi)(\sa)}\leq 
2e^{-\beta\,\eta} (3p^*L_\phi+1)\sum_{\sc\in\PP_1} e^{S_{p_1}(\beta\phi+\psi)(\sa)},
\]
for $\beta$ large enough. Now, using~\eqref{this-is-required-before-first} and
Proposition~\ref{periodic-proposition} we derive the 
\[
\sum_{\sa_0\in \AA\setminus\bigsqcup_{K=1}^{N_\phi}\bAA_K}\mu_{\bpp}[\a_0]\leq 
2e^{-\beta\,\eta}(3p^*L_\phi+1)\exp\left(e^{-\beta(\gamma-s_{\phi})}\right)\leq e^{\beta\phi_g/4}
\]
for $\beta$ large enough, and the result follows.
\end{proof}


\end{document}